\documentclass[11pt]{article}

\usepackage[letterpaper,top=2cm,bottom=2.5cm,left=3cm,right=3cm,marginparwidth=1.75cm]{geometry}
\usepackage{tikz} 
\usetikzlibrary{positioning}

\usepackage{graphicx}
\usepackage{amssymb}

\usepackage{epstopdf}

\usepackage[english]{babel}
\usepackage{latexsym}
\usepackage{mathrsfs}

\usepackage{color}

\usepackage{float}
\usepackage{mathtools}

\usepackage{amsmath}
\usepackage{amsfonts}
\numberwithin{equation}{section}
\usepackage{enumerate}
\usepackage{amsthm}

\usepackage[title]{appendix}
\usepackage{algorithm}
\usepackage{algpseudocode}
\usepackage[dvipsnames]{xcolor}
\usepackage[bookmarks=true,colorlinks=true,linkcolor={blue},urlcolor={blue}, citecolor={blue},pdfstartview={XYZ null null 1.22}]{hyperref}

\providecommand{\keywords}[1]{\textbf{{Keywords.}} #1}
\providecommand{\ams}[1]{\textbf{{MSC}} #1}

\def\be#1\ee{\begin{equation}#1\end{equation}}

\theoremstyle{definition}

\newcommand{\xx}{{\bf x}}

\newcommand{\vv}{{\bf v}}
\newcommand{\BB}{{\bf B}}
\newcommand{\pp}{{\bf p}}
\newcommand{\qq}{{\bf q}}
\newcommand{\ee}{{\bf e}}

\newcommand{\EE}{{\bf E}}
\def\RR{\mathbb R}

\newcommand{\com}[1]{}

\title{Full-horizon magnetic control of the Vlasov--Poisson system for plasma confinement}

\author{G. Albi,\thanks{Department of Computer Science,
		University of Verona, Strada le Grazie 15, 37134 Verona, ITALY. (giacomo.albi@univr.it)}\and
	G. Dimarco,\thanks{Department of Mathematics and Computer Science \& Center for Modeling, Computing and Statistics (CMCS), University of Ferrara, via Machiavelli 30, 44121 Ferrara, ITALY. (giacomo.dimarco@unife.it)} \and
	F. Ferrarese, \thanks{Department of Mathematics and Computer Science \& Center for Modeling, Computing and Statistics (CMCS), University of Ferrara, via Machiavelli 30, 44121 Ferrara, ITALY. (federica.ferrarese@unife.it)} \and
	L. Pareschi \thanks{Maxwell Institute and Department of Mathematics
		School of Mathematical and Computer Sciences
		Heriot-Watt University, Edinburgh, UK (L.Pareschi@hw.ac.uk) \& Department of Mathematics and Computer Science \& Center for Modeling, Computing and Statistics (CMCS), University of Ferrara, via Machiavelli 30, 44121 Ferrara, ITALY. (lorenzo.pareschi@unife.it) 		
	}
}

\begin{document}
\maketitle
\begin{abstract}
We study the full-horizon control of a two-dimensional Vlasov--Poisson
system through an externally imposed magnetic field.  The control is
represented by time-dependent coefficients associated with smooth,
spatially localized basis functions, and is designed to steer selected
position and velocity moments of the distribution function while limiting particle dispersion and the
control effort.  Starting from the kinetic constrained optimization
problem, we derive the continuous first-order optimality conditions,
including a backward Vlasov equation, an elliptic adjoint equation, and the
variational inequality associated with bounded controls.  We then project
this system onto a particle--grid representation and identify the resulting
particle adjoints and control gradient. The full-horizon field is then computed
on a reduced particle ensemble, initialized by an instantaneous control,
and subsequently applied to a larger set of particles.  Numerical experiments
indicate that the instantaneous control field
already provides an effective baseline, while the full-horizon correction
can further reduce the relative thermal energy near the boundaries and
improve confinement over the prescribed time interval.
\end{abstract}


\ams{49M41; 35Q83; 35Q93; 49K20; 65M75; 82D10}

\section{Introduction}\label{sec:intro}
Plasma dynamics is governed by the collective motion of charged particles and by their interaction with both externally imposed and self-consistent electromagnetic fields \cite{bittencourt2013fundamentals,freidberg2008plasma,wesson2011tokamaks}. A particularly important application is magnetic confinement fusion, whose objective is to confine a high-temperature plasma for sufficiently long times to sustain thermonuclear reactions while limiting energy and particle losses \cite{freidberg2008plasma,wesson2011tokamaks,garabedian2003computational}.
Depending on the physical regime, plasmas may be thermal or non-thermal and may exhibit phenomena spanning a wide range of spatial and temporal scales. This multiscale character has led to different levels of description, including fluid, kinetic, and hybrid models \cite{Degond,Sonnen2,dimarco2014numerical,Cafl,dimarco2014asymptotic}. Fluid descriptions, such as magnetohydrodynamics, provide an effective representation of the macroscopic evolution of the plasma and play a central role in the analysis of equilibrium and large-scale stability in magnetic-fusion devices \cite{ravindran2005real,Blum}. However, several processes relevant to tokamaks and stellarators, including non-equilibrium transport, energetic-particle dynamics, velocity-space anisotropies, and losses near the plasma boundary, require a kinetic description. Kinetic models describe the evolution of particle distribution functions in phase space and can therefore represent effects that cannot be adequately captured through a finite set of macroscopic quantities \cite{crestetto2012kinetic,medaglia2025dsmc,Dima,Crous}.
Among kinetic descriptions, the Vlasov equation constitutes a fundamental model for collisionless plasmas. It describes the transport of the particle distribution under the action of electromagnetic forces, generally coupled to the distribution itself through Maxwell’s or Poisson’s equations. The high dimensionality of phase space, together with the nonlinear coupling between particles and fields, makes the numerical approximation of Vlasov-type systems particularly demanding. These challenges have motivated the development of a broad range of numerical techniques, including semi-Lagrangian, conservative, multiscale, and asymptotic-preserving methods \cite{cheng1976integration,crouseilles2010conservative,degond2017asymptotic,dimarco2015numerical,sonnendrucker1999semi,Sonnen,dimarco2014asymptotic,dimarco2014numerical,filbet2016asymptotically,filbet2017asymptotically,filbet2018numerical,zhu2017vlasov,Ross,Gamba}. In this context, particle-in-cell (PIC) methods are especially attractive for high-dimensional problems. They approximate the distribution function through a finite ensemble of numerical particles, whose trajectories are advanced in a Lagrangian fashion, while the electromagnetic fields are computed on a spatial grid \cite{chacon2016optimization,hairer2018energy,gu2022hamiltonian,deluzet2023efficient}. This combination considerably reduces the computational cost associated with a direct discretization of the full phase space while retaining the essential kinetic features of the model.

The control of plasma dynamics is a central issue in magnetic confinement fusion. In both tokamaks and stellarators, externally generated magnetic fields can be used to modify particle trajectories, improve confinement, reduce cross-field transport, and limit the flux of energetic particles toward boundaries \cite{davies2023magnetic,fasoli2016computational,garabedian2003computational,grandgirard2013gyrokinetic}. The design of effective controls is nevertheless challenging because the magnetic field acts on an evolving particle distribution and must account for the nonlinear coupling between the plasma and its self-consistent fields. Furthermore, a successful strategy should not only reduce instantaneous losses but also anticipate the subsequent evolution of the plasma and preserve confinement over relevant time scales. Both the analytical foundations and the numerical realization of kinetic plasma control have therefore attracted increasing attention in recent years \cite{caprino2012time,Knopf,knopf2020optimal,weber2021optimal,einkemmer2024suppressing,einkemmer2025control,chen2025control,lu2025dynamical,lu2025controlling,ferrarese2025control}. In particular, instantaneous magnetic control strategies for deterministic and uncertain Vlasov--Poisson models were developed in \cite{albi2025instantaneous,albi2025robust}. These approaches replace optimization over the entire time horizon with a sequence of local-in-time control problems. As a result, they are computationally inexpensive and can be readily coupled with kinetic solvers. Their local character, however, also represents an important limitation: the control is determined using only information available at the current time and does not explicitly account for the future evolution of the distribution function, the self-consistent field, or the resulting plasma position. This lack of predictive capability may become particularly relevant in magnetic confinement applications, where the plasma response to external actuation may be delayed and the primary objective is to maintain confinement and prevent particle losses over extended time intervals.

Motivated by the limitations of local in time control strategies, in this work we consider a deterministic, collisionless kinetic model and seek an optimal magnetic control over a prescribed time horizon. We describe the evolution of a population of negatively charged particles, identified with electrons, while the ions are treated as a fixed neutralizing background. In dimensionless variables, the governing Vlasov--Poisson system reads
\begin{equation}\label{eq:Vlasov_Poisson_intro}
	\begin{split}
		&\partial_t f(t,\xx,\vv) + \vv \cdot \nabla_\xx f(t,\xx,\vv) + F(t,\xx,\vv) \cdot \nabla_\vv f(t,\xx,\vv) = 0, \\
		&-\Delta_\xx \phi(t,\xx) = \rho(t,\xx)-\bar{\rho},
		\qquad E(t,\xx) = -\nabla_\xx \phi(t,\xx),
	\end{split}
\end{equation}
supplemented with the initial condition \(f(0,\xx,\vv)=f^0(\xx,\vv)\). 
   The electron density is obtained as the velocity marginal of the distribution function
\begin{equation*}
    \rho(t,\mathbf{x}) = \int_{\mathbb{R}^{d_v}} f(t,\mathbf{x},\mathbf{v})\,d\mathbf{v},
\end{equation*}
and the ion density is assumed to be homogeneous and stationary. Under this assumption, the net charge density on the right-hand side of the Poisson equation reduces to $\rho - \bar{\rho}$, where $\bar{\rho}$
is the spatial mean of the electron density, which is conserved in time by the Vlasov equation and plays the role of the neutralizing ion background.
With the sign convention adopted here, the effective Lorentz force acting on the electrons is given by
\begin{equation}\label{eq:force_Vlasov}
	F(t,\xx,\vv) = \EE(t,\xx) + \vv \times \BB_{ext}(t,\xx),
\end{equation}
where \(\BB_{ext}(t,\xx)\) denotes the externally applied magnetic field, which plays the role of the control variable, while \(\EE(t,\xx)\) is the self-consistent electric field determined through the Poisson equation in \eqref{eq:Vlasov_Poisson_intro}.

Our starting point is the instantaneous magnetic control framework introduced in \cite{albi2025instantaneous,albi2025robust}. In the deterministic setting of \cite{albi2025instantaneous}, the control was obtained through a sequence of local in time optimization problems and represented by a spatially piecewise constant, time-dependent magnetic field. This spatial localization reflects the physical constraint that arbitrarily localized magnetic configurations cannot be generated in practice. The approach was subsequently extended in \cite{albi2025robust} to systems affected by uncertainty and collisions. Although these strategies are computationally efficient, their instantaneous nature prevents them from explicitly anticipating the future evolution of the plasma.

The aim of the present work is consequently to move beyond this instantaneous formulation and determine an open loop optimal magnetic control over the entire time interval \([0,T]\). Full-horizon magnetic control of the Vlasov–Poisson system has been studied analytically in \cite{knopf2020optimal}, for fields generated by fixed coils, and numerically in \cite{Knopf}, where state and adjoint equations are solved by a particle-in-cell method. With respect to these works, we consider moment-based tracking and dispersion objectives in a 2D–2V geometry with reflecting walls, derive a particle–grid optimality system consistent with the discretize-then-optimize approach, initialize the optimization with the instantaneous control of \cite{albi2025instantaneous}, thus allowing a direct comparison between the two strategies, and compute the control on a reduced particle ensemble before validating it on a larger one.

To this end, we study the associated first-order optimality conditions through a coupled forward--backward system \cite{lu2025dynamical, Knopf, knopf2020optimal}. Accounting for the full time horizon allows the control to anticipate the subsequent evolution of the particle distribution and to balance immediate confinement with longer-term objectives. This predictive capability comes however at a substantially higher computational cost, since the state equation must be coupled with a backward adjoint equation and an optimality condition for the control. 
More precisely, we consider the constrained optimization problem
\begin{equation}\label{eq:min_problem_continuos_intro}
	\min_{\BB_{\mathrm{ext}}\in\mathcal{B}_{\mathrm{ad}}}
	\mathcal{J}(f,\BB_{\mathrm{ext}}),
	\qquad
	\text{subject to \eqref{eq:Vlasov_Poisson_intro}}.
\end{equation}
Here, \(\mathcal{B}_{\mathrm{ad}}\) denotes the set of admissible
magnetic controls. At this stage, we introduce the cost functional
in the general form
\begin{equation}\label{eq:functional_continuous_intro}
	\mathcal{J}(f,\BB_{\mathrm{ext}})
	=
	\int_0^T
	\left(
	\mathcal{D}(f)(t)
	+
	\frac{\gamma}{2}
	\mathcal{R}(f,\BB_{\mathrm{ext}})(t)
	\right)dt
	+
	\Psi(f(T)),
\end{equation}
where \(\gamma>0\). The running cost \(\mathcal{D}(f)\)
measures the deviation of selected position- and velocity-space
observables from their prescribed targets and controls their
dispersion. The regularization term
\(\mathcal{R}(f,\BB_{\mathrm{ext}})\) penalizes the control effort,
whereas \(\Psi(f(T))\) promotes the desired plasma configuration at
the final time \(T\). The precise definitions of these terms are
introduced in Section~\ref{sec:problem_setting}, after specifying
the spatial parameterization of the external magnetic field.
In the following, the first-order optimality conditions are derived
through a Lagrangian (adjoint) approach.
Solving the complete optimality system directly
at the kinetic scale would, however, be prohibitively expensive.
We therefore synthesize the control using a reduced particle
representation of the plasma dynamics. The reduced system is employed
to compute the open-loop optimal field, which is subsequently
reconstructed and applied to a larger particle simulation. In this
way, the expensive optimization is performed on a representative
particle ensemble, while the effectiveness and transferability of
the resulting control are assessed on the higher-dimensional
dynamics.

The remainder of the paper is organized as follows.
Section~\ref{sec:problem_setting} introduces the model and the optimal
control objective. Section~\ref{sec:derivation_control} derives the
continuous optimality system and its particle--grid representation.
Section~\ref{sec:numerical_optimization} describes the numerical approximation and
the optimization method. Section~\ref{sec:numexp} presents the
computational experiments. A final section summarizes the main
conclusions and discusses possible future developments. The
appendices provide the discretize-then-optimize derivation and the
instantaneous control strategy employed to initialize the full-horizon
optimization.

\section{Problem setting}\label{sec:problem_setting}

We consider a reduced \(d_x\times d_v=2\times 2\) setting, with particle
positions \(\xx=(x,y)\in\Omega_{\xx}\subset\RR^2\) and velocities
\(\vv=(v_x,v_y)\in\Omega_{\vv}\subset\RR^2\). The dynamics is studied
in the plane orthogonal to the dominant component of the externally
applied magnetic field. 
To obtain a physically and computationally realizable control, we
partition the spatial domain into \(N_c^x\times N_c^y\) cells, with
\[
N_g=N_c^xN_c^y.
\]
Denoting the cells by \(C_k\), \(k=1,\ldots,N_g\), the partition
satisfies
\[
\Omega_{\xx}
=
\bigcup_{k=1}^{N_g}C_k,
\qquad
C_k\cap C_j=\emptyset
\quad\text{for }k\neq j.
\]
We then associate with each cell a time-dependent control coefficient
\(b_k(t)\) and introduce the control vector
\[
\mathbf{b}(t)
:=
\bigl(b_1(t),\ldots,b_{N_g}(t)\bigr).
\]
For notational convenience, we set
\[
\mathcal{K}_k(\xx)
:=
\mathcal{K}(\xx,\xx_c^k),
\qquad
k=1,\ldots,N_g,
\]
where \(\xx_c^k\) denotes the center of \(C_k\). The scalar component
of the external magnetic field is reconstructed from the control
coefficients as
\begin{equation}\label{eq:Bi_cells}
	B_{\mathbf{b}}(t,\xx)
	=
	\sum_{k=1}^{N_g}
	b_k(t)\mathcal{K}_k(\xx).
\end{equation}

Each kernel \(\mathcal{K}_k\) is assumed to be nonnegative and provides
a smooth approximation of the indicator function of the corresponding
cell. More precisely, it is approximately constant inside \(C_k\)
and smoothly decays to zero outside it. Consequently, the resulting
magnetic field is nearly piecewise constant over the partition while
remaining globally smooth across cell interfaces.
The corresponding three-dimensional external magnetic field is
defined as
\begin{equation}\label{eq:B_ext_vector}
	\BB_{\mathrm{ext}}[\mathbf{b}](t,\xx)
	:=
	\bigl(0,0,B_{\mathbf{b}}(t,\xx)\bigr).
\end{equation}
The scalar component \(B_{\mathbf{b}}\) is subject to the pointwise
amplitude constraint
\[
\left|B_{\mathbf{b}}(t,\xx)\right|\leq B_{\rm max},
\]
where \(B_{\rm max}>0\) denotes the maximum admissible magnetic-field intensity.
With this convention, the planar component of the magnetic force is
given by
\[
\vv\times\BB_{\mathrm{ext}}[\mathbf{b}]
=
\bigl(v_yB_{\mathbf{b}},-v_xB_{\mathbf{b}}\bigr).
\]

The precise form of the kernel, the choice of the number \(N_g\) of
control cells, and a sensitivity analysis with respect to this number
will be discussed in the numerical experiments section.
The corresponding two-dimensional Vlasov--Poisson system reads
\begin{equation}\label{eq:Vlasov_Poisson}
	\begin{split}
		&\partial_t f(t,\xx,\vv)
		+\vv\cdot\nabla_{\xx}f(t,\xx,\vv)
		+\left(
		\EE(t,\xx)
		+\vv\times\BB_{\mathrm{ext}}[\mathbf{b}](t,\xx)
		\right)\cdot\nabla_{\vv}f(t,\xx,\vv)
		=0,
		\\
		&-\Delta_{\xx}\phi(t,\xx)=\rho(t,\xx)-\bar{\rho},
		\qquad
		\EE(t,\xx)=-\nabla_{\xx}\phi(t,\xx),
		\qquad
		f(0,\xx,\vv)=f_0(\xx,\vv).
	\end{split}
\end{equation}
Here, \(f_0(\xx,\vv)\) is the prescribed initial distribution,
whereas \(f(t,\xx,\vv)\) denotes the electron distribution function
in phase space.   The quantity 
\begin{equation*}
    \bar{\rho} = \frac{1}{|\Omega_{\xx}|}\int_{\Omega_{\xx}} \rho(t,\mathbf{x})\,d\mathbf{x}
\end{equation*} 
is conserved in time by the Vlasov equation and plays the role of the neutralizing ion background. In both equations, periodic boundary
conditions are imposed in the \(x\)-direction, while reflective
boundary conditions are prescribed in the \(y\)-direction to model
the presence of material walls. The control objective is to steer
the plasma by means of the external magnetic field and, in particular,
to reduce particle transport toward and interaction with these walls.

We now reformulate the continuous optimal control problem
\eqref{eq:min_problem_continuos_intro} in terms of the finite set of
control coefficients. Let
\[
\mathcal{I}:=\{\mathrm{x},\mathrm{v}\},
\qquad
\Omega:=\Omega_{\xx}\times\Omega_{\vv}.
\]
The set of admissible controls is defined as
\begin{equation}\label{eq:admissible_controls}
	\mathcal{U}_{\mathrm{ad}}
	:=
	\left\{
	\mathbf{b}\in L^2\bigl(0,T;\RR^{N_g}\bigr):
	\left|B_{\mathbf{b}}(t,\xx)\right|\leq B_{\rm max}
	\ \text{for a.e. }(t,\xx)
	\in(0,T)\times\Omega_{\xx}
	\right\}.
\end{equation}
The reduced optimal control problem is therefore
\begin{equation}\label{eq:min_prob_2d}
	\min_{\mathbf{b}\in\mathcal{U}_{\mathrm{ad}}}
	\mathcal{J}(f,\mathbf{b}),
	\qquad
	\text{subject to \eqref{eq:Vlasov_Poisson}}.
\end{equation}
We first specify the regularization of the control coefficients.
The regularization of the reconstructed magnetic field would
naturally lead to
\begin{equation}\label{eq:exact_control_penalty}
	\begin{split}
		\mathcal{R}_{\mathrm{field}}(f,\mathbf{b})(t)
		:=\;&
		\int_{\Omega}
		\left|B_{\mathbf{b}}(t,\xx)\right|^2
		f(t,\xx,\vv)\,d\xx\,d\vv
		\\
		=\;&
		\sum_{k=1}^{N_g}
		\left|b_k(t)\right|^2
		\int_{\Omega}
		\mathcal{K}_k(\xx)^2
		f(t,\xx,\vv)\,d\xx\,d\vv
		\\
		&+
		2\sum_{1\leq k<j\leq N_g}
		b_k(t)b_j(t)
		\int_{\Omega}
		\mathcal{K}_k(\xx)\mathcal{K}_j(\xx)
		f(t,\xx,\vv)\,d\xx\,d\vv.
	\end{split}
\end{equation}
Since the kernels are smooth approximations of disjoint indicator
functions, they satisfy
\begin{equation}\label{eq:assumption-kernel}
\mathcal{K}_k^2\simeq\mathcal{K}_k,
\qquad
\mathcal{K}_k\mathcal{K}_j\simeq0
\quad\text{for }k\neq j,
\end{equation}
except within the narrow transition regions near cell interfaces.
Motivated by this localization property, we neglect the overlap terms
and adopt the diagonal approximation
\begin{equation}\label{eq:diagonal_control_penalty}
	\mathcal{R}_{\mathrm{diag}}(f,\mathbf{b})(t)
	:=
	\sum_{k=1}^{N_g}
	\left|b_k(t)\right|^2
	\int_{\Omega}
	\mathcal{K}_k(\xx)
	f(t,\xx,\vv)\,d\xx\,d\vv.
\end{equation}
This approximation penalizes the control coefficients independently
and avoids coupling between controls associated with different cells.
It will be used throughout the remainder of the paper.
The cost functional is defined by
\begin{equation}\label{eq:functional_2d}
	\mathcal{J}(f,\mathbf{b})
	:=
	\int_0^T
	\mathcal{L}(f,\mathbf{b})(t)\,dt
	+
	\Psi(f(T)),
\end{equation}
where the running cost is
\begin{equation}\label{eq:running_cost_2d}
	\mathcal{L}(f,\mathbf{b})(t)
	:=
	\sum_{q\in\mathcal{I}}
	\mathcal{D}_q(f)(t)
	+
	\frac{\gamma}{2}
	\mathcal{R}_{\mathrm{diag}}(f,\mathbf{b})(t),
\end{equation}
with \(\gamma>0\). The first contribution measures the deviation
from the prescribed target configuration, whereas the second
penalizes the magnitude of the localized control coefficients.
For each \(q\in\mathcal{I}\), the tracking term is defined as
\begin{equation}\label{eq:tracking_cost_2d}
	\begin{split}
		\mathcal{D}_q(f)(t)
		:=\;&
		\frac{\alpha_q}{2}
		\int_{\Omega}
		\left|
		\psi_q(\xx,\vv)-\widehat{\psi}_q
		\right|^2
		f(t,\xx,\vv)\,d\xx\,d\vv
		\\
		&+
		\frac{\beta_q}{2}
		\int_{\Omega}
		\left|
		\psi_q(\xx,\vv)-m_q[f](t)
		\right|^2
		f(t,\xx,\vv)\,d\xx\,d\vv,
	\end{split}
\end{equation}
where \(\alpha_q,\beta_q\geq0\). The functions \(\psi_q\) represent
the selected physical- or velocity-space observables, while
\(\widehat{\psi}_q\) denotes the corresponding prescribed target.
The first term in \eqref{eq:tracking_cost_2d} promotes the
concentration of the distribution around the target value
\(\widehat{\psi}_q\). The second term penalizes the dispersion of
the observable \(\psi_q\) around its instantaneous mean.
More precisely, the mean of the observable \(\psi_q\) is defined by
\begin{equation}\label{eq:mean_observable}
	m_q[f](t)
	:=
	\frac{1}{\mathcal{M}}
	\int_{\Omega}
	\psi_q(\xx,\vv)
	f(t,\xx,\vv)\,d\xx\,d\vv,
	\qquad
	\mathcal{M}
	:=
	\int_{\Omega}
	f_0(\xx,\vv)\,d\xx\,d\vv.
\end{equation}
Finally, the terminal cost is defined as
\begin{equation}\label{eq:terminal_cost_2d}
	\Psi(f(T))
	:=
	\sum_{q\in\mathcal{I}}
	\mathcal{D}_q(f)(T).
\end{equation}
It has the same structure as the running
cost, but is evaluated at the final time \(T\). Its purpose is to
promote the desired plasma configuration at the end of the control
horizon.

\section{Derivation of the optimality system}\label{sec:derivation_control}

The aim of this section is to introduce and compare two formally
equivalent strategies for solving the minimization problem
\eqref{eq:min_prob_2d}. 
We first derive the
continuous optimality system through a Lagrangian formulation and
then project the resulting forward-backward equations onto a particle
approximation of the distribution function. This procedure yields a
particle-level optimality system involving the characteristic
variables, their corresponding adjoint variables, and the gradient
with respect to the control coefficients. We subsequently relate this construction to the
alternative procedure in which the particle-grid approximation
is introduced before deriving the optimality conditions. When the same
kernel interpolation, quadrature rules, and particle discretization
are employed, the two approaches formally lead to the same
particle-level optimality system, up to the errors introduced by the
discretization.

\subsection{Continuous first-order optimality conditions}
We derive here the continuous first-order optimality conditions associated with the constrained minimization problem \eqref{eq:min_prob_2d} through a Lagrangian approach. Let \(p(t,\xx,\vv)\) and \(\lambda(t,\xx)\) denote the adjoint variables associated with the Vlasov equation and the Poisson equation in \eqref{eq:Vlasov_Poisson}, respectively. For a given control vector \(\mathbf{b}=(b_1,\ldots,b_{N_g})\), we recall the notation
\[
F_{\mathbf{b}}(t,\xx,\vv)
:=
-\nabla_{\xx}\phi(t,\xx)
+
\vv\times
\BB_{\mathrm{ext}}[\mathbf{b}](t,\xx).
\]
To distinguish the Lagrangian functional from the running cost
\(\mathcal{L}\) introduced in \eqref{eq:running_cost_2d}, we denote
the former by \(\mathscr{L}\). It is defined as
\begin{equation}\label{eq:Lagrangian_OtD}
	\begin{split}
		\mathscr{L}(f,\phi,p,\lambda,\mathbf{b})
		:=\;&
		\int_0^T
		\mathcal{L}(f,\mathbf{b})(t)\,dt
		+
		\Psi(f(T))
		\\
		&+
		\int_0^T\int_{\Omega}
		p(t,\xx,\vv)
		\Bigl(
		\partial_t f(t,\xx,\vv)
		+
		\vv\cdot\nabla_{\xx}f(t,\xx,\vv)
		\\
		&\hspace{4.7cm}
		+
		F_{\mathbf{b}}(t,\xx,\vv)
		\cdot\nabla_{\vv}f(t,\xx,\vv)
		\Bigr)
		\,d\xx\,d\vv\,dt
		\\
		&+
		\int_0^T\int_{\Omega_{\xx}}
		\lambda(t,\xx)
		\left(
		-\Delta_{\xx}\phi(t,\xx)
		-
		\int_{\Omega_{\vv}}
		f(t,\xx,\vv)\,d\vv
		+
		\bar{\rho}
		\right)
		\,d\xx\,dt,
	\end{split}
\end{equation}
where
\(\Omega=\Omega_{\xx}\times\Omega_{\vv}\). The running cost \(\mathcal{L}\), the tracking terms \(\mathcal{D}_q\), the diagonal control regularization, and the terminal cost are defined in
\eqref{eq:running_cost_2d},\eqref{eq:tracking_cost_2d},\eqref{eq:diagonal_control_penalty}, and \eqref{eq:terminal_cost_2d}, respectively. In deriving the optimality conditions, we assume that all functions are sufficiently regular and that the boundary terms arising from integration by parts vanish under the prescribed state boundary conditions and the corresponding compatible adjoint boundary
conditions. Let $(f^\star,\phi^\star,p^\star,\lambda^\star,\mathbf{b}^\star)$ denote the optimal state-adjoint-control. Stationarity, together with the variational optimality condition for the constrained control, gives
\begin{equation}\label{eq:optimality_conditions}
	\begin{cases}
		D_p\mathscr{L}=0,\\
		D_\lambda\mathscr{L}=0,\\
		D_f\mathscr{L}=0,\\
		D_\phi\mathscr{L}=0,\\
		\displaystyle
		\left\langle
		D_{\mathbf{b}}\mathscr{L},
		\mathbf{c}-\mathbf{b}^\star
		\right\rangle
		\geq0
		\quad
		\text{for every }
		\mathbf{c}\in\mathcal{U}_{\mathrm{ad}},
	\end{cases}
\end{equation}
where all derivatives are evaluated at the optimal value. If
\(\mathbf{b}^\star\) lies in the interior of
\(\mathcal{U}_{\mathrm{ad}}\), the last condition reduces to
\(D_{\mathbf{b}}\mathscr{L}=0\). The first two conditions in \eqref{eq:optimality_conditions} recover
the forward Vlasov--Poisson system
\begin{equation}\label{eq:forward_continuos}
	\begin{cases}
		\partial_t f
		+
		\vv\cdot\nabla_{\xx}f
		+
		F_{\mathbf{b}}\cdot\nabla_{\vv}f
		=0,
		\\
		-\Delta_{\xx}\phi
		=
		\rho-\bar{\rho},
	\end{cases}
\end{equation}
where, for brevity, the dependence on \(t,\xx,\vv\) has been omitted and
the forward system is supplemented with the initial condition \(f(0,\xx,\vv)=f_0(\xx,\vv)\).
The variations with respect to \(f\) and \(\phi\) yield the backward
system
\begin{equation}\label{eq:backward_continuos}
	\begin{cases}
		-\left(
		\partial_t p
		+
		\vv\cdot\nabla_{\xx}p
		+
		F_{\mathbf{b}}\cdot\nabla_{\vv}p
		\right)
		-
		\lambda
		+
		\mathcal{G}
		=0,
		\\[2mm]
		-\Delta_{\xx}\lambda
		-
		\nabla_{\xx}\cdot
		\left(
		\displaystyle
		\int_{\Omega_{\vv}}
		f\,\nabla_{\vv}p\,d\vv
		\right)
		=0,
	\end{cases}
\end{equation}
where the source term \(\mathcal{G}\) is the derivative of the running cost
with respect to \(f\). It is given by
\begin{equation}\label{eq:G}
	\begin{split}
		\mathcal{G}(t,\xx,\vv)
		=\;&
		\sum_{q\in\mathcal{I}}
		\left[
		\frac{\alpha_q}{2}
		\left|
		\psi_q(\xx,\vv)-\widehat{\psi}_q
		\right|^2
		\right.
		\\
		&\left.
		+
		\frac{\beta_q}{2}
		\left|
		\psi_q(\xx,\vv)-m_q[f](t)
		\right|^2
		\right]
		+
		\frac{\gamma}{2}
		\sum_{k=1}^{N_g}
		|b_k(t)|^2\mathcal{K}_k(\xx).
	\end{split}
\end{equation}
The backward Vlasov equation is supplemented with a terminal
condition obtained from
\[
p(T,\xx,\vv)
+
\frac{\delta\Psi(f(T))}
{\delta f(T,\xx,\vv)}
=
0.
\] 
Therefore,
\begin{equation}\label{eq:terminal_continuos}
	p(T,\xx,\vv)
	=
	-\mathcal{G}_T(\xx,\vv),
\end{equation}
where
\begin{equation}\label{eq:G_terminal}
	\mathcal{G}_T(\xx,\vv)
	:=
	\sum_{q\in\mathcal{I}}
	\left[
	\frac{\alpha_q}{2}
	\left|
	\psi_q(\xx,\vv)-\widehat{\psi}_q
	\right|^2
	+
	\frac{\beta_q}{2}
	\left|
	\psi_q(\xx,\vv)-m_q[f](T)
	\right|^2
	\right].
\end{equation}
The multiplier \(\lambda(t,\xx)\), on the other hand, does not satisfy a terminal condition, since it is determined at each time by the elliptic adjoint Poisson equation.
 In addition, periodic boundary conditions are imposed on $p$, while reflective boundary conditions are applied to $q$. 
Finally, we compute the derivative with respect to each control
coefficient. 
For any $k=1,\ldots,N_g$, the \(k\)-th component of the control gradient is
\begin{equation}\label{eq:gradient_continuos}
	\begin{split}
		&g_k(t)
		:=
		D_{b_k}\mathscr{L}(f,\phi,p,\lambda,\mathbf{b})
		\\
		&=
		\int_{\Omega_{\xx}}
		\left[
		\gamma b_k(t)\rho(t,\xx)
		+
		\int_{\Omega_{\vv}}
		\left(
		\vv\times\nabla_{\vv}p(t,\xx,\vv)
		\right)
		f(t,\xx,\vv)\,d\vv
		\right]
		\mathcal{K}_k(\xx)\,d\xx,
	\end{split}
\end{equation}
where, for two planar vectors $\mathbf{a}=(a_x,a_y)$ and $\mathbf{c}=(c_x,c_y)$, we denote by
$\mathbf{a}\times\mathbf{c} := a_x c_y - a_y c_x$ the scalar (out-of-plane) component of their
cross product.

The control optimality condition can thus be written as
\begin{equation}\label{eq:control_variational_inequality}
	\sum_{k=1}^{N_g}
	\int_0^T
	g_k^\star(t)
	\left(
	c_k(t)-b_k^\star(t)
	\right)\,dt
	\geq0
	\qquad
	\text{for every }
	\mathbf{c}\in\mathcal{U}_{\mathrm{ad}}.
\end{equation}
For an interior optimum, this condition reduces to
\(g_k^\star(t)=0\) for every \(k=1,\ldots,N_g\) and almost every
\(t\in(0,T)\). Altogether, equations \eqref{eq:forward_continuos},
\eqref{eq:backward_continuos}, \eqref{eq:terminal_continuos}, and
\eqref{eq:control_variational_inequality} define the continuous
forward--backward optimality system. The state equations evolve
forward from the prescribed initial distribution, whereas the adjoint
Vlasov equation evolves backward from the terminal condition. The
elliptic state and adjoint Poisson equations are solved at each time,
and the control coefficients are characterized by the variational
inequality \eqref{eq:control_variational_inequality}.

\subsection{From the continuous to the particle--grid optimality system}

We now derive a particle--grid counterpart of the continuous
optimality system
\eqref{eq:forward_continuos}--\eqref{eq:control_variational_inequality}.
Throughout this subsection, \(N_g\) continues to denote the number of
control cells, whereas \(N_h\) denotes the number of nodes of the
spatial grid used by the Particle-in-Cell (PIC) discretization \cite{PIC, filbet2016asymptotically,filbet2017asymptotically}. 
Let \(N\) be the number of particles and consider the regularized
empirical measure
\begin{equation}\label{eq:empirical_density}
	f^N(t,\xx,\vv)
	:=
	\sum_{i=1}^N
	\omega_i\,
	\mathcal{S}_h(\xx-\xx_i(t))
	\delta(\vv-\vv_i(t)),
\end{equation}
where the particle weights satisfy \(\omega_i>0\),
\(\sum_{i=1}^N\omega_i=\mathcal{M}^N\), and
\(\mathcal{S}_h\) is a normalized, compactly supported shape function
of width \(h\) \cite{hesthaven2024adaptive, PIC, filbet2016asymptotically}. The particle
positions and velocities are denoted by \((\xx_i(t),\vv_i(t))\), with
initial values chosen consistently with \(f_0\).

We introduce spatial grid nodes
\(\{\xx_\ell\}_{\ell=1}^{N_h}\), quadrature weights
\(\{\eta_\ell\}_{\ell=1}^{N_h}\), a discrete Laplacian \(A_h\), and a
discrete gradient \(G_h\). We then have
\begin{equation}\label{eq:discrete_rho}
	\rho_\ell^N(t)
	:=
	\sum_{i=1}^N
	\omega_i\mathcal{S}_h(\xx_\ell-\xx_i(t)),
	\qquad
	\ell=1,\ldots,N_h,
\end{equation}
while the discrete Poisson problem becomes
\begin{equation}\label{eq:discretize_Poisson}
	-(A_h\boldsymbol{\phi})_\ell
	=
	\rho_\ell^N-\bar{\rho},
	\qquad
	\ell=1,\ldots,N_h.
\end{equation}
The grid electric field and its interpolation to the particle positions are defined by
\begin{equation}\label{eq:grid_electric_field}
	\boldsymbol{E}_h
	:=
	-G_h\boldsymbol{\phi},
	\qquad
	\EE_i^h
	:=
	\sum_{\ell=1}^{N_h}
	\eta_\ell(\boldsymbol{E}_h)_\ell
	\mathcal{S}_h(\xx_\ell-\xx_i).
\end{equation}
The particle characteristics associated with the Vlasov equation are
\begin{equation}\label{eq:char_curves}
	\dot{\xx}_i=\vv_i,
	\qquad
	\dot{\vv}_i=F_i^h
	:=
	\EE_i^h
	+\vv_i\times B_i,
	\qquad
	B_i
	:=
	B_{\mathbf b}(t,\xx_i)
	=
	\sum_{k=1}^{N_g}b_k(t)\mathcal{K}_k(\xx_i),
\end{equation}
for \(i=1,\ldots,N\) while
\(\vv_i\times B_i=(v_{y,i}B_i,-v_{x,i}B_i)\).
To project the continuous adjoint equation onto the particle
trajectories, set
\[
\mathcal{P}:=\nabla_{\xx}p,
\qquad
\mathcal{Q}:=\nabla_{\vv}p.
\]
Taking the spatial and velocity gradients of the first equation in
\eqref{eq:backward_continuos} yields
\begin{equation}\label{eq:P_Q}
	\begin{split}
		&\partial_t\mathcal{P}
		+\vv\cdot\nabla_{\xx}\mathcal{P}
		+F_{\mathbf b}\cdot\nabla_{\vv}\mathcal{P}
		+(\nabla_{\xx}F_{\mathbf b})^T\mathcal{Q}
		=\nabla_{\xx}\mathcal{G}-\nabla_{\xx}\lambda,
		\\
		&\partial_t\mathcal{Q}
		+\vv\cdot\nabla_{\xx}\mathcal{Q}
		+F_{\mathbf b}\cdot\nabla_{\vv}\mathcal{Q}
		+\mathcal{P}
		+(\nabla_{\vv}F_{\mathbf b})^T\mathcal{Q}
		=\nabla_{\vv}\mathcal{G}.
	\end{split}
\end{equation}
Defining the particle adjoints by
\begin{equation}\label{eq:particle_adjoints}
	\pp_i(t)
	:=
	\mathcal{P}(t,\xx_i(t),\vv_i(t)),
	\qquad
	\qq_i(t)
	:=
	\mathcal{Q}(t,\xx_i(t),\vv_i(t)),
\end{equation}
and applying the chain rule along \eqref{eq:char_curves}, we obtain
\begin{equation}\label{eq:backward_equivalence}
	\begin{split}
		\dot{\pp}_i
		&=
		-(\nabla_{\xx_i}F_i^h)^T\qq_i
		+\nabla_{\xx_i}\mathcal{G}_i^N
		-\nabla_{\xx_i}\lambda_i^h,
		\\
		\dot{\qq}_i
		&=
		-\pp_i
		+\nabla_{\vv_i}\mathcal{G}_i^N
		+\qq_i\times B_i.
	\end{split}
\end{equation}
Now we set
\[
\psi_{\mathrm{x}}(\xx,\vv)=y,
\qquad
\psi_{\mathrm{v}}(\xx,\vv)=v_y,
\qquad
\widehat{\psi}_{\mathrm{x}}=\widehat y,
\qquad
\widehat{\psi}_{\mathrm{v}}=\widehat v_y,
\]
and
\[
\overline y^{\,N}
:=
\frac{1}{\mathcal{M}^N}
\sum_{j=1}^N\omega_jy_j,
\qquad
\overline v_y^{\,N}
:=
\frac{1}{\mathcal{M}^N}
\sum_{j=1}^N\omega_jv_{y,j}.
\]
The gradients of the particle approximation of the source
\(\mathcal{G}\) in \eqref{eq:G} become then
\begin{equation}\label{eq:nabla_G}
	\begin{split}
		\nabla_{\xx_i}\mathcal{G}_i^N
		={}&
		\left[
		\alpha_{\mathrm{x}}(y_i-\widehat y)
		+\beta_{\mathrm{x}}(y_i-\overline y^{\,N})
		\right]\ee_y
		+
		\frac{\gamma}{2}
		\sum_{k=1}^{N_g}
		|b_k(t)|^2\nabla\mathcal{K}_k(\xx_i),
		\\
		\nabla_{\vv_i}\mathcal{G}_i^N
		={}&
		\left[
		\alpha_{\mathrm{v}}(v_{y,i}-\widehat v_y)
		+\beta_{\mathrm{v}}(v_{y,i}-\overline v_y^{\,N})
		\right]\ee_y,
	\end{split}
\end{equation}
with \(\ee_y=(0,1)^T\). The discretized elliptic adjoint equation becomes
\begin{equation}\label{eq:backward_discrete_2}
	-A_h^*\widetilde{\boldsymbol{\lambda}}
	+G_h^*\boldsymbol{Q}
	=0.
\end{equation}
where \begin{equation}\label{eq:adjoint_deposition}
	\boldsymbol{Q}_\ell
	:=
	\sum_{i=1}^N
	\omega_i\mathcal{S}_h(\xx_\ell-\xx_i)\qq_i,
	\qquad
	\ell=1,\ldots,N_h,
\end{equation}
and 
\(\widetilde{\boldsymbol{\lambda}}
=(\widetilde\lambda_1,\ldots,\widetilde\lambda_{N_h})^T\).
The adjoints \(A_h^*\) and \(G_h^*\) use the
same quadrature-weighted grid inner product used in
\eqref{eq:discretize_Poisson}--\eqref{eq:grid_electric_field}, meaning that
on a uniform grid one has \(A_h^*=A_h^T=A_h\) and \(G_h^*=G_h^T\).
Thus, the particle adjoint equations in backward form are
\begin{equation}\label{eq:backward_discrete}
	\begin{split}
		-\dot{\pp}_i
		&=
		(\nabla_{\xx_i}F_i^h)^T\qq_i
		-\nabla_{\xx_i}\mathcal{G}_i^N
		+\nabla_{\xx_i}\lambda_i^h,
		\\
		-\dot{\qq}_i
		&=
		\pp_i
		-\nabla_{\vv_i}\mathcal{G}_i^N
		-\qq_i\times B_i,
	\end{split}
	\qquad i=1,\ldots,N.
\end{equation}
The terminal condition follows from \eqref{eq:terminal_continuos}:
\begin{equation}\label{eq:particle_terminal_conditions}
	\pp_i(T)
	=-\nabla_{\xx_i}\mathcal{G}_{T,i}^N,
	\qquad
	\qq_i(T)
	=-\nabla_{\vv_i}\mathcal{G}_{T,i}^N,
\end{equation}
where \(\mathcal{G}_{T,i}^N\) is obtained from the tracking part of
\(\mathcal{G}_i^N\), evaluated at \(T\); in particular, it contains no
control-regularization term.
Finally, inserting \eqref{eq:empirical_density} into the continuous
gradient \eqref{eq:gradient_continuos} gives, for
\(k=1,\ldots,N_g\),
\begin{equation}\label{eq:grid-gradient}
	\begin{split}
		g_k^{N,h}(t)
		={}&
		\sum_{i=1}^N\omega_i
		\int_{\Omega_{\xx}}
		\mathcal{S}_h(\xx-\xx_i)
		\left[
		\gamma b_k(t)
		+\vv_i\times\qq_i
		\right]
		\mathcal{K}_k(\xx)\,d\xx.
	\end{split}
\end{equation}
If \(\mathcal{S}_h\) has approximation order \(r\), then, for every
sufficiently smooth function \(\varphi\),
\[
\int_{\Omega_{\xx}}
\mathcal{S}_h(\xx-\xx_i)\varphi(\xx)\,d\xx
=
\varphi(\xx_i)+\mathcal{O}(h^r).
\]
Hence the particle control gradient is
\begin{equation}\label{eq:gradient_discrete}
	g_k^{N}(t)
	=
	\sum_{i=1}^N
	\omega_i
	\left[
	\gamma b_k(t)
	+\vv_i(t)\times\qq_i(t)
	\right]
	\mathcal{K}_k(\xx_i(t))
	+
	\mathcal{O}(h^r).
\end{equation}

The discrete control optimality condition is therefore
\begin{equation}\label{eq:discrete_control_variational_inequality}
	\sum_{k=1}^{N_g}
	\int_0^T
	(g_k^{N})^\star(t)
	\bigl(c_k(t)-b_k^\star(t)\bigr)\,dt
	\geq0
	\qquad
	\text{for every }
	\mathbf c\in\mathcal{U}_{\mathrm{ad}}^h,
\end{equation}
where $\mathcal{U}_{\mathrm{ad}}^h$ is defined in \eqref{eq:Uadm_h}.

In the unconstrained case, or when the discrete amplitude constraint
is inactive, this reduces formally to \(g_k^N(t)=0\).
Equations \eqref{eq:char_curves}, \eqref{eq:discretize_Poisson},
\eqref{eq:backward_discrete}, \eqref{eq:backward_discrete_2},
\eqref{eq:particle_terminal_conditions}, \eqref{eq:discrete_control_variational_inequality} form the
particle--grid optimality system. If
the same shape function and quadrature are used, and if the discrete
differential operators are paired with their adjoints, this projected
optimize-then-discretize system agrees formally with the system
obtained by discretizing the particle--grid dynamics first and then
differentiating it, up to particle, grid, and time-discretization
errors. A detailed comparison is provided in
Appendix~\ref{app:computations_optimality}.

\section{Numerical approximation and optimization method}
\label{sec:numerical_optimization}

In this section, we describe the fully discrete strategy used to solve
the particle--grid optimality system
\eqref{eq:char_curves}, \eqref{eq:discretize_Poisson},
\eqref{eq:backward_discrete}, and \eqref{eq:backward_discrete_2},
with terminal conditions \eqref{eq:particle_terminal_conditions} and
control optimality condition
\eqref{eq:discrete_control_variational_inequality}. Optimisation 
is performed using a projected-gradient method with Armijo backtracking,
following standard approaches for particle and mean-field control
problems \cite{albi2017mean,bailo2018optimal}.
Let
\[
j_{N,h,\Delta t}(\mathbf b)
:=
\mathcal{J}_{N,h,\Delta t}
\bigl(f_{N,h,\Delta t}[\mathbf b],\mathbf b\bigr)
\]
denote the reduced, fully discrete cost functional, where
\(f_{N,h,\Delta t}[\mathbf b]\) is the particle approximation obtained
by solving the discrete state equations for the control \(\mathbf b\).
The time integrals in \(\mathcal{J}\) are approximated using the same
time quadrature used to define the discrete gradient below.
We initialize the optimization with
\[
b_k^{(0)}(t)=b_k^0,
\qquad k=1,\ldots,N_g,
\]
where the coefficients \(b_k^0\) are obtained from the instantaneous
control described in Appendix~\ref{app:instantaneous_control}.
 In optimisation iteration \(m\), the current coefficients
\(\mathbf b^{(m)}\) are kept fixed while the particle--grid state
system is integrated forward in time. This produces
\((\xx_i^{(m)},\vv_i^{(m)})\). The adjoint equations are then solved
backward along the same trajectory, producing
\((\pp_i^{(m)},\qq_i^{(m)})\), and the gradient
\(\mathbf g^{(m)}\) is evaluated from
\eqref{eq:gradient_discrete}. In particular, the magnetic field acting
on particle \(i\) is
\begin{equation}\label{eq:Bi}
	B_i^{(m)}(t)
	=
	B_{\mathbf b^{(m)}}(t,\xx_i^{(m)}(t))
	=
	\sum_{k=1}^{N_g}
	b_k^{(m)}(t)\mathcal K_k(\xx_i^{(m)}(t)).
\end{equation}
After time discretization, let
\(\mathcal U_{\mathrm{ad}}^{h,\Delta t}\) denote the discrete-time
counterpart of \(\mathcal U_{\mathrm{ad}}^h\), namely
\[
\mathcal U_{\mathrm{ad}}^{h,\Delta t}
:=
\left\{
\bigl(\mathbf b^n\bigr)_{n=0}^{N_t-1}:
\left|
\sum_{k=1}^{N_g}b_k^n\mathcal K_k(\xx_\ell)
\right|
\leq B_{\rm max},\quad
\ell=1,\ldots,N_h,\quad n=0,\ldots,N_t-1
\right\}.
\]
We denote by
\(\Pi_{\mathcal U_{\mathrm{ad}}^{h,\Delta t}}\) the orthogonal
projection onto this closed convex set with respect to the discrete
\(L^2\) inner product. For a trial step size \(\alpha>0\), define
\begin{equation}\label{eq:projected_trial_control}
	\mathbf b_\alpha^{(m)}
	:=
	\Pi_{\mathcal U_{\mathrm{ad}}^{h,\Delta t}}
	\left(
	\mathbf b^{(m)}-\alpha\mathbf g^{(m)}
	\right).
\end{equation}
The accepted step size is
\begin{equation}\label{eq:Armijo_rule}
	\alpha^{(m)}
	=
	\alpha_0\eta^{r_m},
	\qquad
	\eta\in(0,1),
\end{equation}
where \(r_m\) is the smallest nonnegative integer for which
\begin{equation}\label{eq:Armijo_condition}
	\begin{split}
		j_{N,h,\Delta t}
		\bigl(\mathbf b_{\alpha^{(m)}}^{(m)}\bigr)
		\leq{}&
		j_{N,h,\Delta t}\bigl(\mathbf b^{(m)}\bigr)
		\\
		&+
		c\left\langle
		\mathbf g^{(m)},
		\mathbf b_{\alpha^{(m)}}^{(m)}
		-\mathbf b^{(m)}
		\right\rangle_{\Delta t},
		\qquad c\in(0,1),
	\end{split}
\end{equation}
is satisfied,  where
\begin{equation*}
    \langle \mathbf{g}, \mathbf{b} \rangle_{\Delta t} = \Delta t \sum_{n=0}^{N_t-1} \mathbf{g}^n \cdot \mathbf{b}^n
\end{equation*}  
denotes the discrete time inner product approximating the continuous $L^2(0,T)$ inner product via a Riemann sum, with $\mathbf{g}^n, \mathbf{b}^n \in \mathbb{R}^{N_c^x \times N_c^y}$ denoting the gradient and control arrays at time $t^n$ respectively. The control update is then
\begin{equation}\label{eq:B_update}
	\mathbf b^{(m+1)}
	=
	\mathbf b_{\alpha^{(m)}}^{(m)}.
\end{equation}
Each trial evaluation in the line search requires a new forward solve,
but the adjoint system and the gradient are computed only once at the
beginning of iteration \(m\) \cite{conn2000toint,kelley1999iterative}.
The iteration stops either when the maximum number of iterations  is reached, or when the cost functional stagnates, that is if
\begin{equation*}
    \left| j_{N,h,\Delta t}\bigl(\mathbf{b}^{(m)}\bigr) - j_{N,h,\Delta t}\bigl(\mathbf{b}^{(m-1)}\bigr) \right| \leq \varepsilon_{\mathrm{stall}},
\end{equation*} 
holds for $n_{\mathrm{stall}}$ consecutive iterations. 

To reduce the computational cost, the state, adjoint, and gradient are evaluated on a representative subset
\(\mathcal I_s\subset\{1,\ldots,N\}\), with
\(|\mathcal I_s|=N_s\ll N\). Assuming equal particle weights and
relabeling the particles in the reduced ensemble by
\(j=1,\ldots,N_s\), we assign the rescaled weights
\[
\widehat\omega_j
:=
{\mathcal M^N}/{N_s},
\qquad j=1,\ldots,N_s,\quad\textrm{ with }\quad \sum_{j=1}^{N_s}\widehat\omega_j=\mathcal M^N.
\]
Equivalently, if the full ensemble has weights
\(\omega_i=\mathcal M^N/N\), then $
\widehat\omega_j={N}/{N_s}\omega_j$. Omitting for simplicity in the following the dependence on optimization iteration $(m)$, we partition \([0,T]\) into \(N_t\) intervals of length
\(\Delta t=T/N_t\). The controls are then represented as piecewise constant
functions in time, with  $b_k^{n} \simeq b_k^{}(t^n)$ while
the Poisson problem is solved according to
\begin{equation}\label{eq:rho_E_discrete}
	\begin{split}
		\rho_\ell^{n}
		&=
		\sum_{j=1}^{N_s}
		\widehat\omega_j
		\mathcal S_h
		(\xx_\ell-\xx_j^{n}),
		\\
		-(A_h\boldsymbol\phi^{n})_\ell
		&=
		\rho_\ell^{n}-1,
		\qquad
		\ell=1,\ldots,N_h,
		\\
		\boldsymbol E_h^{n}
		&=
		-G_h\boldsymbol\phi^{n},
		\quad
		\EE_i^{h,n}
		=
		\sum_{\ell=1}^{N_h}
		\eta_\ell
		(\boldsymbol E_h^{n})_\ell
		\mathcal S_h
		(\xx_\ell-\xx_i^{n}).
	\end{split}
\end{equation}
The controlled magnetic field at the particle position becomes
\begin{equation}\label{eq:Bi_discrete}
B_i^{n}=
\sum_{k=1}^{N_g}
b_k^{n}
\mathcal K_k(\xx_i^{n}).
\end{equation}
 
The state characteristics in the reduced system ($N_s \ll N$) are advanced using a first-order semi-implicit scheme \cite{filbet2016asymptotically}:
\begin{equation}\label{eq:forward_discretization}
	\begin{cases}
		\displaystyle
		\vv_i^{n+1}
		=
		\vv_i^{n}
		+
		\Delta t
		\left(
		\EE_i^{h,n}
		+
		\vv_i^{n+1}\times B_i^{n}
		\right),
		\\[0.8em]
		\displaystyle
		\xx_i^{n+1}
		=
		\xx_i^{n}
		+
		\Delta t\,\vv_i^{n+1},
	\end{cases}
	\qquad
	i=1,\ldots,N_s.
\end{equation}
In the same way, the particle adjoints are initialized by
\begin{equation}\label{eq:backward_terminal_discrete_time}
	\pp_i^{N_t}
	=
	-\nabla_{\xx_i}\mathcal G_{T,i}^{N},
	\qquad
	\qq_i^{N_t}
	=
	-\nabla_{\vv_i}\mathcal G_{T,i}^{N},
	\qquad i=1,\ldots,N_s,
\end{equation}
while for \(n=N_t-1,\ldots,0\), we employ a semi-implicit backward step for the
particle adjoint equations:
\begin{equation}\label{eq:backward_discretization}
	\begin{cases}
		\displaystyle
		\pp_i^{n}
		=
		\pp_i^{n+1}
		+
		\Delta t
		\left[
		J_{F,i}^{n+1}\qq_i^{n}
		-
		\nabla_{\xx_i}\mathcal G_i^{N,n+1}
		+
		L_i^{n+1}
		\right],
		\\[0.9em]
		\displaystyle
		\qq_i^{n}
		=
		\qq_i^{n+1}
		+
		\Delta t
		\left[
		\pp_i^{n}
		-
		\nabla_{\vv_i}\mathcal G_i^{N,n+1}
		-
		\qq_i^{n}\times B_i^{n+1}
		\right],
		\\[0.9em]
		\displaystyle
		-A_h^*
		\widetilde{\boldsymbol\lambda}^{n+1}
		+
		G_h^*\boldsymbol Q^{n+1}
		=0.
	\end{cases}
\end{equation}
with
\[
J_{F,i}^{n+1}
:=
\left(
\nabla_{\xx_i}F_i^{h,n+1}
\right)^T,
\qquad
F_i^{h,n+1}
=
\EE_i^{h,n+1}
+
\vv_i^{n+1}\times B_i^{n+1}.
\]
and
\begin{align}
	\boldsymbol Q_\ell^{n+1}
	&=
	\sum_{i=1}^{N_s}
	\widehat\omega_i
	\mathcal S_h(\xx_\ell-\xx_i^{n+1})
	\qq_i^{n+1},
	\qquad
	\ell=1,\ldots,N_h,
	\label{eq:Q_discrete_time}
	\\
	L_i^{n+1}
	&:=
	\nabla_{\xx_i}\lambda_i^{h,n+1}
	=
	\sum_{\ell=1}^{N_h}
	\eta_\ell
	\widetilde\lambda_\ell^{n+1}
	\nabla_{\xx_i}
	\mathcal S_h(\xx_\ell-\xx_i^{n+1}).
	\label{eq:L_discrete_time}
\end{align}
The discrete gradient with respect to the control coefficients
is finally computed as
\begin{equation}\label{eq:gradient_discrete_time}
	g_{k,n}^{}
	:=
	\sum_{i=1}^{N_s}
	\widehat\omega_i
	\left[
	\gamma b_k^{n}
	+
	\vv_i^{n}\times\qq_i^{n}
	\right]
	\mathcal K_k(\xx_i^{n}),
\end{equation}

We emphasize that \eqref{eq:gradient_discrete_time} is obtained by
time-discretizing the semidiscrete particle adjoint system. Its relation
with the gradient obtained by differentiating the fully discrete state map
is discussed in Appendix~\ref{app:computations_optimality}.
 
While the first-order scheme \eqref{eq:forward_discretization} yields a simple and efficient formulation for the iterative optimization phase, it exhibits numerical dissipation over long times.  To avoid artificial damping when computing physical quantities, all validation tests on the full ensemble are performed using an energy-preserving second-order semi-implicit scheme \cite{filbet2016asymptotically}.  This reads
\begin{equation}\label{eq:second_order_first_stage}
	  	\begin{cases}
        &\vv_i^{(1)} = \vv_i^n + \frac{\Delta t}{2} \left(\EE_i^{h,n} +  \vv_i^{(1)} \times B_i^n   \right),\\
		  		&\xx_i^{(1)} = \xx_i^n + \frac{\Delta t}{2} \vv_i^{(1)},
		  	\end{cases}
	\qquad
	i=1,\ldots,N,
	  \end{equation}
where the apex $(1)$ denotes the first stage, while the second stage reads 
\begin{equation}\label{eq:second_order_second_stage}
		\begin{cases}	
        	&\vv_i^{(2)} = \vv_i^n + \frac{\Delta t}{2} \left( \EE_i^{h,\star} + \vv_i^{(2)} \times B_i^\star   \right),\\
				&\xx_i^{(2)} = \xx_i^n + \frac{\Delta t}{2} \vv_i^{(2)}, 	 
			\end{cases}
	\qquad
	i=1,\ldots,N,
	\end{equation}
    where the electric field $\EE_i^{h,\star}$ is computed as in \eqref{eq:rho_E_discrete} by solving the Poisson problem, whereas the magnetic field $B_i^\star$ is computed as in \eqref{eq:Bi_discrete} by replacing $\xx_i^n$ with $\xx_i^\star = 2 \xx_i^{(1)} - \xx_i^n$.
	Finally, the numerical solution at time $t^{n+1}$ reads,
\begin{equation}\label{eq:second_order_scheme}
	\begin{cases}
		&\vv_i^{n+1} = \vv_i^{(1)}+\vv_i^{(2)}-\vv_i^n,\\
		&\xx_i^{n+1} = \xx_i^{(1)}+\xx_i^{(2)}-\xx_i^n, 
	\end{cases}
	\qquad
	i=1,\ldots,N.
	\end{equation} 
 
 \section{Numerical experiments}\label{sec:numexp}

In this section, we assess the performance of the proposed control strategy on a set of standard benchmark problems. For the Armijo line search in \eqref{eq:Armijo_rule}–\eqref{eq:Armijo_condition}, we set the parameters to $\alpha_0 = 10$, $\eta = 0.5$, and $c = 10^{-4}$.  The stopping tolerance for the optimization algorithm is set to $\varepsilon_{\mathrm{stall}}=10^{-5}$ over a stall limit of $n_{stall} = 50$ iterations.  A full ensemble of $N=10^7$ particles is sampled from the initial distribution, whereas the forward--backward optimality system and the corresponding gradient are evaluated using a fixed representative random subset of $N_s=10^3$ particles at each optimization iteration. The rest of the parameters are specified separately for each test case. As described in Section~\ref{sec:problem_setting}, the control is represented on a coarse rectangular partition consisting of
$
N_g=N_c^xN_c^y
$
cells, with $N_c^x$ and $N_c^y$ cells along the $x$- and $y$-directions, respectively.
%
We denote the corresponding cell centers by
$
\xx_c^{j,k}=(x_c^j,y_c^k).
$
The control-cell widths are denoted by $h_x^c$ and $h_y^c$. They are chosen substantially larger than the spacings $\Delta x$ and $\Delta y$ of the finer grid used in the Particle-in-Cell reconstruction, namely
$
h_x^c\gg\Delta x,
\
h_y^c\gg\Delta y.
$
Using a double-index representation
$b_{j,k}(t)$ of the control coefficients introduced in
Section~\ref{sec:problem_setting}, the magnetic field evaluated at the
$i$-th particle position is given by
\begin{equation}\label{eq:Bi_cells_experiments}
	B_i(t)
	:=
	B_{\mathbf b}(t,\xx_i(t))
	=
	\sum_{j=1}^{N_c^x}
	\sum_{k=1}^{N_c^y}
	b_{j,k}(t)
	\mathcal K_{j,k}(\xx_i(t)).
\end{equation}
The two-dimensional control kernels are defined as tensor products,
\begin{equation}\label{eq:Phi_control}
	\mathcal K_{j,k}(x,y)
	=
	\kappa_{j}^{x}(x)
	\kappa_{k}^{y}(y),
\end{equation}
where
\begin{align}
	\kappa_j^x(x)
	&=
	\sigma\left(
	\frac{x-x_{j-\frac12}^c}{\varepsilon_c}
	\right)
	-
	\sigma\left(
	\frac{x-x_{j+\frac12}^c}{\varepsilon_c}
	\right),
	\\
	\kappa_k^y(y)
	&=
	\sigma\left(
	\frac{y-y_{k-\frac12}^c}{\varepsilon_c}
	\right)
	-
	\sigma\left(
	\frac{y-y_{k+\frac12}^c}{\varepsilon_c}
	\right),
\end{align}
with
$
\sigma(z)=\frac{1}{1+\exp(-z)}
$
and $\varepsilon_c = 10^{-10}$. 
Once the optimal coefficients $\mathbf b^\star(t)$ have been computed
over the time interval $[0,T]$ using the reduced ensemble of $N_s$
particles, the resulting control is validated by evolving the full
ensemble of $N=10^7$ particles using the second order scheme in \eqref{eq:second_order_first_stage}-\eqref{eq:second_order_second_stage}-\eqref{eq:second_order_scheme}.

To assess the effectiveness of the control strategy near the material
walls, we introduce a boundary strip
$$
\Omega_b
=
[a_x,b_x]
\times
\left(
[a_y,a_y+\delta_b]
\cup
[b_y-\delta_b,b_y]
\right)
\subset\Omega_{\xx},
$$
where $\delta_b$ denotes the thickness of the strip. Unless otherwise
specified, we set $\delta_b=\Delta y$.
Let $\{\xx_\ell\}_{\ell=1}^{N_h}$ denote the nodes of the fine spatial
grid used for the PIC reconstruction, and define

$$
\mathcal I_b
:=
\left\{
\ell\in\{1,\ldots,N_h\}:
\xx_\ell\in\Omega_b
\right\},
\qquad
N_b:=|\mathcal I_b|.
$$

The particle mass contained in the boundary strip is approximated by
\begin{equation}\label{eq:boundary_mass}
	\mathcal M_b(t)
	=
	\Delta x \Delta y
	\sum_{\ell\in\mathcal I_b}
	\sum_{i=1}^{N}
	\omega_i
	\mathcal S_h(\xx_\ell-\xx_i(t)),
\end{equation}
where $S_h$ is the two-dimensional tensor-product quadratic B-spline shape function scaled by the spatial resolution parameter $h = \Delta x \Delta y$.
The corresponding mean velocity is
\begin{equation}\label{eq:boundary_mean_velocity}
	U_b(t)
	=
	\frac{\Delta x \Delta y}{\mathcal M_b(t)}
	\sum_{\ell\in\mathcal I_b}
	\sum_{i=1}^{N}
	\omega_i\vv_i(t)
	\mathcal S_h(\xx_\ell-\xx_i(t)).
\end{equation}
We then define the thermal energy contained in the boundary strip as
\begin{equation}\label{eq:boundary_energy}
	\mathcal E_b(t)
	=
	\frac{\Delta x \Delta y}{2}
	\sum_{\ell\in\mathcal I_b}
	\sum_{i=1}^{N}
	\omega_i
	\left|
	\vv_i(t)-U_b(t)
	\right|^2
	\mathcal S(\xx_\ell-\xx_i(t)).
\end{equation}
For comparison, the total reconstructed mass and the global mean
velocity are
\begin{align}
	\mathcal M(t)
	&=
	\Delta x \Delta y
	\sum_{\ell=1}^{N_h}
	\sum_{i=1}^{N}
	\omega_i
	\mathcal S(\xx_\ell-\xx_i(t)),
	\\
	U(t)
	&=
	\frac{\Delta x\Delta y}{\mathcal M(t)}
	\sum_{\ell=1}^{N_h}
	\sum_{i=1}^{N}
	\omega_i\vv_i(t)
	\mathcal S(\xx_\ell-\xx_i(t)).
\end{align}
The total thermal energy is consequently defined as
\begin{equation}\label{eq:total_thermal_energy}
	\mathcal E(t)
	=
	\frac{\Delta x\Delta y}{2}
	\sum_{\ell=1}^{N_h}
	\sum_{i=1}^{N}
	\omega_i
	\left|
	\vv_i(t)-U(t)
	\right|^2
	\mathcal S(\xx_\ell-\xx_i(t)).
\end{equation}
Finally, the percentage of the total thermal energy localized within
the boundary strip is also measured in the sequel, and denoted by $\mathcal{E}_b^{rel}(t)$.

\subsection{Two-dimensional double Sod test}

We consider a two-dimensional extension of the double Sod shock-tube
test. The spatial and velocity variables are denoted by
$\xx=(x,y)$ and $\vv=(v_x,v_y)$, respectively. The spatial domain is
$\Omega_{\xx}=[0,1.5]^2$, with periodic boundary conditions in the
$x$-direction and reflective boundary conditions in the $y$-direction.
The initial distribution is given by
\begin{equation}\label{eq:sod_2D_initial_distr}
	f_0(\xx,\vv)
	=
	\frac{\rho_0(\xx)}{2\pi T_0(\xx)}
	\exp\left(
	-\frac{|\vv|^2}{2T_0(\xx)}
	\right),
\end{equation}
where the initial density and temperature are defined as
\begin{equation}\label{eq:rho0_T0_2D} \begin{aligned} \rho_0(\xx) &= \begin{cases} 0.125, & y\in[0,0.5)\cup(1,1.5],\\ 1, & y\in[0.5,1], \end{cases} \\[0.5em] T_0(\xx) &= \begin{cases} 0.1, & y\in[0,0.5)\cup(1,1.5],\\ 1, & y\in[0.5,1]. \end{cases} \end{aligned} \end{equation}
Thus, the initial configuration consists of a dense, high-temperature
central layer surrounded by two low-density, low-temperature regions.
Although the problem is posed in a two-dimensional domain, the initial
data are homogeneous in the $x$-direction.

\begin{figure}[!htb]
	\centering
	\includegraphics[width=0.35\linewidth]
	{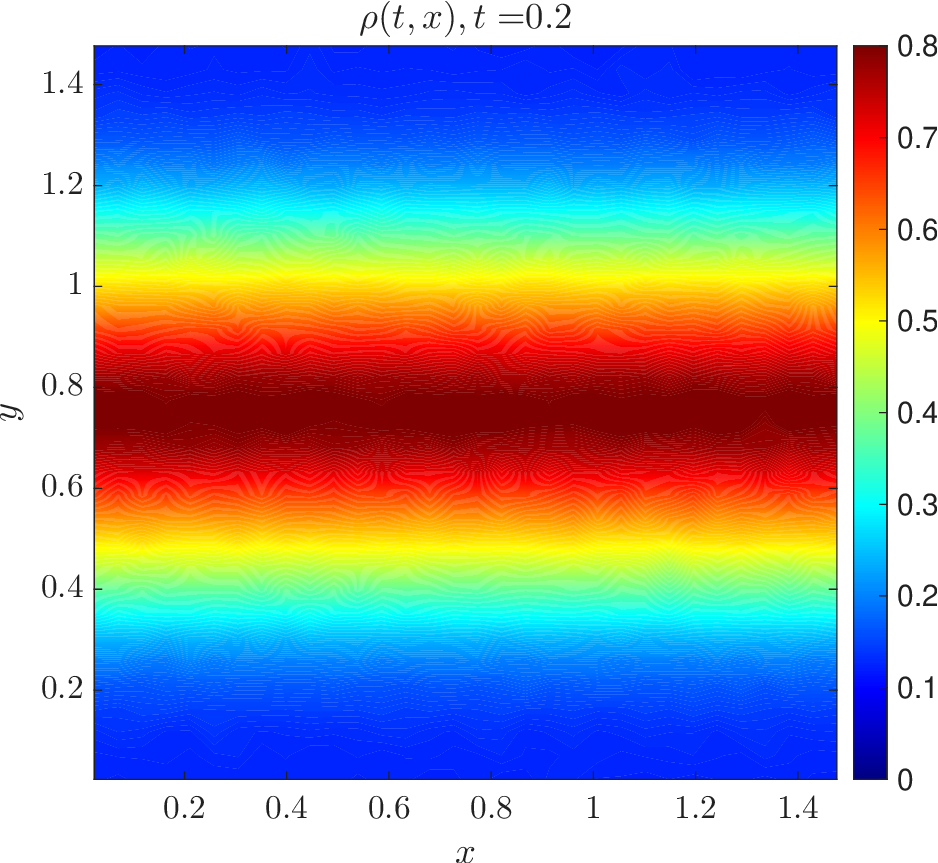}
	\includegraphics[width=0.35\linewidth]
	{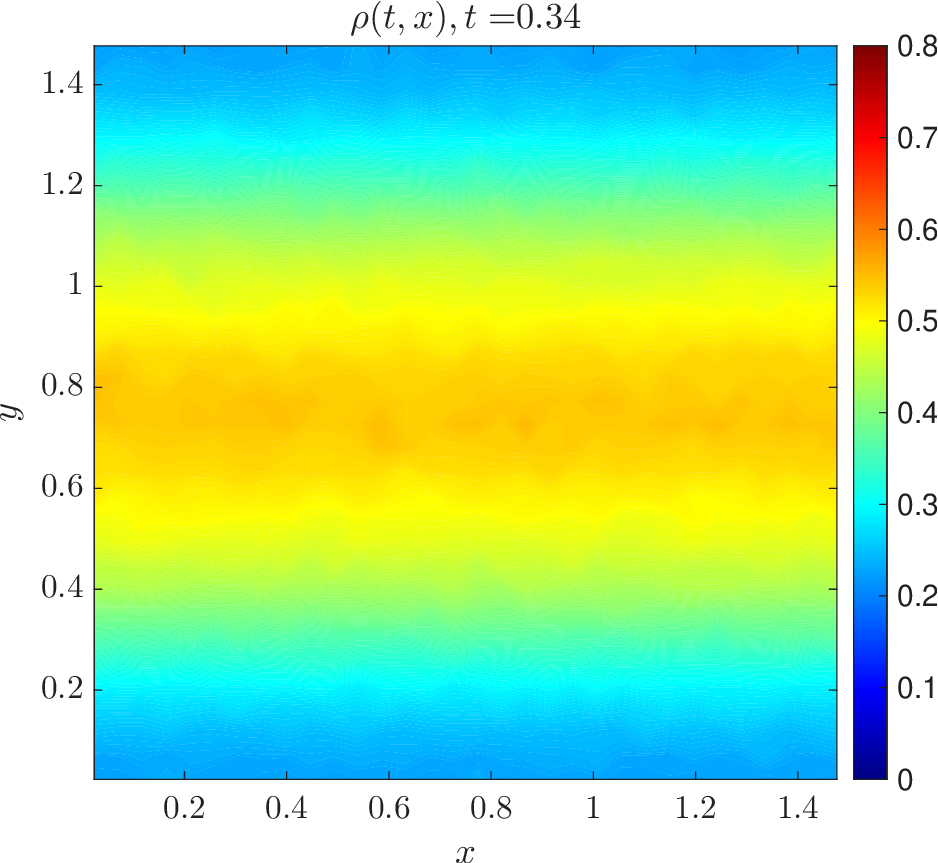}
	\includegraphics[width=0.35\linewidth]
	{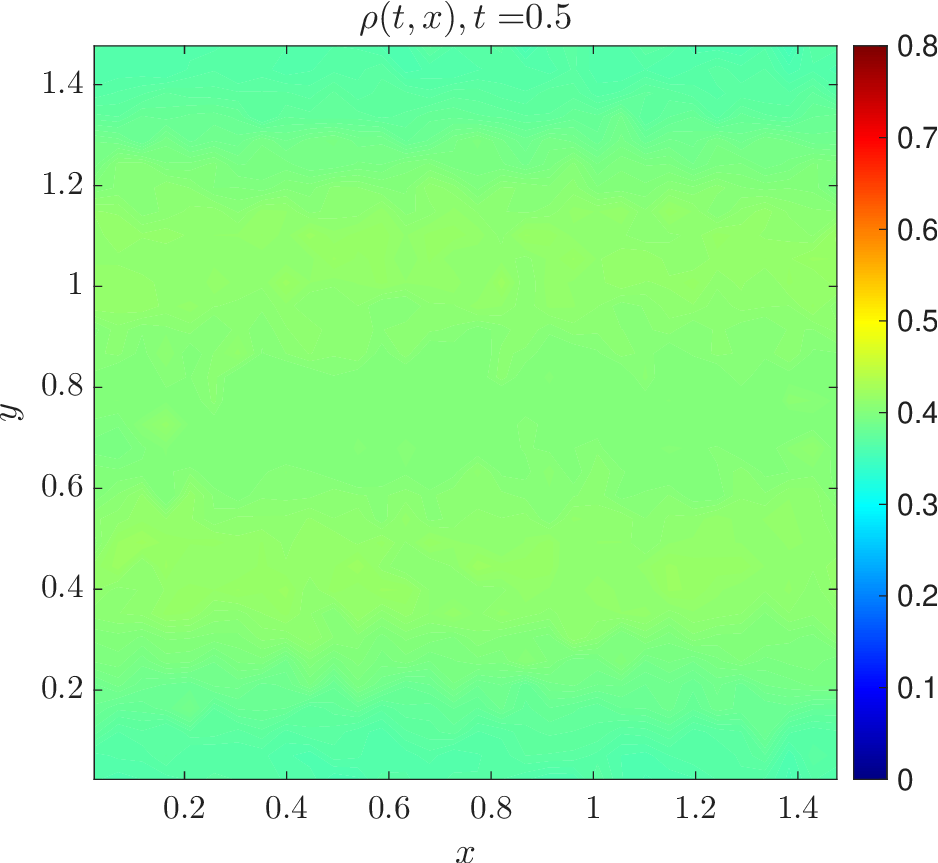}	
	\includegraphics[width=0.35\linewidth]
	{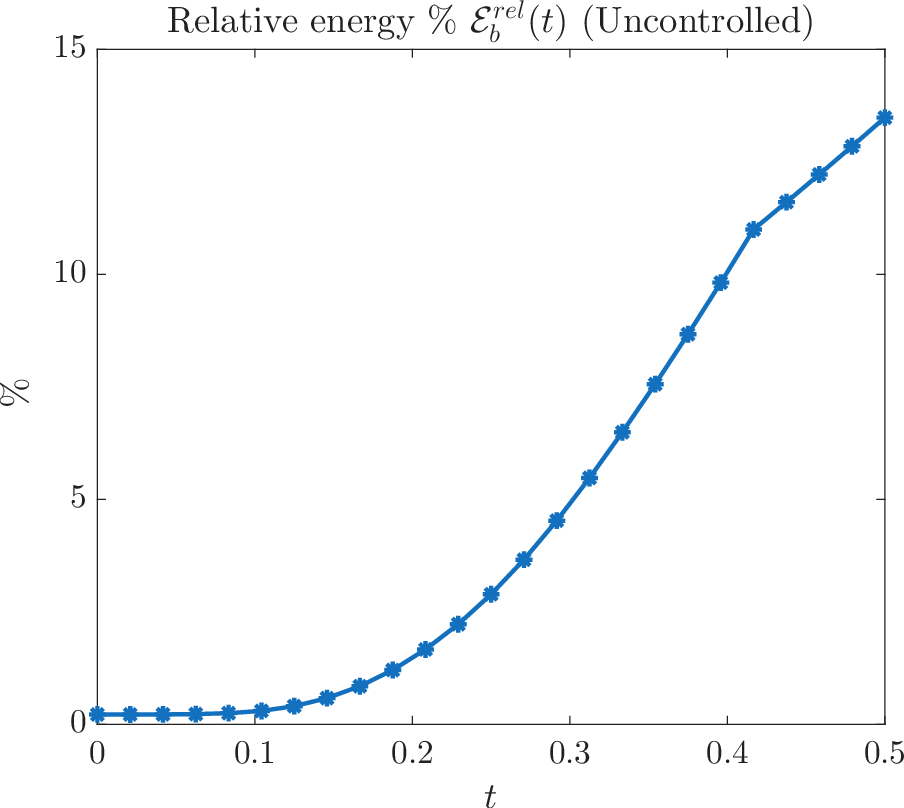}
	\caption{Two-dimensional double Sod test with the prescribed
		uniform magnetic field $B_{\mathrm{ref}}=0.75$. Snapshots of
		the uncontrolled dynamics at $t=0.2$, $t=0.34$, and $t=0.5$ and evolution of the
		percentage of total thermal energy contained in the boundary strip,
		$\mathcal E_b^{\mathrm{rel}}(t)$.}
	\label{fig:sod_uncontrolled_energy}
\end{figure}

The dynamics is evolved up to the final time $T=0.5$, using the time
step $\Delta t=0.01$. The density and the self-consistent electric
field are reconstructed on a uniform PIC grid with
$m_x=m_y=32$ nodes, so that $N_h=m_xm_y$. For both the instantaneous
and optimal controls, we set
\[ 
\alpha_{\mathrm{x}}
=
\beta_{\mathrm{x}}
=
\alpha_{\mathrm{v}}
=
\beta_{\mathrm{v}}
=
10,
\qquad
\gamma=10^{-3} , \qquad M = 50,
\]
with targets $\hat{y} = 0.75$ and $\hat{v}_y = 0$. 
We first consider a reference configuration in which no optimization
is performed and the magnetic field is prescribed to be uniform in
space and constant in time: $
B_{\mathrm{ref}}(t,\xx)=0.75.
$
Figure~\ref{fig:sod_uncontrolled_energy} shows the percentage of the
total thermal energy contained in the boundary strip,
$\mathcal E_b^{\mathrm{rel}}(t)$ together with three snapshots of the density profiles. The increase observed during the
simulation indicates that an increasing fraction of the thermal
energy becomes localized near the upper and lower boundaries.

We next consider the controlled dynamics, with the control computed
using the optimization procedure described in the previous sections
and initialized by the instantaneous control. We set $N_c^x=1$ and
consider $N_c^y\in\{5,7,9\}$. Consequently, each control cell extends
over the entire $x$-direction and the reconstructed magnetic field
depends only on $y$ and $t$.

\begin{figure}[!htb]
	\centering
	\includegraphics[width=0.35\linewidth]
	{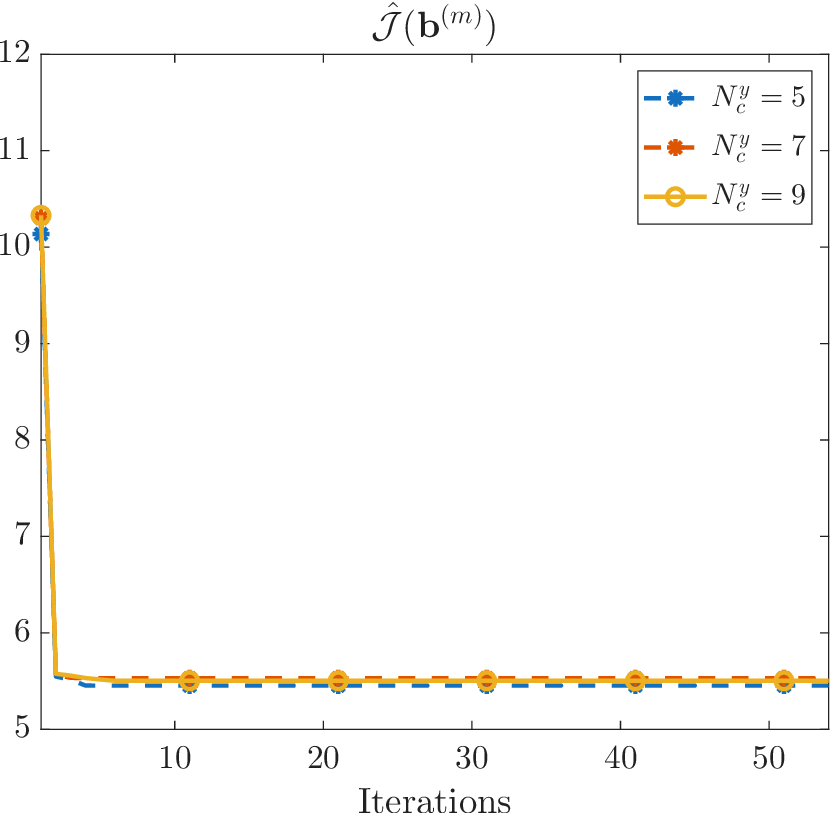}\hspace{1cm}
	\includegraphics[width=0.36\linewidth]
	{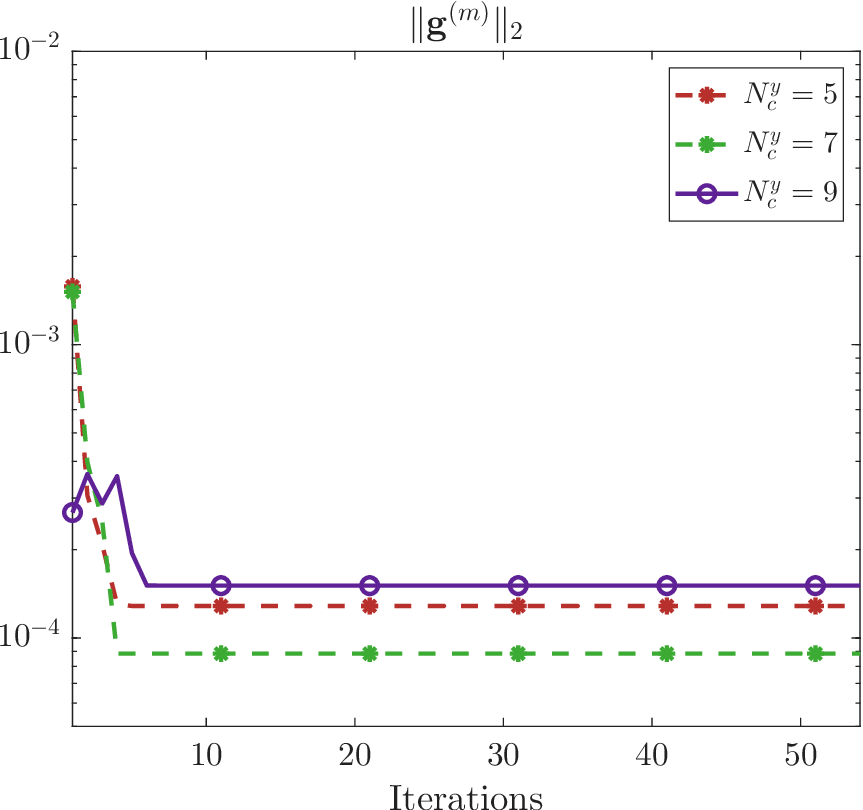}
	\caption{Two-dimensional double Sod test. Left: evolution of the
		reduced cost functional
		$\widehat{\mathcal J}(\mathbf b^{(m)})$ over the optimization
		iterations. Right: evolution of the discrete gradient norm
		$\Vert \mathbf g^{(m)}\Vert_2$. The curves correspond to
		$N_c^y\in\{5,7,9\}$.}
	\label{fig:functional}
\end{figure}

\begin{figure}[!htb]
	\centering
	\includegraphics[width=0.327\linewidth]
	{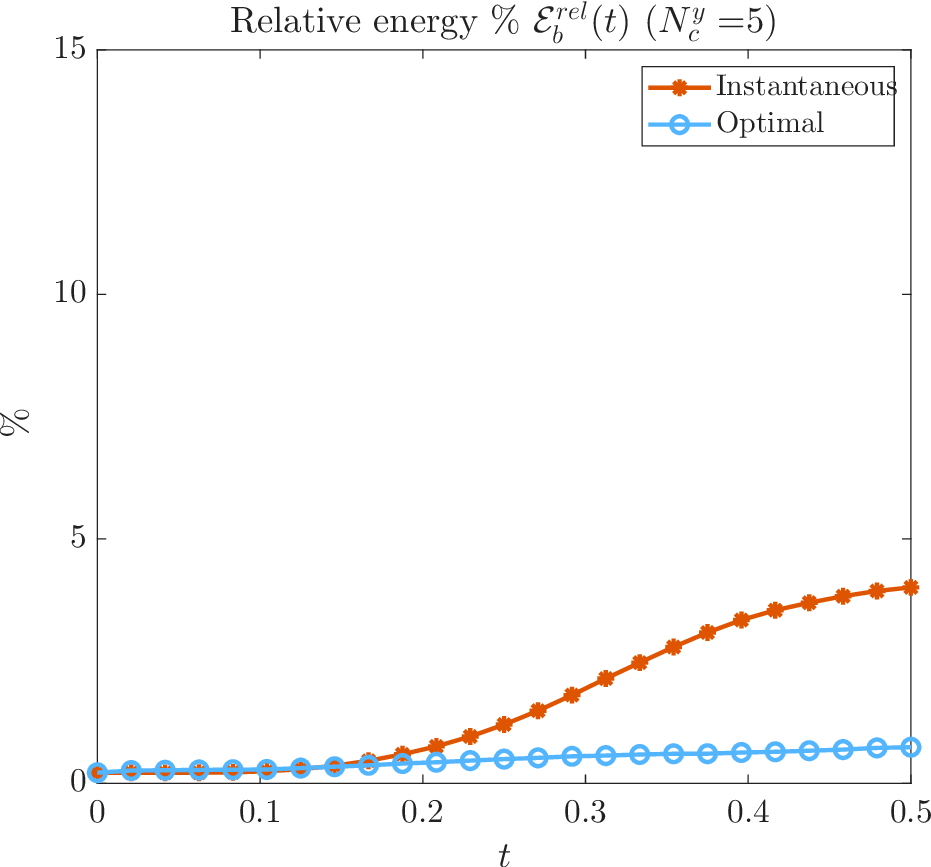}
	\includegraphics[width=0.327\linewidth]
	{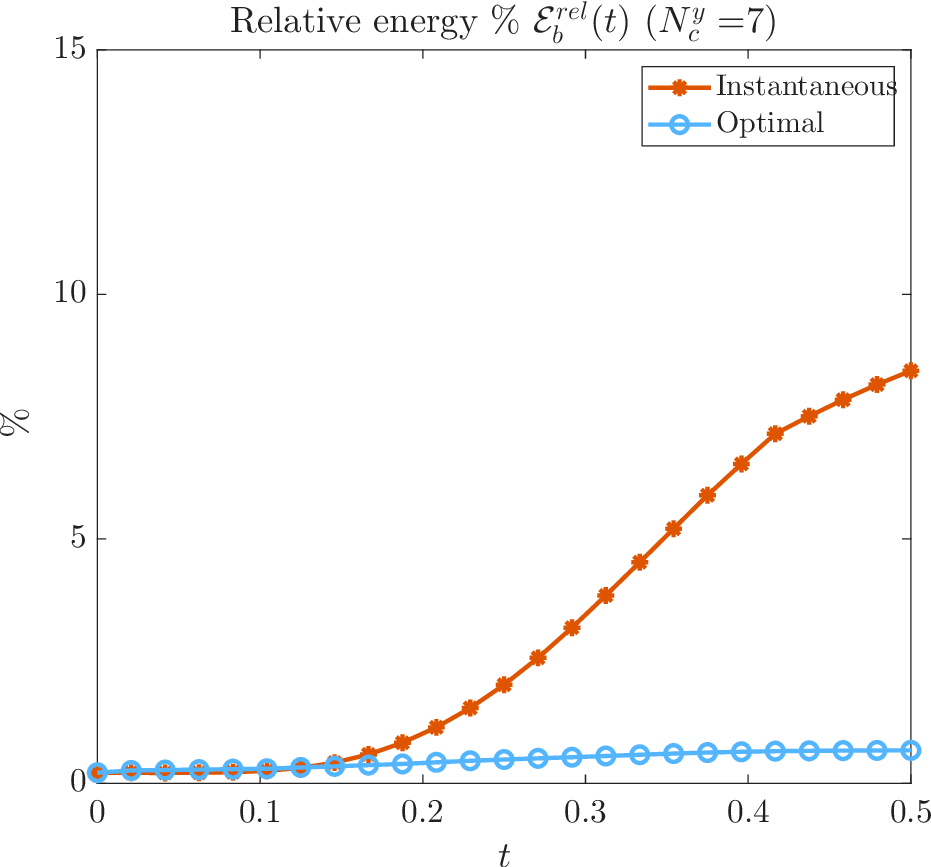}
	\includegraphics[width=0.327\linewidth]
	{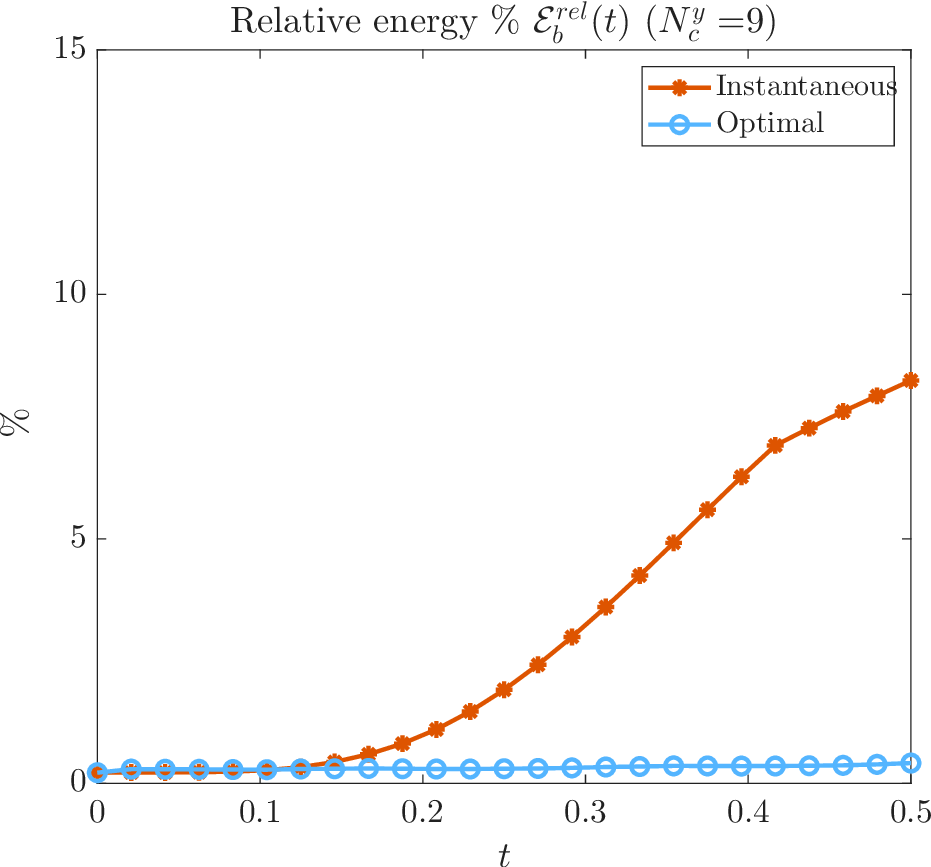}
	\vspace*{5pt}
	\caption{Two-dimensional double Sod test. Comparison of the
		percentage of thermal energy contained in the boundary strip,
		$\mathcal E_b^{\mathrm{rel}}(t)$, obtained with the instantaneous
		and optimal controls. From left to right:
		$N_c^y=5$, $N_c^y=7$, and $N_c^y=9$.}
	\label{fig:energy}
\end{figure}

\begin{figure}[!htb]
	\centering
	\includegraphics[width=0.327\linewidth]
	{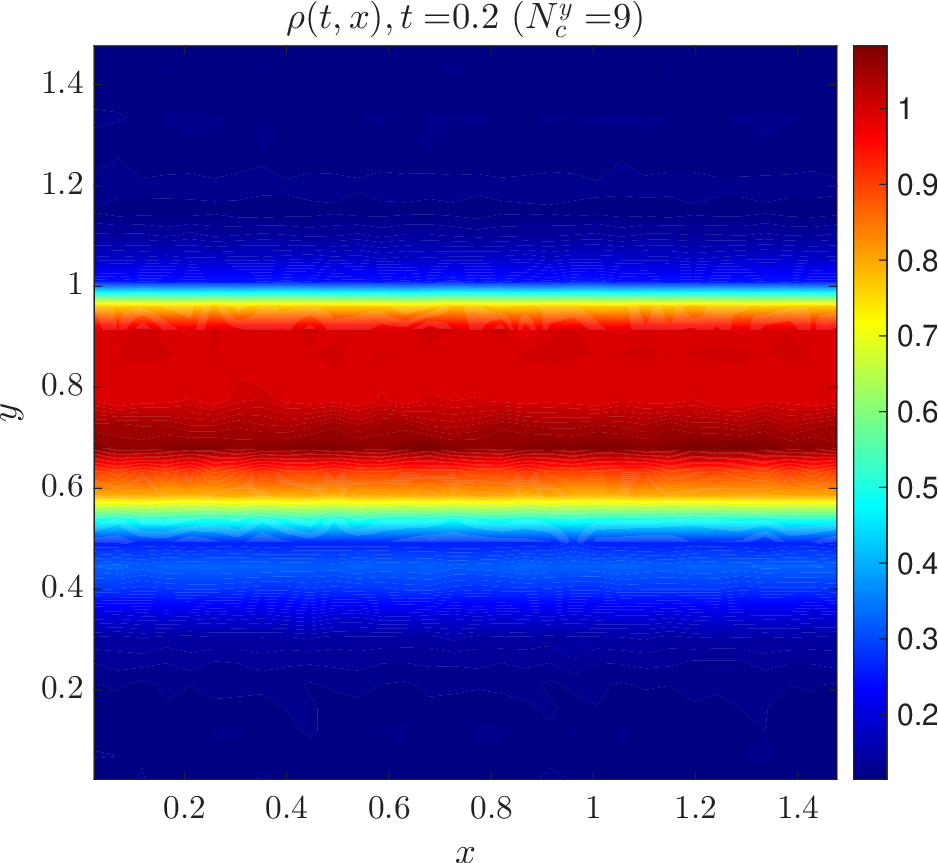}
	\includegraphics[width=0.327\linewidth]
	{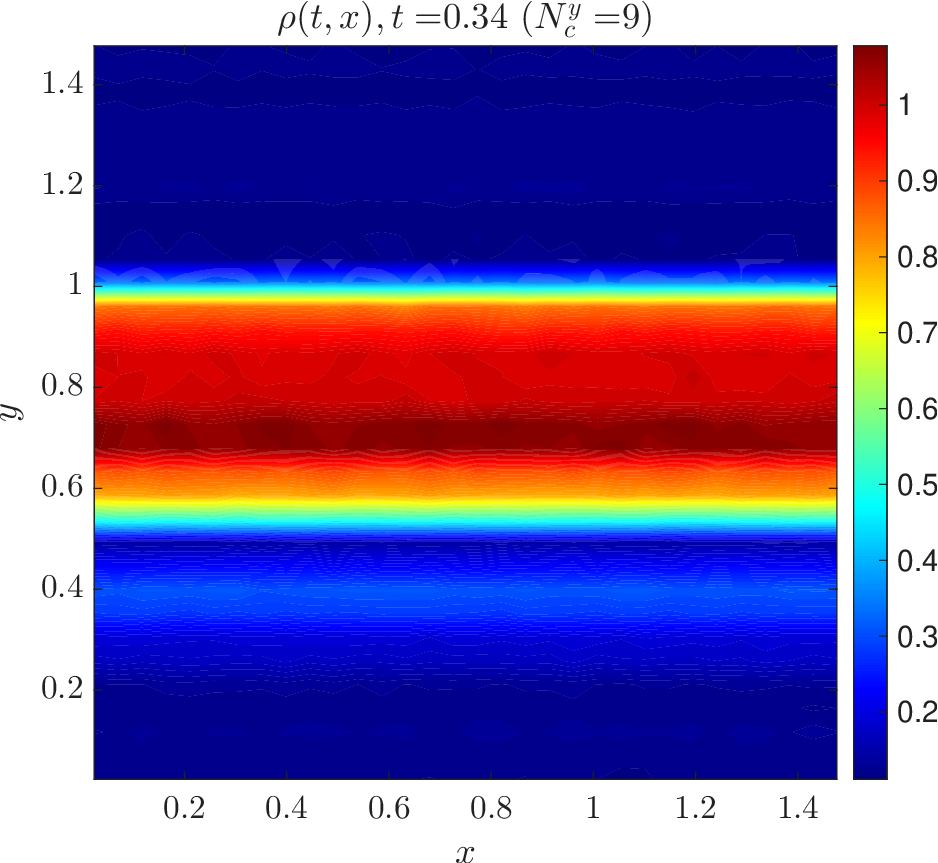}
	\includegraphics[width=0.327\linewidth]
	{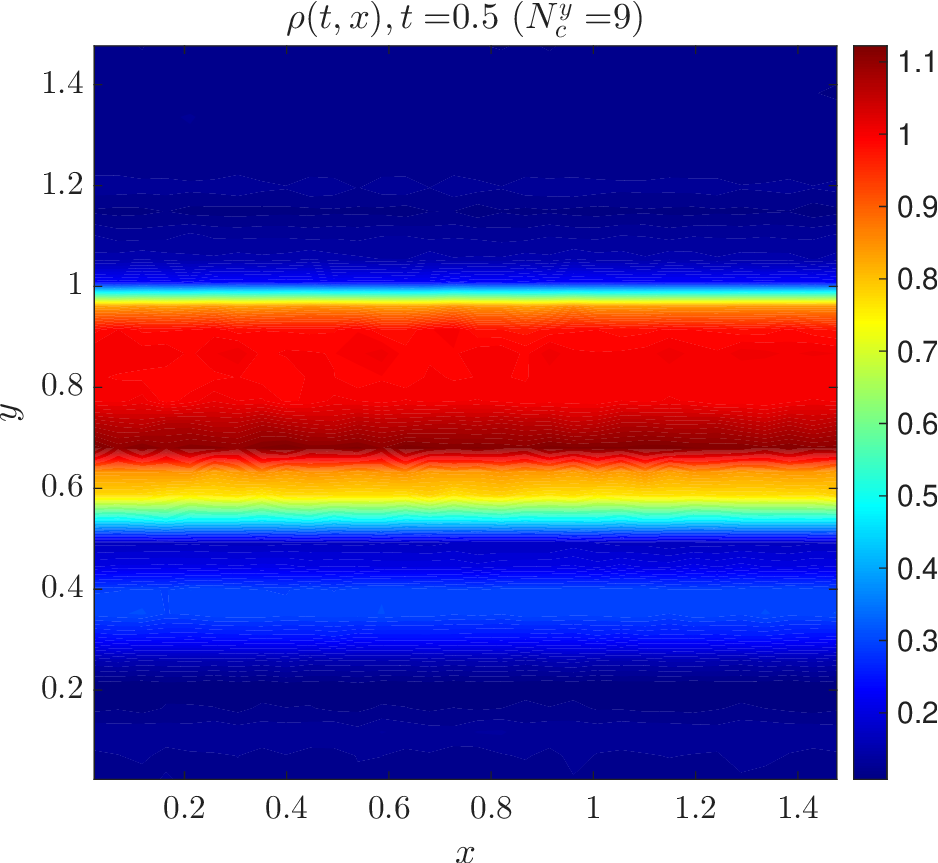}
	\
	\vspace*{5pt}
	\includegraphics[width=0.327\linewidth]
	{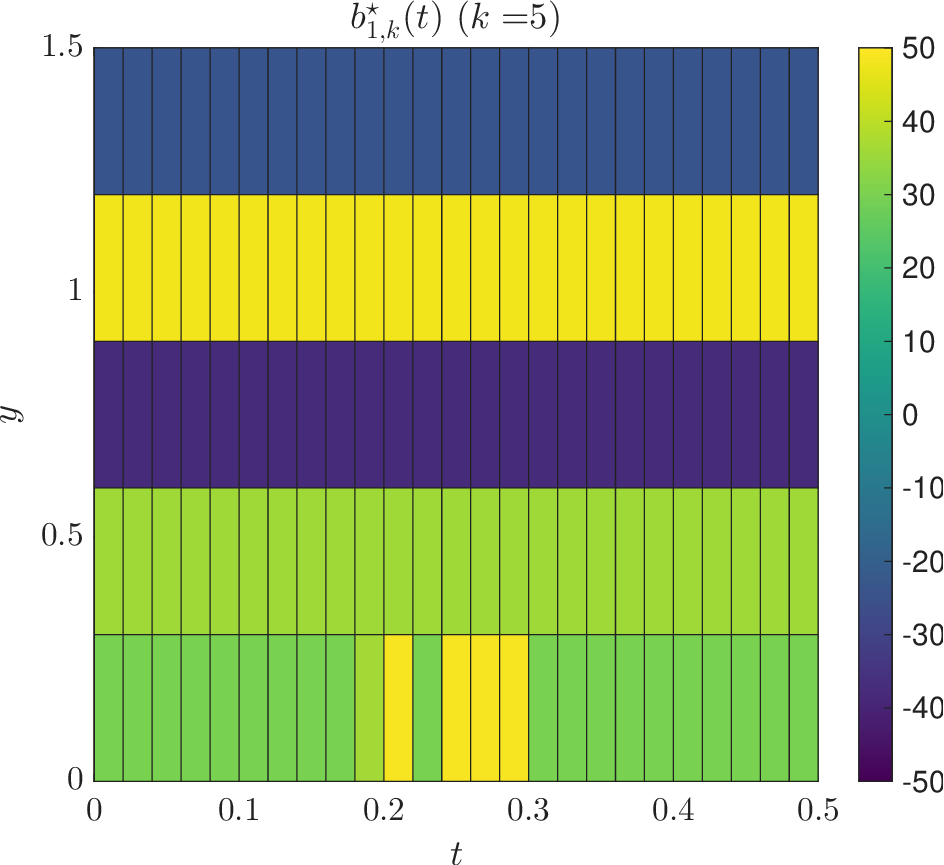}
	\includegraphics[width=0.327\linewidth]
	{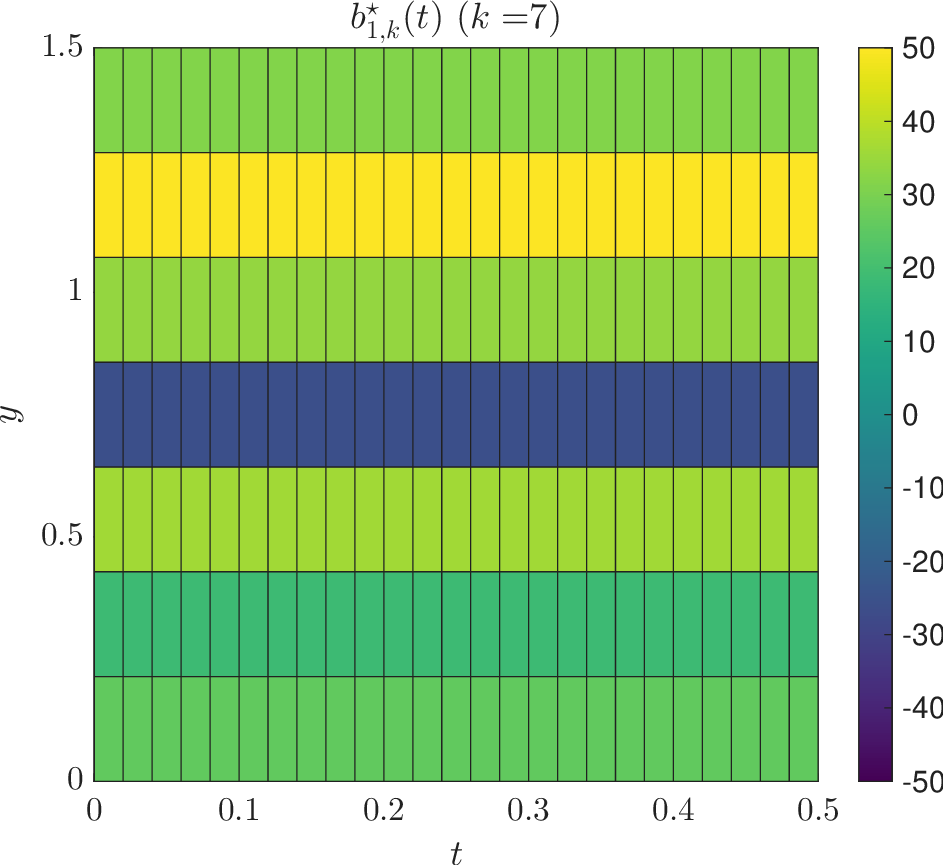}
	\includegraphics[width=0.327\linewidth]
	{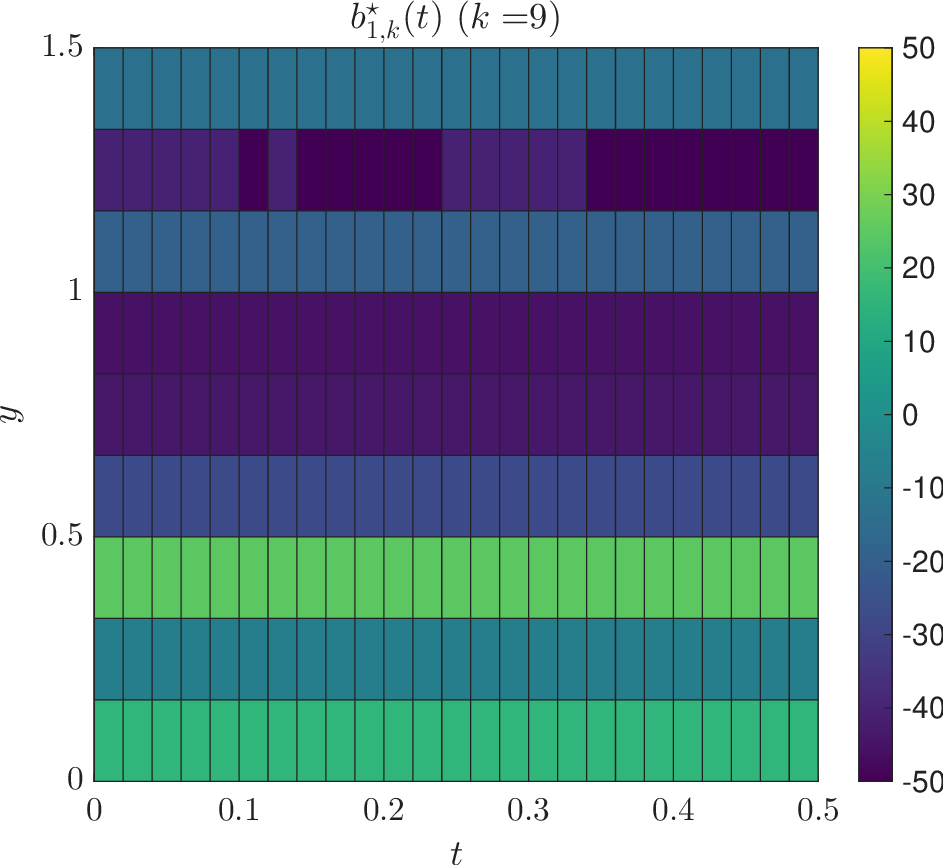}
	\caption{Two-dimensional double Sod test. First row: snapshots of
		the controlled dynamics at $t=0.2$, $t=0.34$, and $t=0.5$ for
		$N_c^y=9$. Second row: optimal control coefficients
		$b_{1,k}^\star(t)$.
		The	columns in the second row correspond, from left to right, to
		$N_c^y=5$, $N_c^y=7$, and $N_c^y=9$.}
	\label{fig:dynamics_B}
\end{figure}

To describe the convergence of the optimization algorithm, we
introduce the reduced cost functional
$$
\widehat{\mathcal J}(\mathbf b)
:=
\mathcal J(f[\mathbf b],\mathbf b),
$$
where $f[\mathbf b]$ denotes the solution of the forward problem
associated with the control $\mathbf b$. Figure~\ref{fig:functional}
reports the value of $\widehat{\mathcal J}(\mathbf b^{(m)})$ and the
Euclidean norm of its discrete gradient as functions of the
optimization iteration, for the different values of $N_c^y$.
The results indicate that the instantaneous control provides an
effective initial guess and is already close to a stationary point of
the discrete optimization problem: the initial gradient norm is
relatively small and the variation of the cost functional over the iterations is limited.
Nevertheless, the optimization produces a non-negligible correction
to the instantaneous control and, in most cases, further reduces the
fraction of thermal energy localized near the boundaries.
This comparison is presented in Figure~\ref{fig:energy}, which shows
$\mathcal E_b^{\mathrm{rel}}(t)$ for the instantaneous and optimal
controls. The three panels correspond to
$N_c^y=5$, $N_c^y=7$, and $N_c^y=9$, respectively. Overall, the
optimal control reduces the concentration of thermal energy in the
boundary strip compared with the instantaneous strategy.
The first row of Figure~\ref{fig:dynamics_B} shows instead three snapshots of
the controlled dynamics at $t=0.2$, $t=0.34$, and $t=0.5$ for
$N_c^y=9$. Similar qualitative behavior is observed for
$N_c^y=5$ and $N_c^y=7$ and they are not reported. The second row displays the optimal control
coefficients $b_{1,k}^\star(t)$, for
$N_c^y\in\{5,7,9\}$. 

\subsection{Two-stream plasma test}\label{subsec:two_stream}
\begin{figure}[!htb]
	\centering
	\includegraphics[width=0.35\linewidth]
	{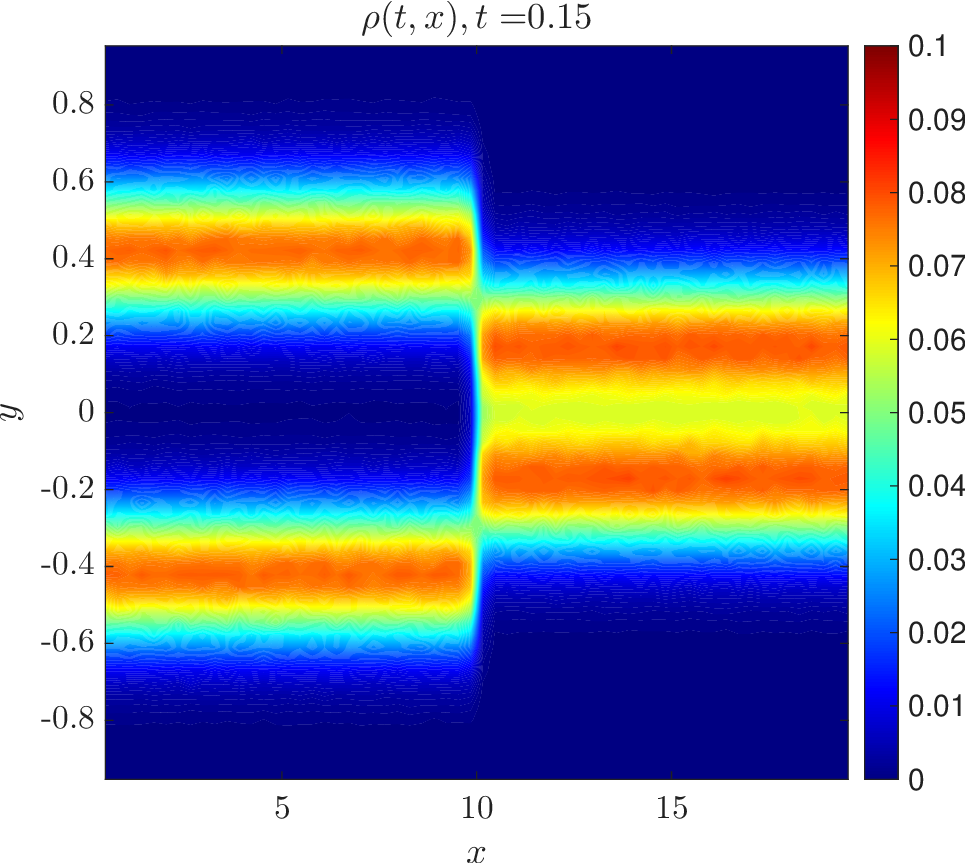}
	\includegraphics[width=0.35\linewidth]
	{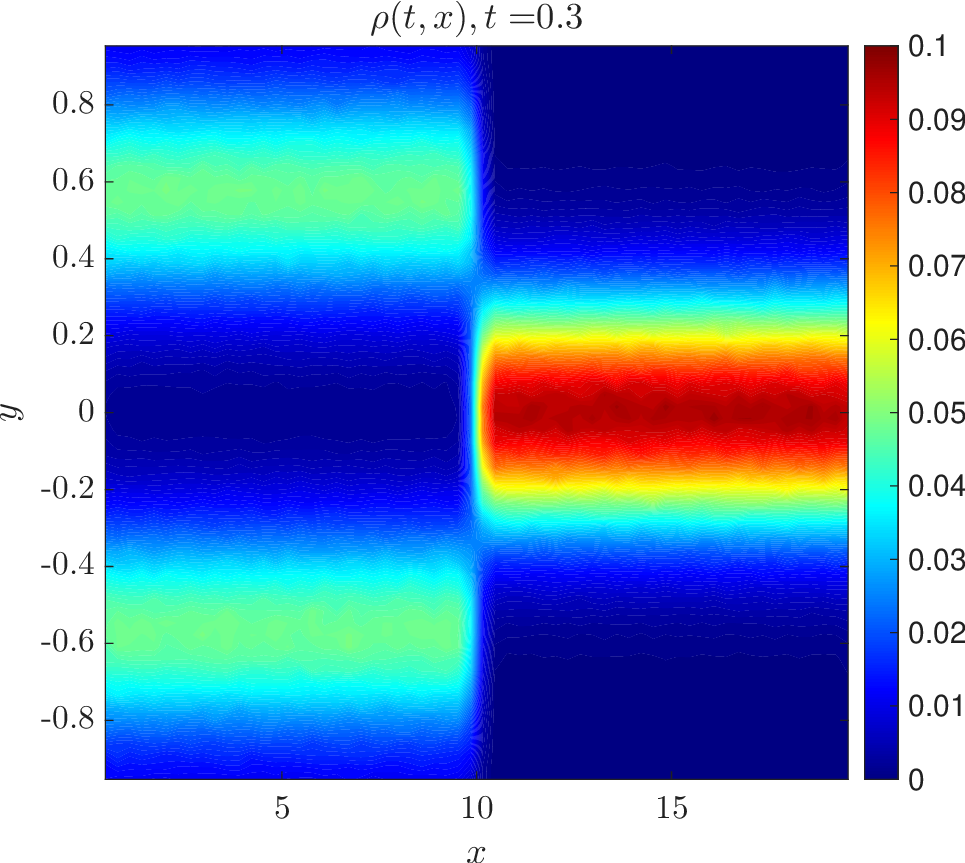}
	\includegraphics[width=0.35\linewidth]
	{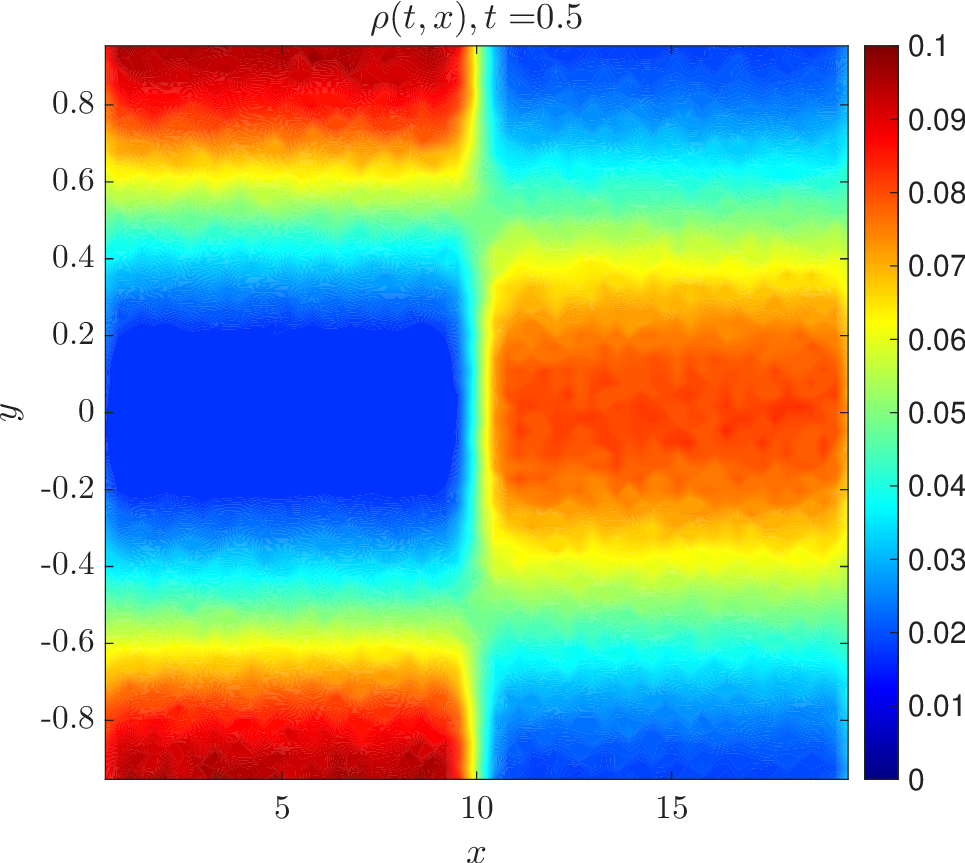}
	\includegraphics[width=0.33\linewidth]
	{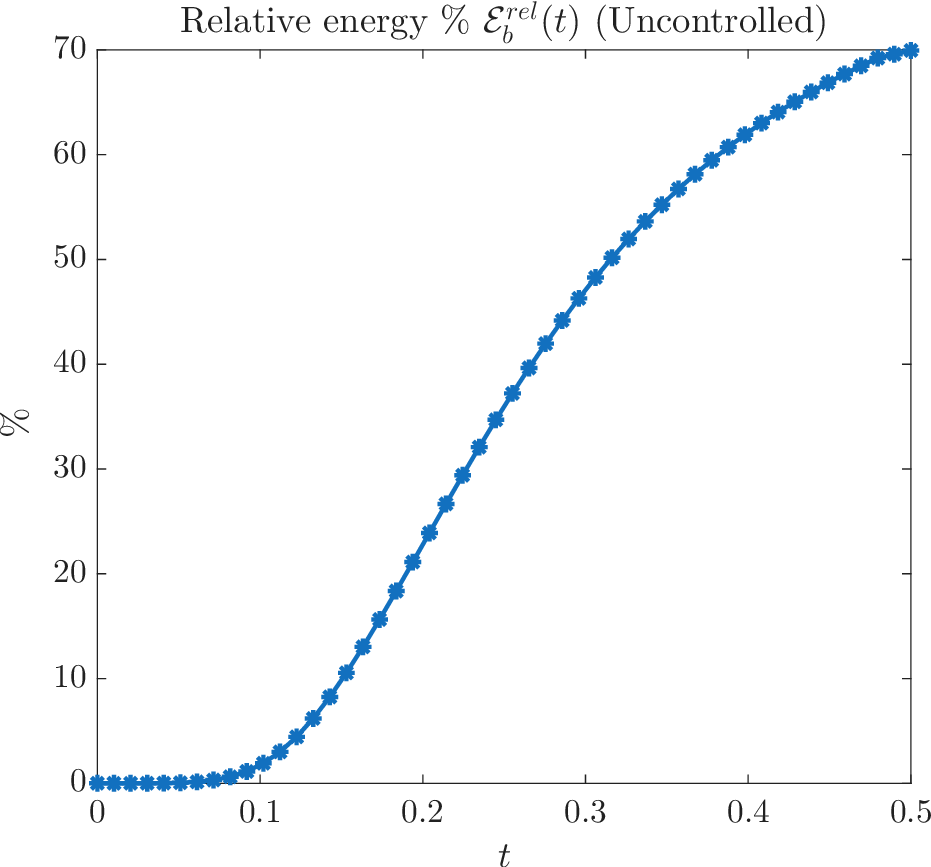}
	\vspace*{5pt}
	\caption{Two-stream plasma test with the prescribed uniform
		magnetic field $B_{\mathrm{ref}}=0.75$. Density
		snapshots at $t=0.15$, $t=0.30$, and $t=0.50$ and
		evolution of the percentage of thermal energy contained in the
		boundary strip, $\mathcal E_b^{\mathrm{rel}}(t)$.}
	\label{fig:two_stream_uncontrolled}
\end{figure}

We consider the spatial domain
\[
\Omega_{\xx}=[0,20]\times[-1,1],
\]
with periodic boundary conditions in the $x$-direction and reflective
boundary conditions in the $y$-direction. The initial distribution is
defined by
\begin{equation}\label{eq:ic_two_stream}
	f_0(\xx,\vv)
	=
	f_0^+(\xx,\vv)\chi_{\{y\geq0\}}
	+
	f_0^-(\xx,\vv)\chi_{\{y<0\}},
\end{equation}
where
\begin{equation}\label{eq:ic_two_stream_pm}
	f_0^\pm(\xx,\vv)
	=
	\frac{\rho_0^\pm(y)}{2\pi T_0}
	\exp\left(
	-\frac{v_x^2+\left(v_y-u_y(x)\right)^2}{2T_0}
	\right).
\end{equation}
The transverse mean velocity is prescribed as
\begin{equation}\label{eq:two_stream_velocity}
	u_y(x)
	=
	\begin{cases}
		-3, & x<10,\\
		3,  & x\geq10,
	\end{cases}
\end{equation}
while the two density profiles are given by
\begin{equation}\label{eq:two_stream_density}
	\rho_0^\pm(y)
	=
	\frac{1}{\sqrt{2\pi}\sigma}
	\exp\left(
	-\frac{(y\mp c_y)^2}{2\sigma^2}
	\right),
\end{equation}
with
\[
T_0=2,
\qquad
c_y=0.3,
\qquad
\sigma=0.1.
\]
With this convention, $\rho_0^+$ is centered at $y=c_y$ and
$\rho_0^-$ at $y=-c_y$. The regions $x<10$ and $x\geq10$ are assigned
opposite transverse velocities, thus generating counter-propagating
particle streams. The dynamics is evolved up to the final time $T=0.5$, using the time
step $\Delta t=0.01$. The density and the electric field are
reconstructed on a uniform PIC grid with $m_x=m_y=64$ nodes, so that
$N_h=m_xm_y$. For both the instantaneous and optimal controls, we set
\[
\alpha_{\mathrm{x}}= \alpha_{\mathrm{v}} = 10,
\qquad
\beta_{\mathrm{x}} =
\beta_{\mathrm{v}}=1,
\qquad
\gamma=10^{-3}, \qquad M=50,
\]
with targets $\hat{y} = \hat{v}_y =0$.
We first consider a reference configuration in which no optimization
is performed and the magnetic field is prescribed to be uniform in
space and constant in time:
$
B_{\mathrm{ref}}(t,\xx)=0.75.
$ Figure~\ref{fig:two_stream_uncontrolled} shows three
density snapshots at $t=0.15$, $t=0.30$, and $t=0.50$. Owing to the
opposite initial transverse velocities, particles are transported
toward the upper and lower boundaries and are subsequently reflected.
The lower panel on the right reports the corresponding percentage of thermal
energy contained in the boundary strip,
$\mathcal E_b^{\mathrm{rel}}(t)$.

\begin{figure}[!htb]
	\centering
	\includegraphics[width=0.35\linewidth]
	{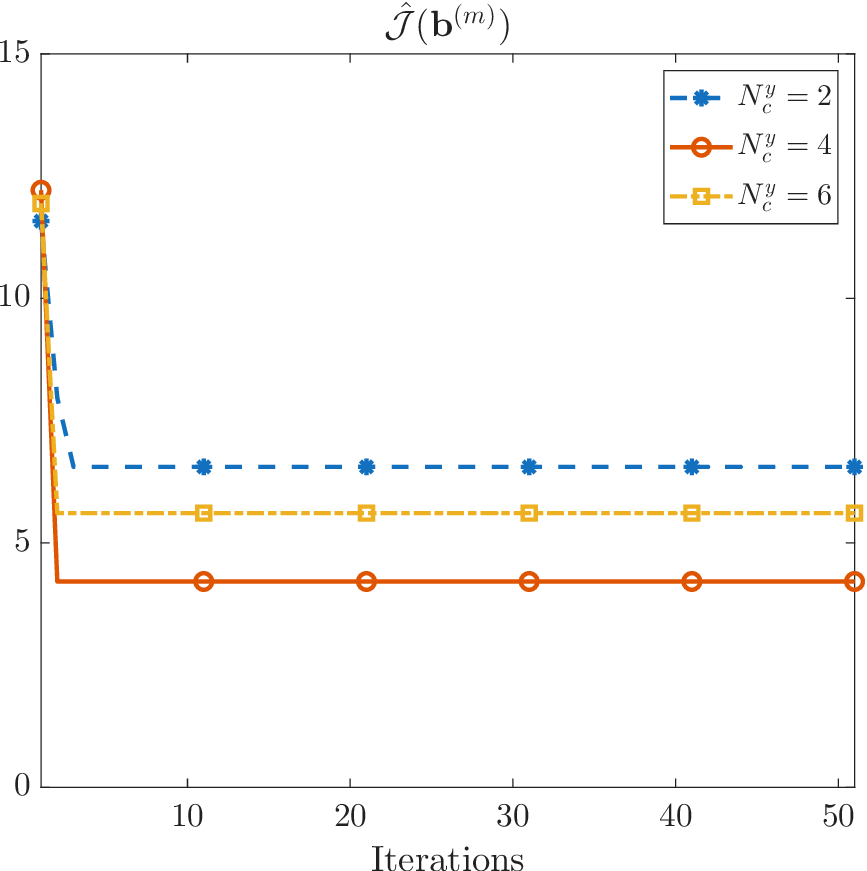}
	\includegraphics[width=0.355\linewidth]
	{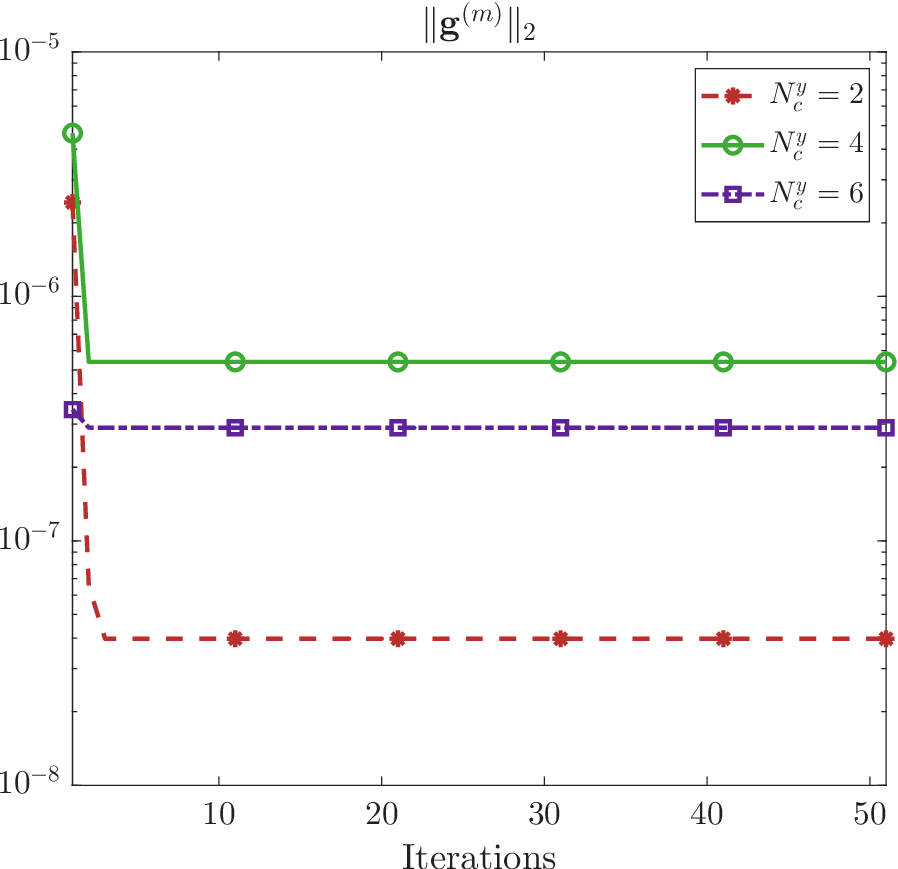}
	\caption{Two-stream plasma test. Left: evolution of the reduced
		cost functional $\widehat{\mathcal J}(\mathbf b^{(m)})$.
		Right: discrete gradient norm $\|\mathbf g^{(m)}\|_2$.
	The curves correspond to
		$N_c^x=1$ and $N_c^y\in\{2,4,6\}$.}
	\label{fig:functional_two_stream}
\end{figure}

\begin{figure}[!htb]
	\centering
	\includegraphics[width=0.327\linewidth]
	{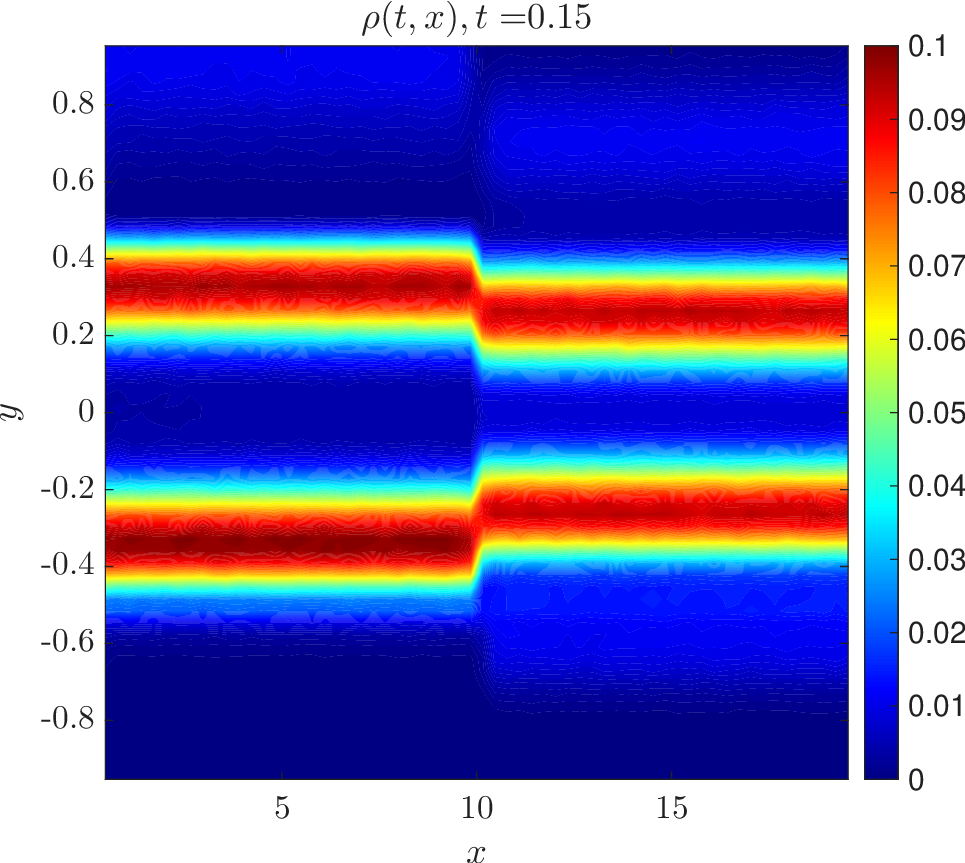}
	\includegraphics[width=0.327\linewidth]
	{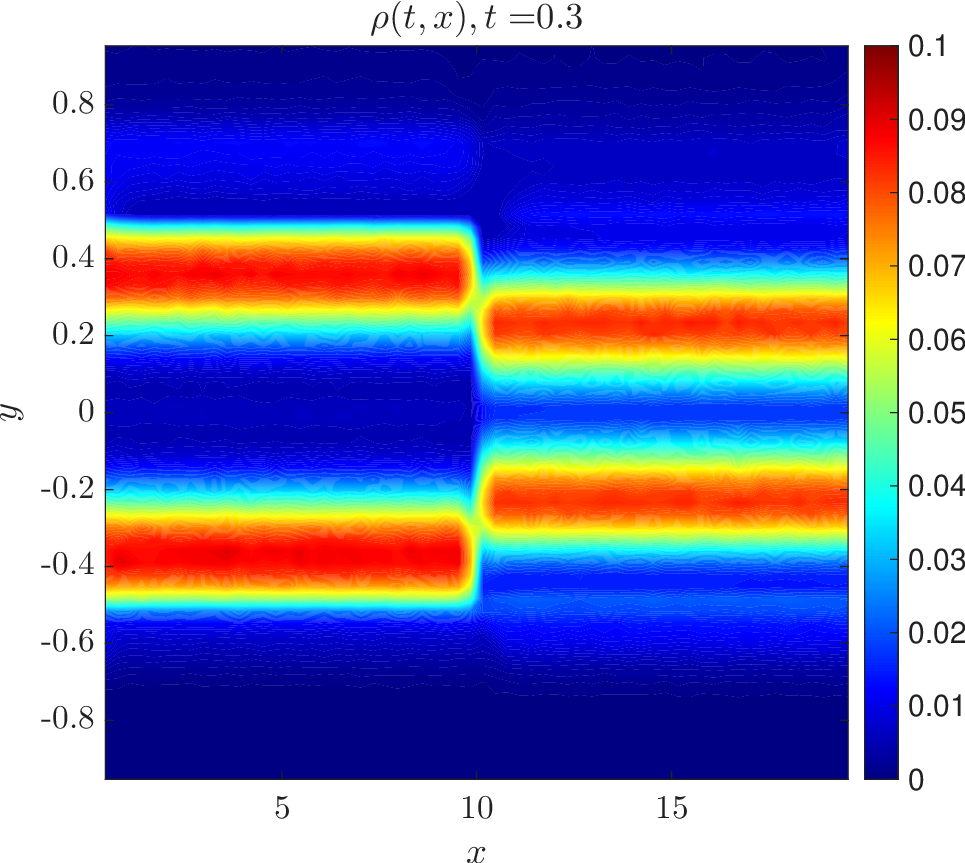}
	\includegraphics[width=0.327\linewidth]
	{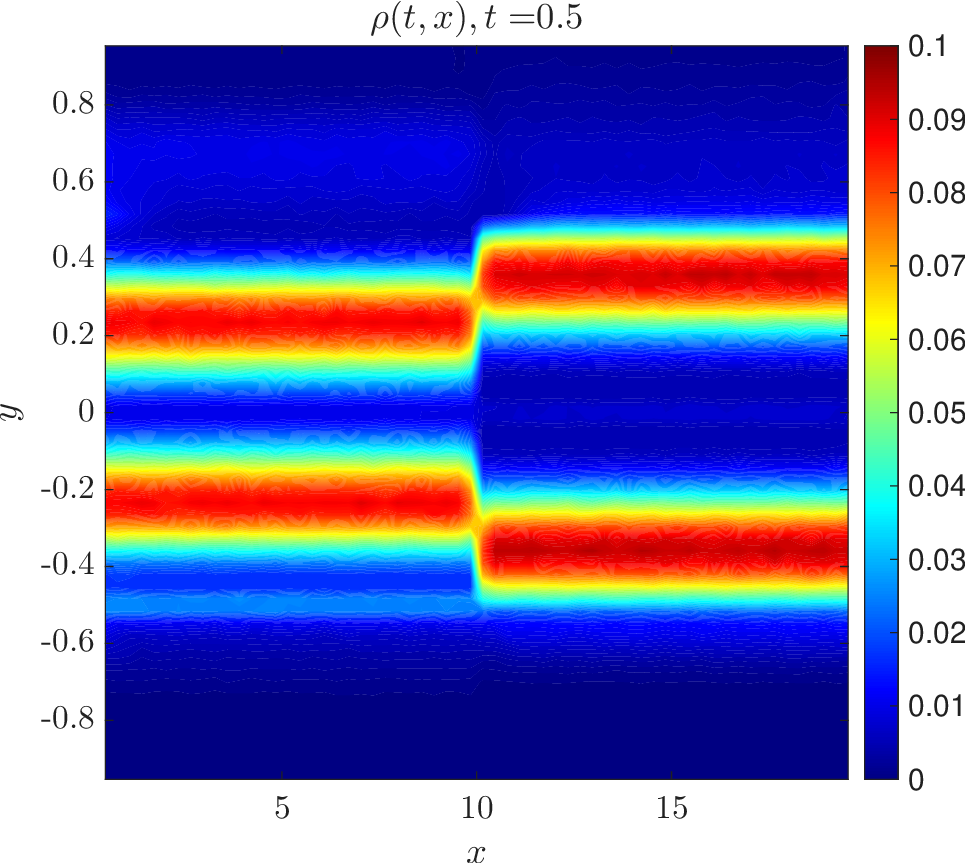}
		\includegraphics[width=0.327\linewidth]
	{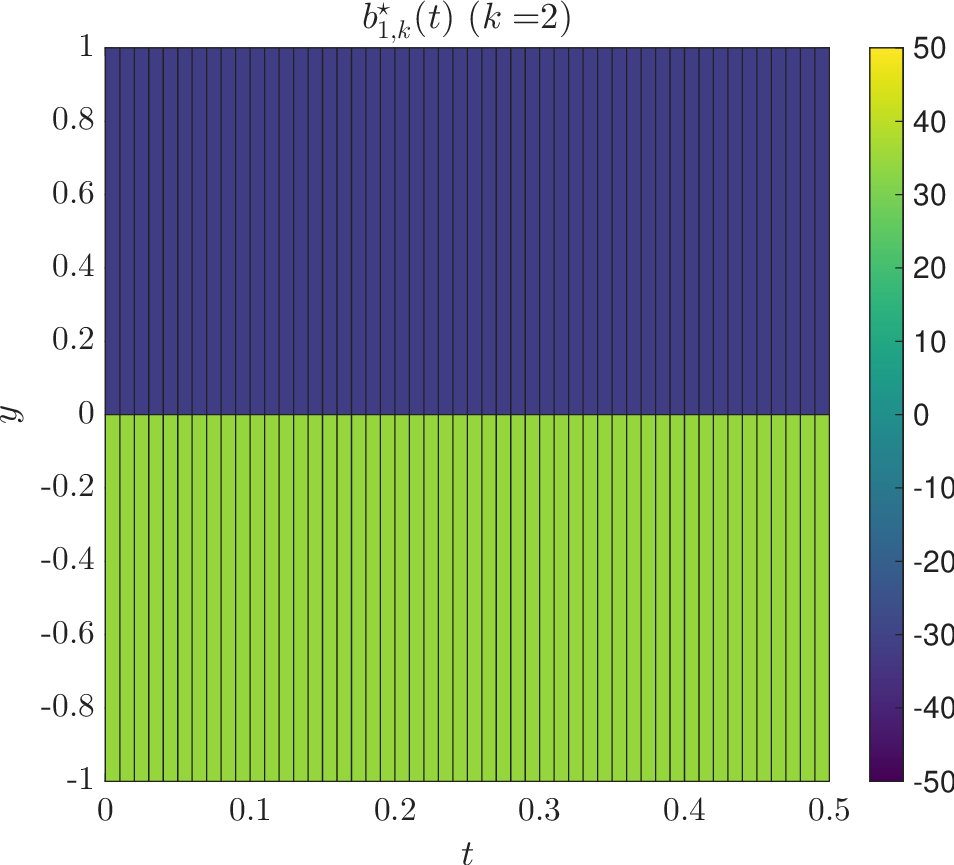}
	\includegraphics[width=0.327\linewidth]
	{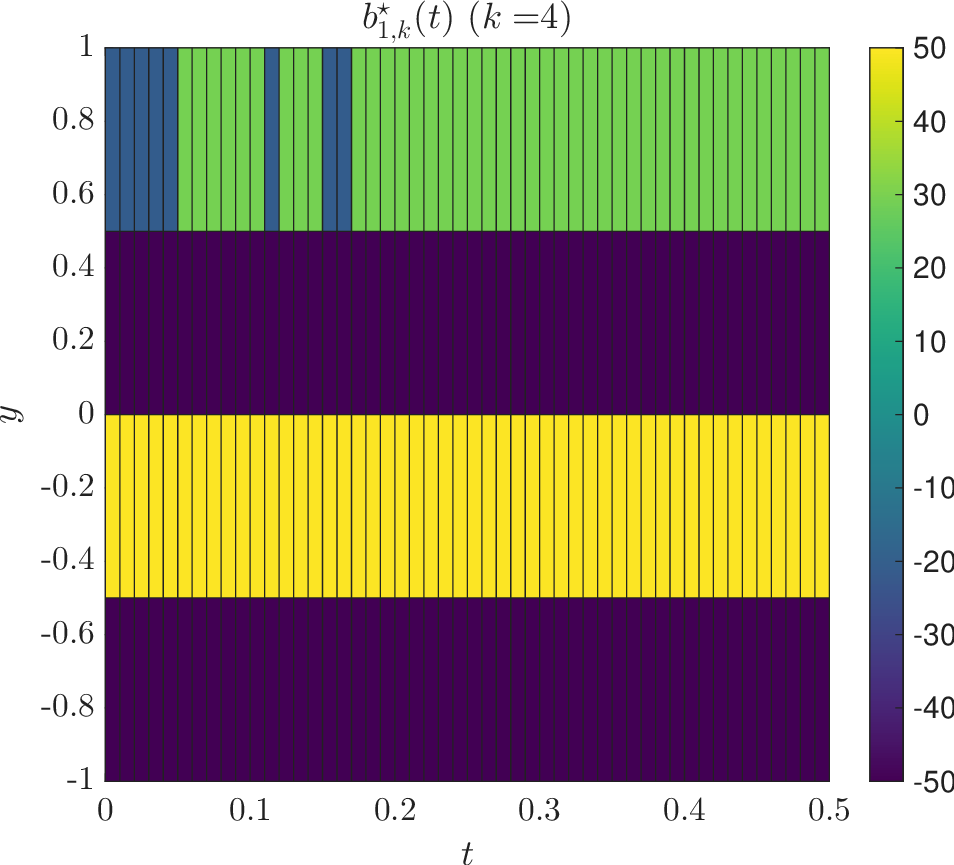}
	\includegraphics[width=0.327\linewidth]
	{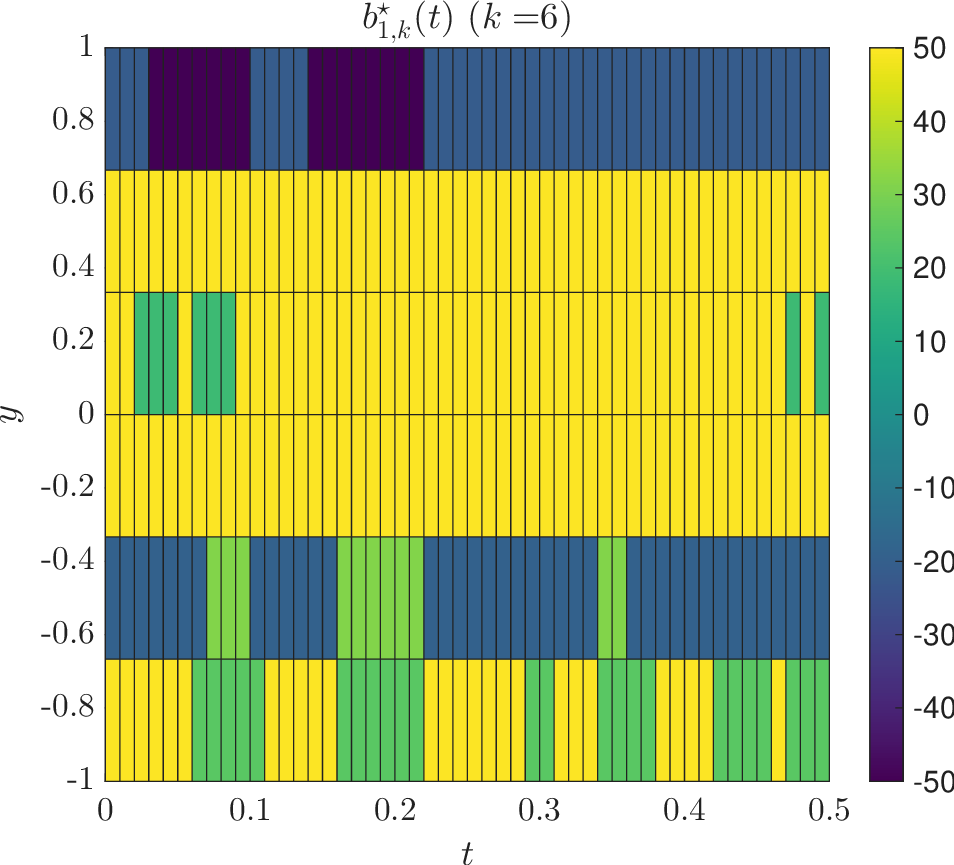}
	\vspace*{5pt}
	\caption{Two-stream plasma test. First row: density snapshots of
		the controlled dynamics at $t=0.15$, $t=0.30$, and $t=0.50$ for
		$N_c^y=4$. Second row: optimal control
		coefficients $b_{1,k}^\star(t)$. The columns in the second row correspond, from left to right, to
		$N_c^y=2$, $N_c^y=4$, and $N_c^y=6$.}
	\label{fig:two_stream_dynamics_control_Nc4}
\end{figure}

We next compute the controlled dynamics. We set $N_c^x=1$ and consider
$N_c^y\in\{2,4,6\}$. Hence, the control cells extend over the entire
$x$-direction and the reconstructed magnetic field depends only on
$y$ and $t$. Figure~\ref{fig:functional_two_stream} shows the convergence history
of the optimization procedure. The two panels report the reduced
cost functional $\widehat{\mathcal J}(\mathbf b^{(m)})$ and the
Euclidean norm of its discrete gradient
$\|\mathbf g^{(m)}\|_2$.
As in the previous test,  the
iterative procedure produces a
decrease of the reduced cost functional.
Figure~\ref{fig:two_stream_dynamics_control_Nc4} shows three density
snapshots of the controlled dynamics at $t=0.15$, $t=0.30$, and
$t=0.50$ for $N_c^y=4$. The control substantially limits transport
toward the upper and lower boundaries and retains a larger fraction
of the density in the central region. Similar qualitative behavior
is observed for $N_c^y=2$ and $N_c^y=6$. The second row displays the optimal control
coefficients over time.

\begin{figure}[!htb]
	\centering
	\includegraphics[width=0.327\linewidth]
	{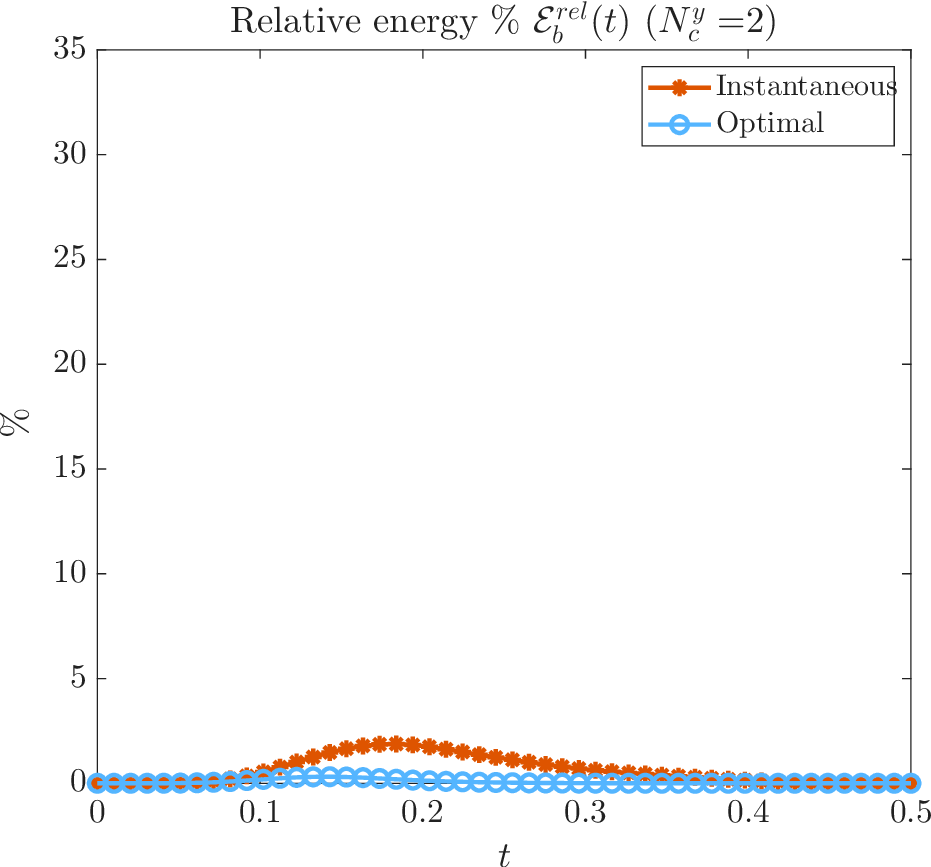}
	\includegraphics[width=0.327\linewidth]
	{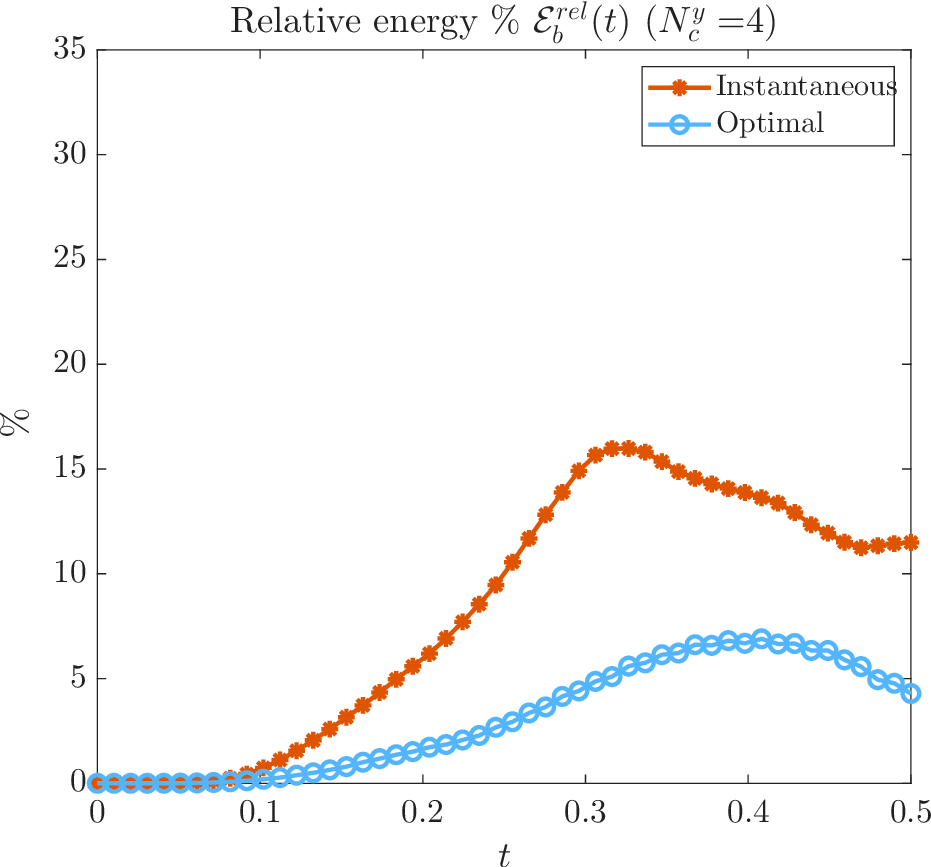}
	\includegraphics[width=0.327\linewidth]
	{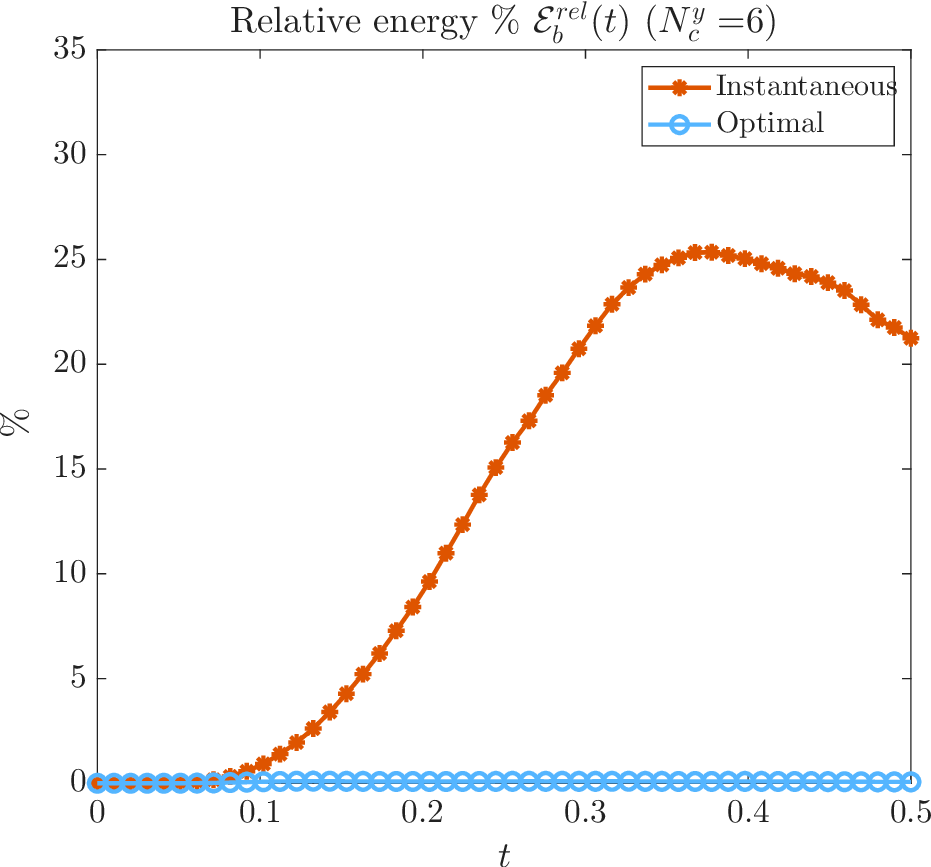}
	\vspace*{5pt}
	\caption{Two-stream plasma test. Comparison of
		$\mathcal E_b^{\mathrm{rel}}(t)$ obtained with the instantaneous
		and optimal controls. From left to right:
		$N_c^y=2$, $N_c^y=4$, and $N_c^y=6$.}
	\label{fig:two_stream_energy_comparison}
\end{figure}

Finally, Figure~\ref{fig:two_stream_energy_comparison} provides a
direct comparison between the instantaneous and optimal controls.
The optimized
strategy yields a further reduction in the fraction of thermal energy
localized near the boundaries.

\subsection{Kelvin--Helmholtz instability}\label{subsec:kelvin_helmholtz}

\begin{figure}[!htb]
	\centering
	\includegraphics[width=0.327\linewidth]
	{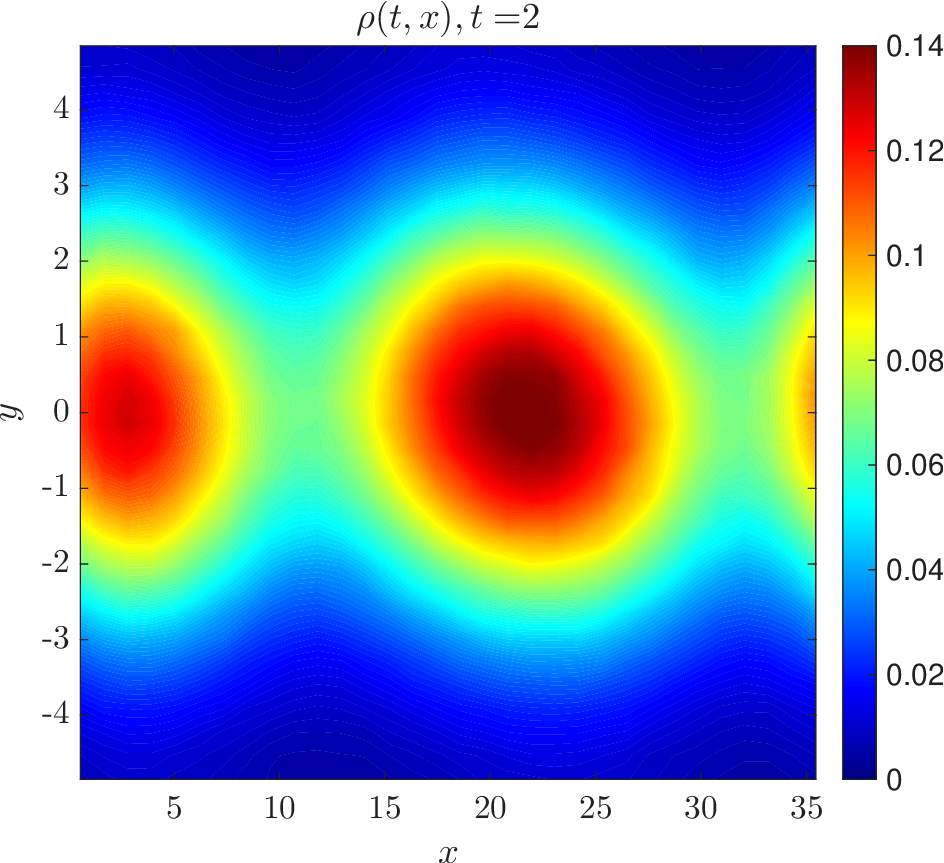}
	\includegraphics[width=0.327\linewidth]
	{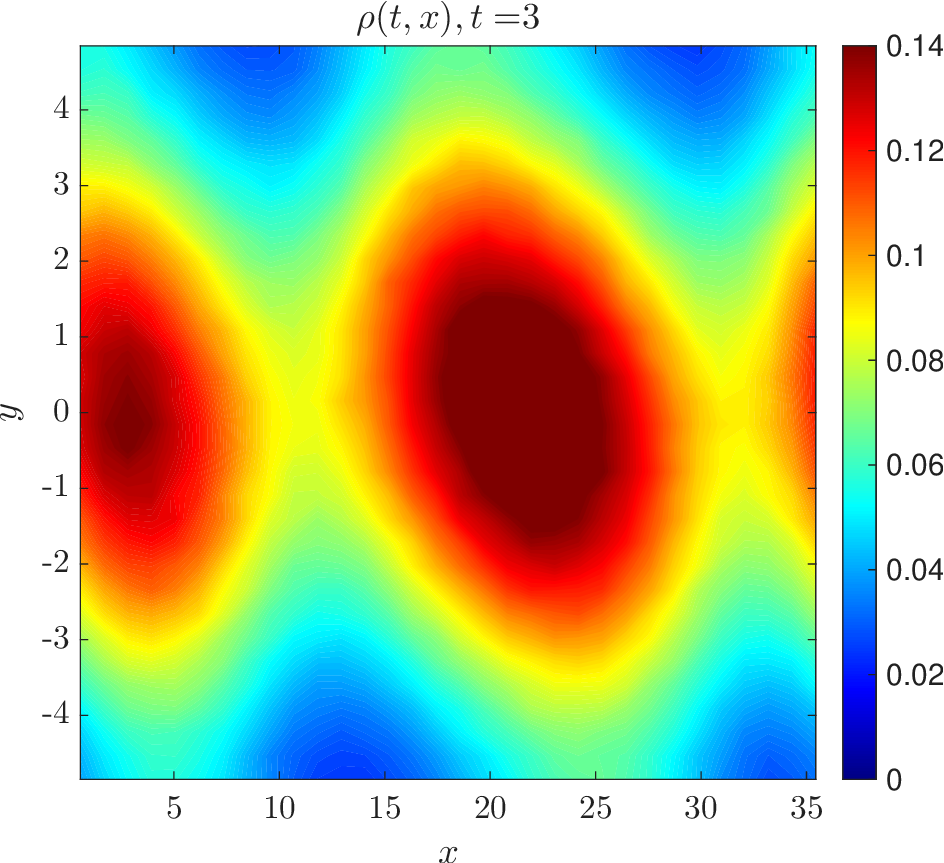}\\
	\includegraphics[width=0.327\linewidth]
	{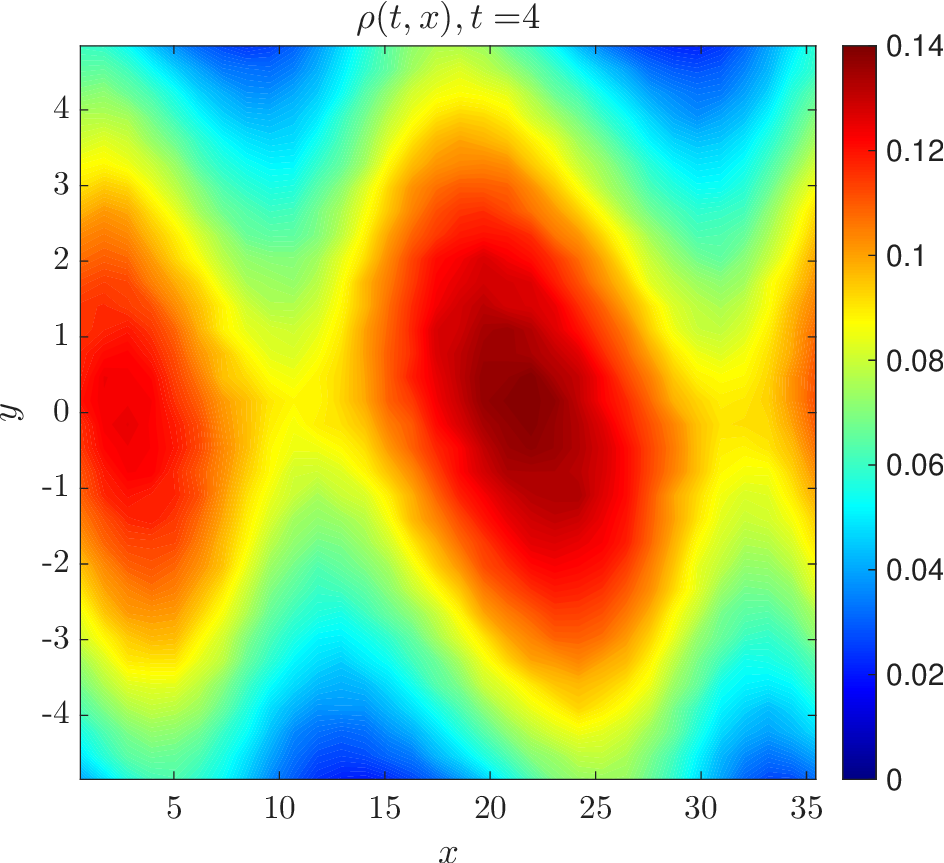}
	\includegraphics[width=0.327\linewidth]
	{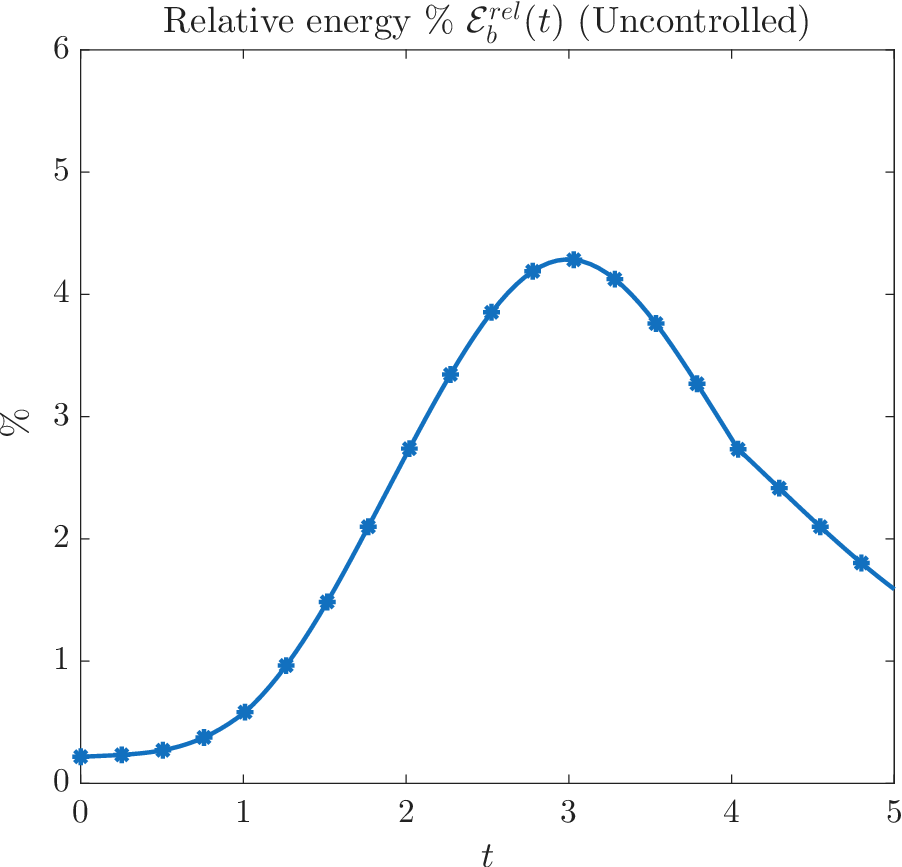}
	\vspace*{5pt}
	\caption{Kelvin--Helmholtz instability with the prescribed uniform
		magnetic field $B_{\mathrm{ref}}=0.75$. Density
		snapshots at $t=2$, $t=3$, and $t=4$ and evolution
		of the percentage of thermal energy contained in the boundary
		strip, $\mathcal E_b^{\mathrm{rel}}(t)$.}
	\label{fig:kelvin_dynamics_nocontrol}
\end{figure}

We next consider the development and control of a
Kelvin--Helmholtz instability. The spatial domain is
\[
\Omega_{\xx}=[0,36]\times[-5,5],
\]
with periodic boundary conditions in the $x$-direction and reflective
boundary conditions in the $y$-direction. The initial spatial density
is prescribed as
\begin{equation}\label{eq:kelvin_initial_distr}
	\rho_0(\xx)
	=
	\frac{1.5}{2\pi}
	\operatorname{sech}\left(\frac{y}{0.9}\right)
	\left[
	1+\varepsilon_0\cos(3k_0x)
	+\varepsilon_1\sin(k_0x)
	\right],
\end{equation}
where
\[
k_0=0.3,
\qquad
\varepsilon_0=0.5,
\qquad
\varepsilon_1=0.1.
\]
The initial phase-space distribution is
\begin{equation}\label{eq:kelvin_f0}
	f_0(\xx,\vv)
	=
	f_0^+(\xx,\vv)\chi_{\{y\geq0\}}
	+
	f_0^-(\xx,\vv)\chi_{\{y<0\}},
\end{equation}
with
\begin{equation}\label{eq:kelvin_f0_pm}
	f_0^\pm(\xx,\vv)
	=
	\frac{\rho_0(\xx)}{2\pi T_0}
	\exp\left(
	-\frac{(v_x\pm u_x)^2+v_y^2}{2T_0}
	\right),
\end{equation}
where
\[
T_0=1,
\qquad
u_x=1.
\]
The dynamics is evolved up to the final time $T=5$, using
$\Delta t=0.05$. The density and the electric field are reconstructed
on a uniform PIC grid with $m_x=m_y=32$ nodes, so that
$N_h=m_xm_y$. For both control strategies, we set
\[
\alpha_{\mathrm{x}}
=
\alpha_{\mathrm{v}}
=
10,
\qquad
\beta_{\mathrm{x}}
=
\beta_{\mathrm{v}}
=
1, \qquad \gamma=10^{-3}, \qquad M=50,
\]  with targets $\hat{y} = \hat{v}_y =0$.
We first consider the reference dynamics generated by the uniform
and time-independent magnetic field
$
B_{\mathrm{ref}}(t,\xx)=0.75.
$
Figure~\ref{fig:kelvin_dynamics_nocontrol} shows
density snapshots at $t=2$, $t=3$, and $t=4$. The deformation of
the initially localized layer illustrates the development of the
instability. The lower panel on the right reports the corresponding percentage of
thermal energy localized within the boundary strip,
$\mathcal E_b^{\mathrm{rel}}(t)$.

\begin{figure}[!htb]
	\centering
	\includegraphics[width=0.35\linewidth]
	{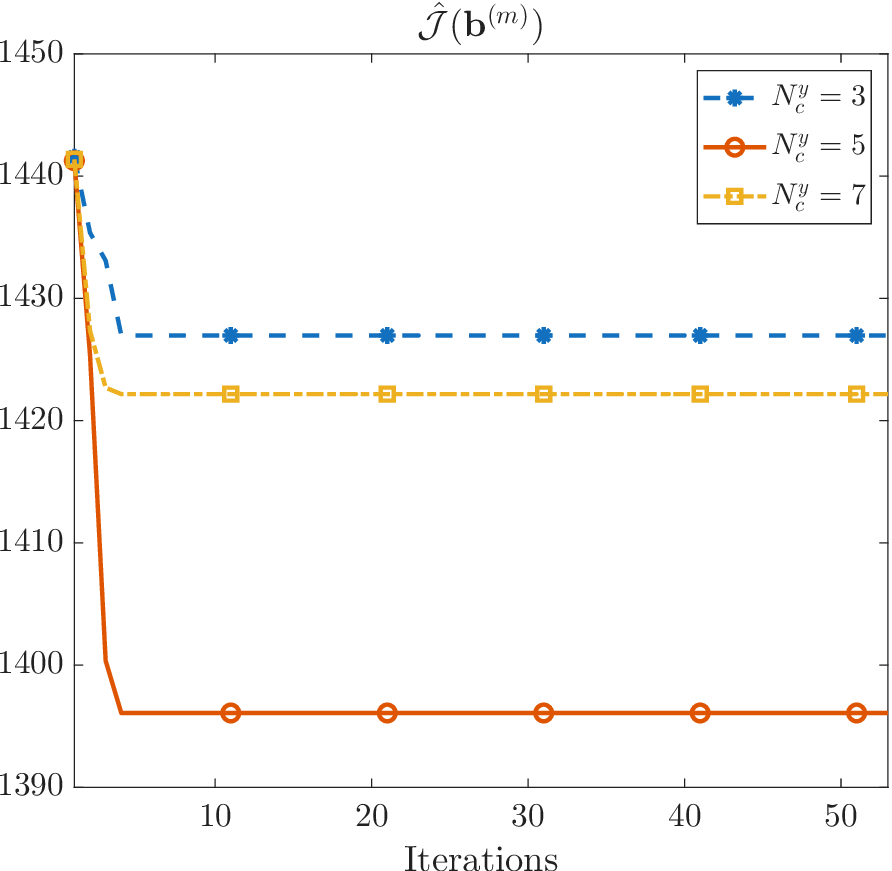}
	\includegraphics[width=0.35\linewidth]
	{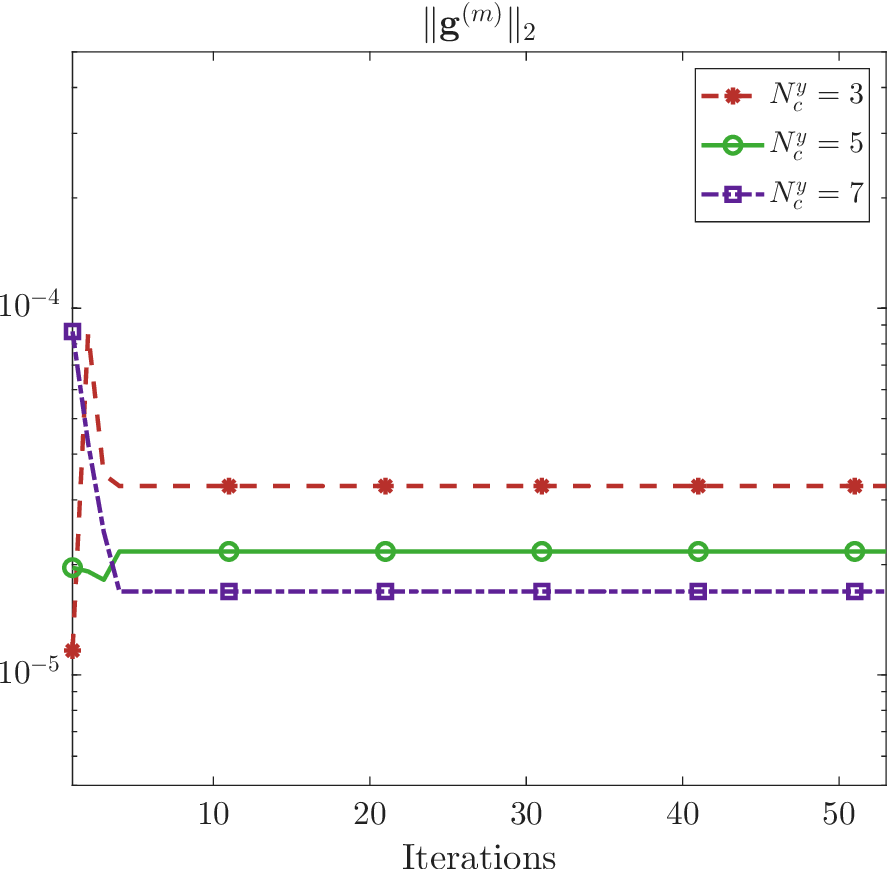}
	\caption{Kelvin--Helmholtz instability. Left: evolution of the
		reduced cost functional
		$\widehat{\mathcal J}(\mathbf b^{(m)})$. Right: discrete
		gradient norm $\|\mathbf g^{(m)}\|_2$. The
		curves correspond to $N_c^x=1$ and
		$N_c^y\in\{3,5,7\}$.}
	\label{fig:functional_kelvin}
\end{figure}

\begin{figure}[!htb]
	\centering
	\includegraphics[width=0.327\linewidth]
	{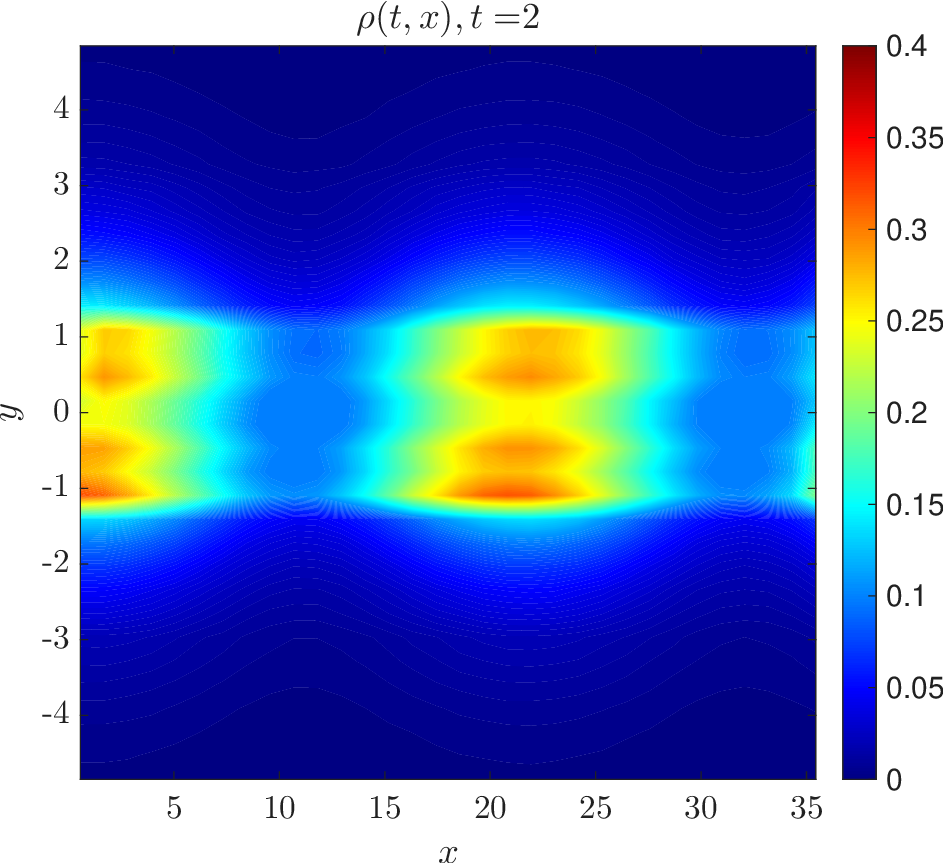}
	\includegraphics[width=0.327\linewidth]
	{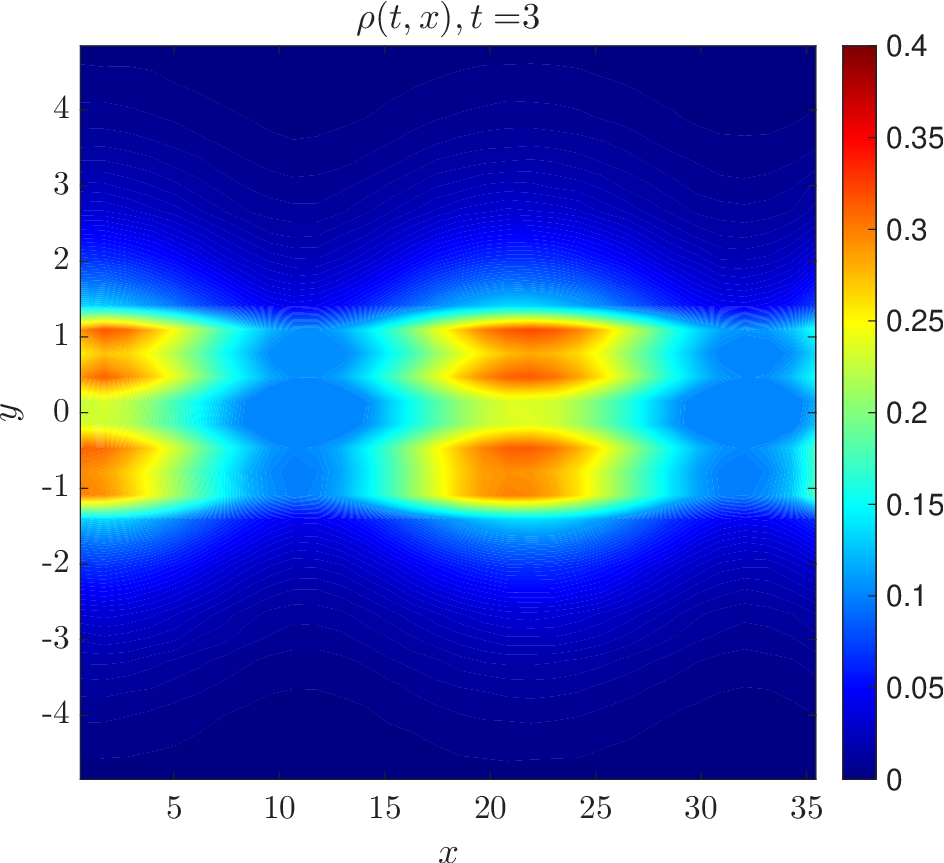}
	\includegraphics[width=0.327\linewidth]
	{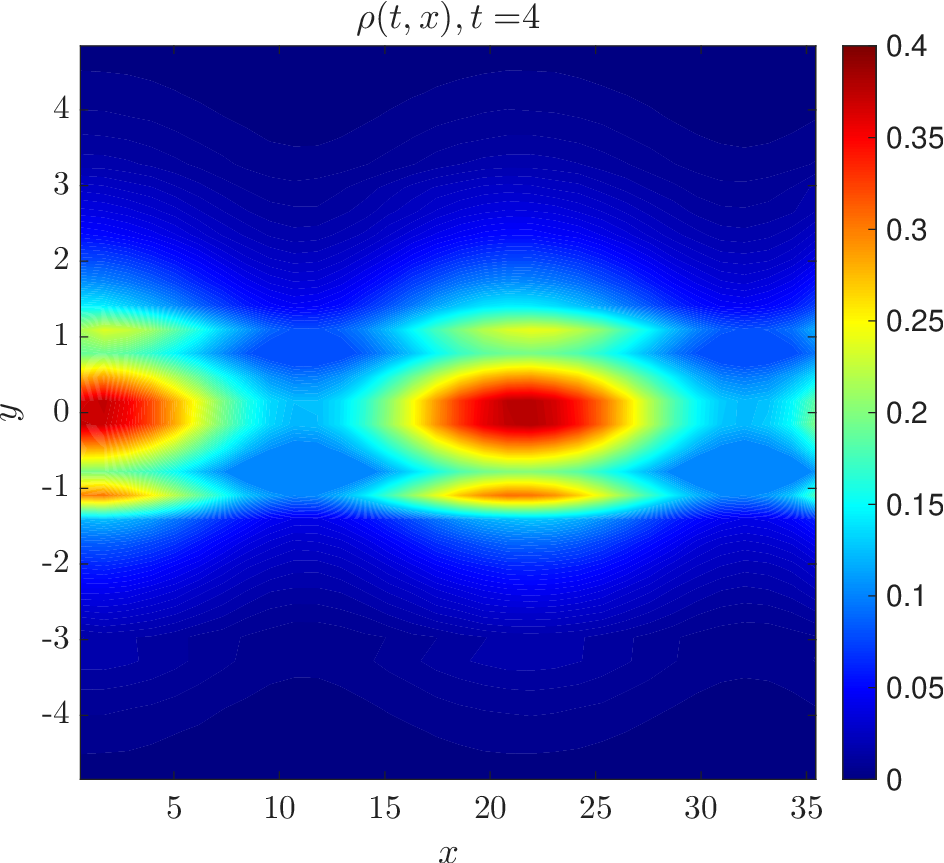}\\
	\includegraphics[width=0.327\linewidth]
	{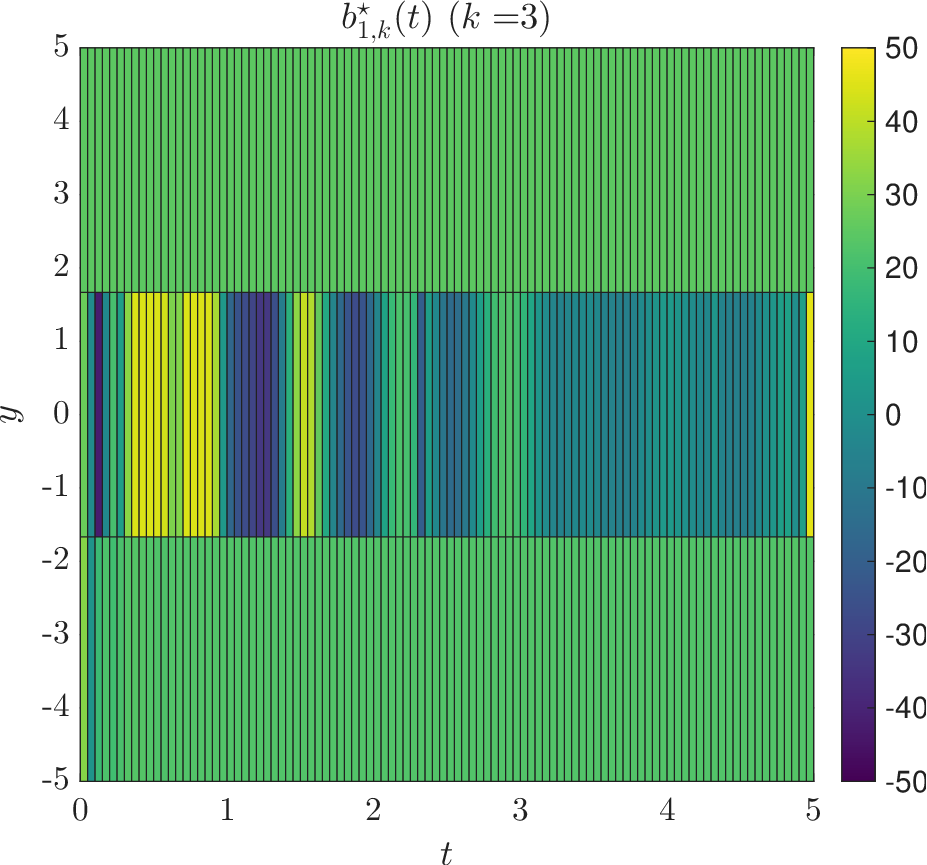}
	\includegraphics[width=0.327\linewidth]
	{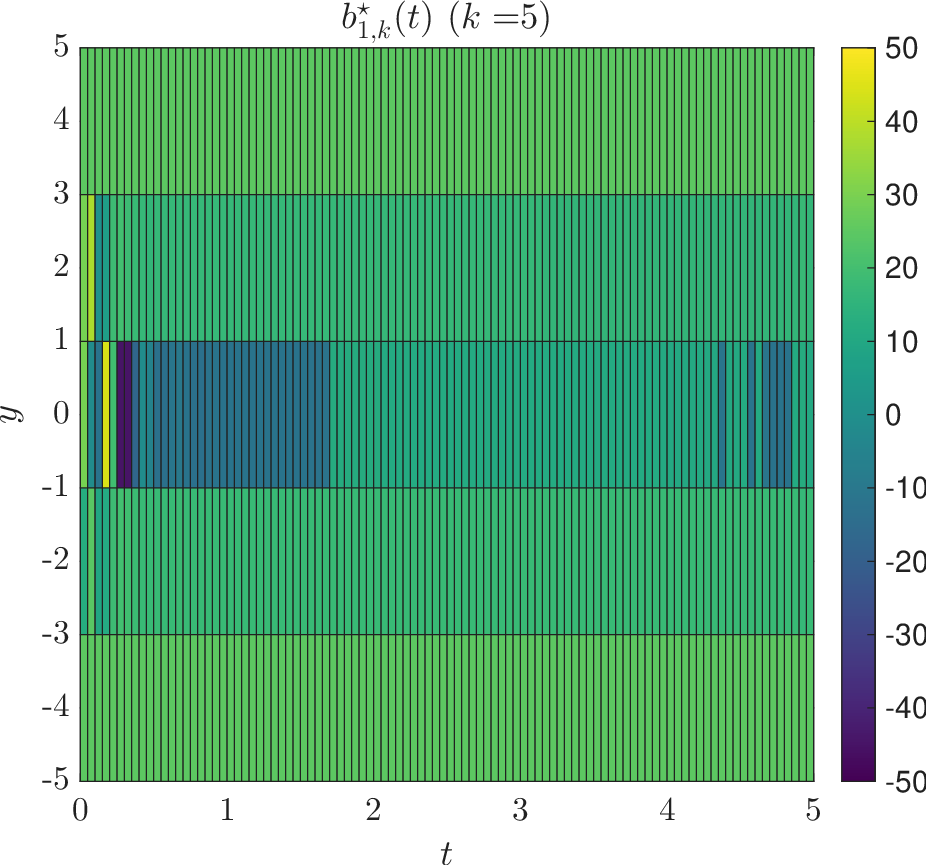}
	\includegraphics[width=0.327\linewidth]
	{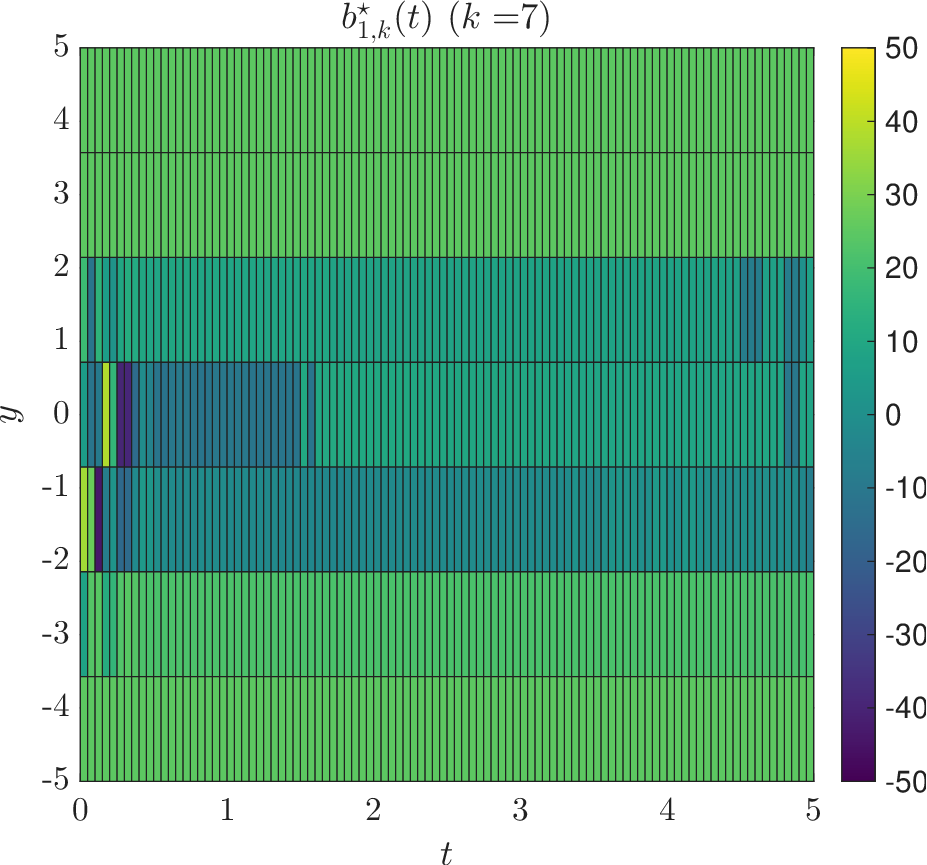}
	\caption{Kelvin--Helmholtz instability. Top row: density
		snapshots of the controlled dynamics at $t=2$, $t=3$, and
		$t=4$ for $N_c^y=5$. Bottom row: optimal control
		coefficients $b_{1,k}^\star(t)$. The columns in the second row correspond, from left to right, to
		$N_c^y=3$, $N_c^y=5$, and $N_c^y=7$.}
	\label{fig:kelvin_dynamics_control_Nc5}
\end{figure}

\begin{figure}[!htb]
	\centering
	\includegraphics[width=0.327\linewidth]
	{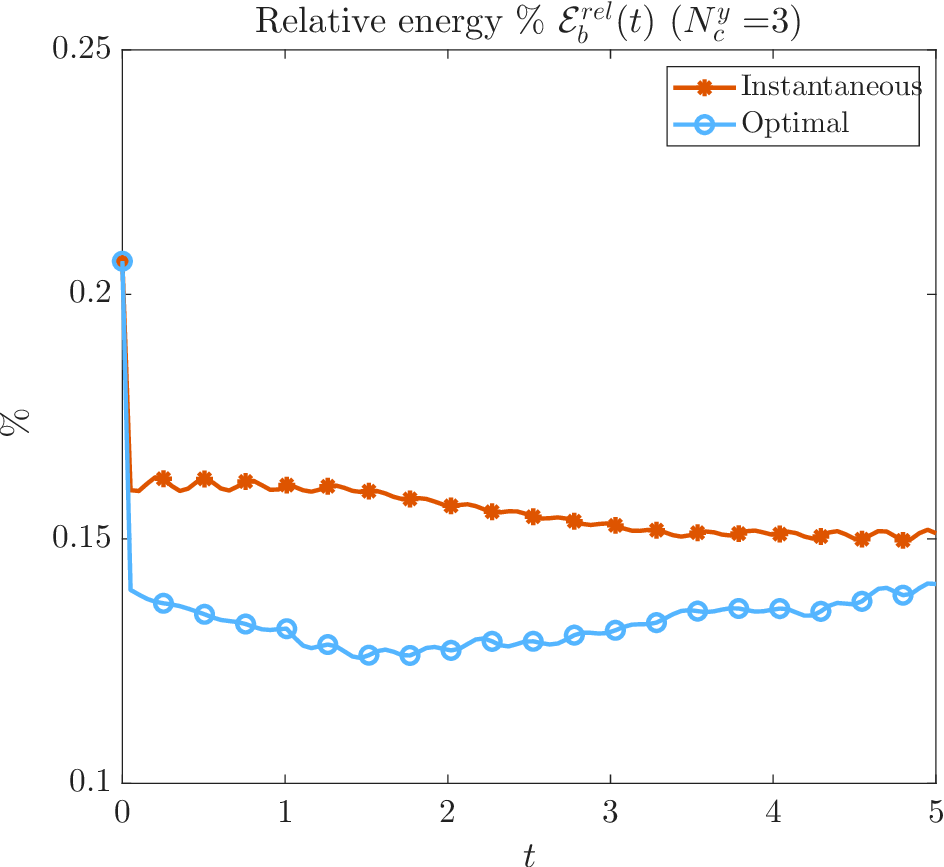}
	\includegraphics[width=0.327\linewidth]
	{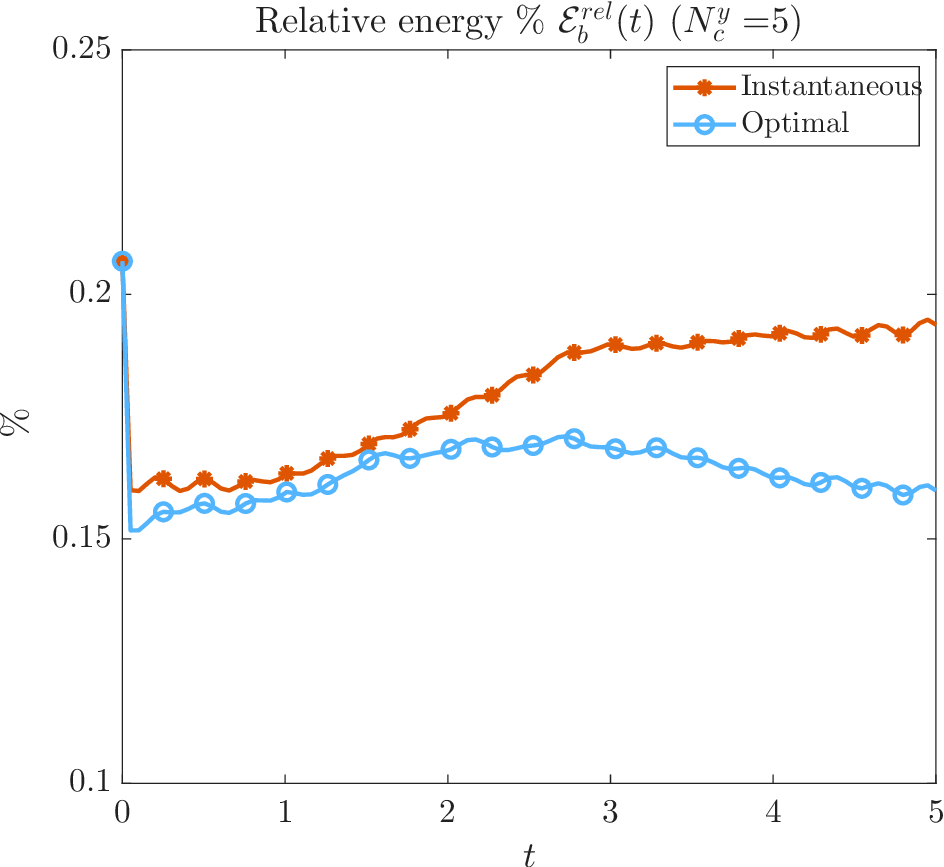}
	\includegraphics[width=0.327\linewidth]
	{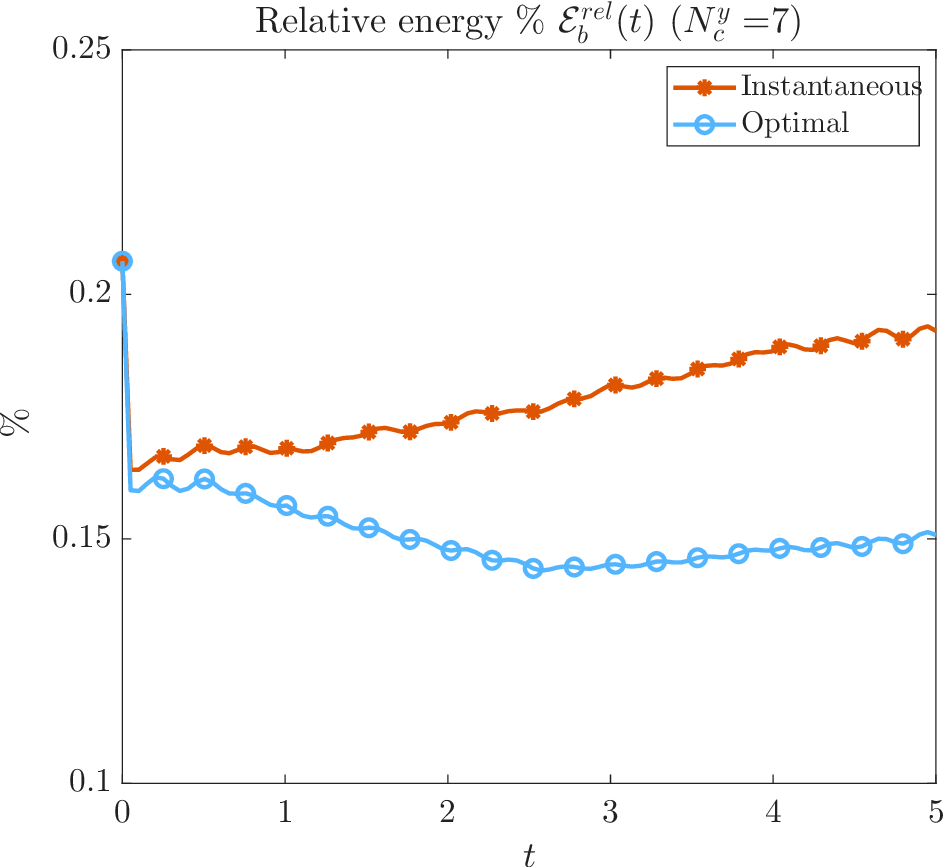}
	\vspace*{5pt}
	\caption{Kelvin--Helmholtz instability. Comparison of
		$\mathcal E_b^{\mathrm{rel}}(t)$ obtained with the instantaneous
		and optimal controls. From left to right:
		$N_c^y=3$, $N_c^y=5$, and $N_c^y=7$.}
	\label{fig:kelvin_comparison_Binst_Bopt}
\end{figure}

For the controlled simulations, we set $N_c^x=1$ and consider
$N_c^y\in\{3,5,7\}$. Consequently, the control cells extend over the
entire $x$-direction and the reconstructed magnetic field depends
only on $y$ and $t$. Figure~\ref{fig:functional_kelvin} reports the convergence history of
the optimization procedure. The panels show the reduced cost
functional $\widehat{\mathcal J}(\mathbf b^{(m)})$ and the discrete
gradient norm $\|\mathbf g^{(m)}\|_2$.
 Figure~\ref{fig:kelvin_dynamics_control_Nc5} shows
three density snapshots of the controlled dynamics at $t=2$,
$t=3$, and $t=6$ for $N_c^y=5$. Similar qualitative results are
obtained for $N_c^y=3$ and $N_c^y=7$. The control attenuates the
development of the instability and substantially reduces transport
toward the upper and lower boundaries. The lower panel show the
optimal controls coefficients for
$N_c^y\in\{3,5,7\}$. Finally, Figure~\ref{fig:kelvin_comparison_Binst_Bopt} provides a
direct comparison in terms of relative thermal energy between the instantaneous and optimal controls.
As for the previous cases, in most of the cases, the optimized
strategy yields a further reduction in the fraction of thermal energy
localized near the boundaries.


\section{Conclusions}\label{sec:conclusions}

We have developed a full-horizon, open-loop strategy for controlling a
2D--2D Vlasov--Poisson system by means of an externally imposed magnetic
field. The field is parameterized by time-dependent coefficients
associated with smooth, spatially localized kernels, thereby reducing the
infinite-dimensional control problem to the determination of a finite family
of control functions. The objective combines the tracking and dispersion of
selected position and velocity observables with a penalty on the control
effort, while a pointwise bound limits the intensity of the reconstructed
magnetic field.

Starting from the kinetic optimization problem, we derived the continuous
first-order optimality conditions.  We then projected this system onto a particle--grid
representation and obtained particle adjoint equations. The
discretize-then-optimize calculation in
Appendix~\ref{app:computations_optimality} shows that the same semidiscrete
system is recovered when the particle deposition and interpolation operators,
the grid quadrature, and the adjoint discrete operators are paired
consistently, up to the localization error associated with the particle shape
function. The resulting optimization problem is solved by a projected gradient
method with Armijo backtracking on a reduced particle ensemble. An instantaneous
feedback control, adapted to the same smooth spatial parameterization, provides
the initial guess for the full-horizon iteration.

The numerical experiments illustrate the behavior of the method. In the tests considered, the instantaneous control already gives an effective initialization and is often close to a
stationary point of the discrete full-horizon problem. Nevertheless, the
full-horizon optimization produces a further correction that decreases the
reduced objective and generally lowers the normalized thermal-energy
diagnostic in the boundary strips. 

Several extensions are natural such as the use of
collisional or electromagnetic kinetic models, more realistic geometries and
coil-generated fields. Moreover, objectives based directly on particle or energy
fluxes at material boundaries would bring the framework closer to magnetic
confinement applications.

\appendix

\section{Discretize-then-optimize derivation of the particle--grid optimality system}
\label{app:computations_optimality}

We provide a discretize-then-optimize derivation of the semidiscrete
particle--grid optimality system introduced in
Section~\ref{sec:derivation_control}. Time remains continuous throughout this
appendix. Let
\[
S_{\ell i}(t)
:=
\mathcal S_h\bigl(\xx_\ell-\xx_i(t)\bigr),
\qquad
\ell=1,\ldots,N_h,\qquad i=1,\ldots,N,
\]
the density and the electric field interpolated to particle \(i\)
are
\begin{equation}\label{eq:app_particle_grid_state}
	\rho_\ell^N
	=
	\sum_{i=1}^N\omega_iS_{\ell i},
	\qquad
	-(A_h\boldsymbol\phi)_\ell
	=
	\rho_\ell^N-\bar{\rho},
	\qquad
	\EE_i^h
	=
	-\sum_{\ell=1}^{N_h}\eta_\ell
	(G_h\boldsymbol\phi)_\ell S_{\ell i}.
\end{equation}
The particle dynamics is therefore
\begin{equation}\label{eq:app_particle_dynamics}
	\dot{\xx}_i=\vv_i,
	\qquad
	\dot{\vv}_i=F_i^h
	:=
	\EE_i^h+\vv_i\times B_i,
	\qquad
	B_i
	=
	\sum_{k=1}^{N_g}b_k(t)\mathcal K_k(\xx_i).
\end{equation}
Set
\[
\mathcal M^N:=\sum_{i=1}^N\omega_i,
\qquad
m_q^N(t)
:=
\frac{1}{\mathcal M^N}
\sum_{i=1}^N\omega_i
\psi_q(\xx_i(t),\vv_i(t)).
\]
For \(q\in\mathcal I=\{\mathrm{x},\mathrm{v}\}\), define
\begin{equation}\label{eq:app_particle_tracking}
	\begin{split}
		\mathcal D_q^N(t)
		:={}&
		\frac{\alpha_q}{2}
		\sum_{i=1}^N\omega_i
		\left|
		\psi_q(\xx_i,\vv_i)-\widehat\psi_q
		\right|^2+
		\frac{\beta_q}{2}
		\sum_{i=1}^N\omega_i
		\left|
		\psi_q(\xx_i,\vv_i)-m_q^N
		\right|^2,
	\end{split}
\end{equation}
and
\begin{equation}\label{eq:app_particle_regularization}
	\mathcal R_{\mathrm{diag}}^N(t)
	:=
	\sum_{k=1}^{N_g}|b_k(t)|^2
	\sum_{i=1}^N\omega_i\mathcal K_k(\xx_i(t)).
\end{equation}
The cost functional at the particle level becomes
\begin{equation}\label{eq:app_particle_objective}
	\mathcal J_N(\mathbf X,\mathbf V,\mathbf b)
	:=
	\int_0^T
	\left[
	\sum_{q\in\mathcal I}\mathcal D_q^N(t)
	+\frac{\gamma}{2}\mathcal R_{\mathrm{diag}}^N(t)
	\right]dt
	+
	\sum_{q\in\mathcal I}\mathcal D_q^N(T),
\end{equation}
with
\(
\mathbf X=(\xx_1,\ldots,\xx_N)
\)
and
\(
\mathbf V=(\vv_1,\ldots,\vv_N)
\).
Let \(\mathcal G_i^N\) denote the particle contribution
\begin{equation}\label{eq:app_particle_source}
	\begin{split}
		\mathcal G_i^N(t)
		:={}&
		\sum_{q\in\mathcal I}
		\left[
		\frac{\alpha_q}{2}
		\left|\psi_q(\xx_i,\vv_i)-\widehat\psi_q\right|^2
		+
		\frac{\beta_q}{2}
		\left|\psi_q(\xx_i,\vv_i)-m_q^N\right|^2
		\right]
		\\
		&+
		\frac{\gamma}{2}
		\sum_{k=1}^{N_g}|b_k(t)|^2\mathcal K_k(\xx_i).
	\end{split}
\end{equation}
Thus
\[
\nabla_{\xx_i}\mathcal G_i^N
:=
\frac{1}{\omega_i}
\nabla_{\xx_i}
\left(
\sum_{j=1}^N\omega_j\mathcal G_j^N
\right),
\qquad
\nabla_{\vv_i}\mathcal G_i^N
:=
\frac{1}{\omega_i}
\nabla_{\vv_i}
\left(
\sum_{j=1}^N\omega_j\mathcal G_j^N
\right),
\]
and the choices
\[
\psi_{\mathrm{x}}=y,
\qquad
\psi_{\mathrm{v}}=v_y,
\qquad
\widehat\psi_{\mathrm{x}}=\widehat y,
\qquad
\widehat\psi_{\mathrm{v}}=\widehat v_y
\]
give precisely the expressions in \eqref{eq:nabla_G}.
The particle--grid Lagrangian is
\begin{equation}\label{eq:app_particle_lagrangian}
	\begin{split}
		\mathscr L_N
		:={}&
		\mathcal J_N(\mathbf X,\mathbf V,\mathbf b)
		\\
		&+
		\sum_{i=1}^N\omega_i
		\int_0^T
		\left[
		\pp_i\cdot(\dot{\xx}_i-\vv_i)
		+
		\qq_i\cdot
		\bigl(\dot{\vv}_i-\EE_i^h-\vv_i\times B_i\bigr)
		\right]dt
		\\
		&+
		\int_0^T
		\left\langle
		\widetilde{\boldsymbol\lambda},
		-A_h\boldsymbol\phi-\boldsymbol\rho^N+\bar{\boldsymbol \rho}_h
		\right\rangle_hdt.
	\end{split}
\end{equation}
Here \(\bar{\boldsymbol \rho}_h\) denotes the grid vector whose entries are all equal
to $\bar{\rho}$ and 
$
\langle \boldsymbol a,\boldsymbol c\rangle_h
:=
\sum_{\ell=1}^{N_h}\eta_\ell
\boldsymbol a_\ell\cdot\boldsymbol c_\ell.
$
Stationarity with respect to \(\pp_i\), \(\qq_i\), and
\(\widetilde{\boldsymbol\lambda}\) recovers
\eqref{eq:app_particle_grid_state}--\eqref{eq:app_particle_dynamics}.
To compute the variation with respect to the grid potential, define
\begin{equation}\label{eq:app_Q_lambda}
	\boldsymbol Q_\ell
	:=
	\sum_{i=1}^N\omega_iS_{\ell i}\qq_i,
	\qquad
	\lambda_i^h
	:=
	\sum_{\ell=1}^{N_h}\eta_\ell
	\widetilde\lambda_\ell S_{\ell i}.
\end{equation}
Using the grid inner product and the adjoint operators \(A_h^*\) and
\(G_h^*\), stationarity with respect to \(\boldsymbol\phi\) gives
\begin{equation}\label{eq:app_adjoint_poisson}
	-A_h^*\widetilde{\boldsymbol\lambda}
	+G_h^*\boldsymbol Q
	=0,
\end{equation}
which coincides with \eqref{eq:backward_discrete_2}.

Integrating the terms containing \(\dot{\xx}_i\) and
\(\dot{\vv}_i\) by parts in time, and using that the initial particle
states are fixed, stationarity with respect to \(\xx_i\) and \(\vv_i\)
yields
\begin{equation}\label{eq:app_particle_adjoints}
	\begin{split}
		-\dot{\pp}_i
		&=
		(\nabla_{\xx_i}F_i^h)^T\qq_i
		-\nabla_{\xx_i}\mathcal G_i^N
		+\nabla_{\xx_i}\lambda_i^h,
		\\
		-\dot{\qq}_i
		&=
		\pp_i
		-\nabla_{\vv_i}\mathcal G_i^N
		-\qq_i\times B_i,
	\end{split}
	\qquad i=1,\ldots,N.
\end{equation}
Moreover,
\[
\nabla_{\xx_i}\lambda_i^h
=
\sum_{\ell=1}^{N_h}\eta_\ell
\widetilde\lambda_\ell
\nabla_{\xx_i}\mathcal S_h(\xx_\ell-\xx_i),
\]
so that \eqref{eq:app_particle_adjoints} is exactly
\eqref{eq:backward_discrete}. The terminal conditions are
\begin{equation}\label{eq:app_particle_terminal_conditions}
	\pp_i(T)
	=
	-\nabla_{\xx_i}\mathcal G_{T,i}^N,
	\qquad
	\qq_i(T)
	=
	-\nabla_{\vv_i}\mathcal G_{T,i}^N,
\end{equation}
where \(\mathcal G_{T,i}^N\) contains only the terminal tracking and
dispersion terms, and no control regularization.
Finally, the derivative of \eqref{eq:app_particle_lagrangian} with respect to
\(b_k\) is
\begin{equation}\label{eq:app_particle_gradient}
	g_k^N(t)
	=
	\sum_{i=1}^N\omega_i
	\left[
	\gamma b_k(t)+\vv_i(t)\times\qq_i(t)
	\right]
	\mathcal K_k(\xx_i(t)),
	\qquad k=1,\ldots,N_g.
\end{equation}
At the particle--grid level, we enforce the field bound at the spatial
grid nodes and set
\begin{equation}\label{eq:Uadm_h}
   \mathcal U_{\mathrm{ad}}^h
:=
\left\{
\mathbf b\in L^2(0,T;\RR^{N_g}):
\left|
\sum_{k=1}^{N_g}b_k(t)\mathcal K_k(\xx_\ell)
\right|
\leq B_{\rm max}
\text{ for a.e. }t,\quad
\ell=1,\ldots,N_h
\right\}.
\end{equation}
For a bounded control, stationarity is then expressed by the variational
inequality
\begin{equation}\label{eq:app_particle_variational_inequality}
	\sum_{k=1}^{N_g}
	\int_0^T
	(g_k^N)^\star(t)
	\bigl(c_k(t)-b_k^\star(t)\bigr)dt
	\geq0
	\qquad
	\text{for every }\mathbf c\in\mathcal U_{\mathrm{ad}}^h.
\end{equation}
Thus, the discretize-then-optimize calculation recovers
\eqref{eq:backward_discrete}, \eqref{eq:backward_discrete_2},
\eqref{eq:particle_terminal_conditions}, and
\eqref{eq:gradient_discrete}. For the particle objective
\eqref{eq:app_particle_objective}, the calculation is exact at the
particle--grid semidiscrete level when the same discretization are employed. Its comparison with the
projected continuous system additionally involves replacing convolutions
against \(\mathcal S_h\) by evaluations along the particle trajectories. For a shape function of order
\(r\), these localization discrepancies are of order \(\mathcal O(h^r)\);
in particular, this accounts for the remainder displayed in
\eqref{eq:gradient_discrete}.

\subsection{Relation with the fully discrete adjoint}

The derivation above is semidiscrete, with time kept continuous.
Indeed, in Section~\ref{sec:numerical_optimization}, the gradient
\eqref{eq:gradient_discrete_time} is obtained by discretizing in time the
semidiscrete optimality system. We now compare it with the gradient obtained
by differentiating the fully discrete state problem.
The semi-implicit forward scheme \eqref{eq:forward_discretization} can then
be written as
\begin{equation}\label{eq:fully_discrete_forward_map}
\begin{cases}
\displaystyle
\vv_i^{n+1}
=
\mathcal V_i^n
:=
(A_i^n)^{-1}
\left(
\vv_i^n+\Delta t\,\EE_i^{h,n}
\right),
\\[1mm]
\displaystyle
\xx_i^{n+1}
=
\mathcal X_i^n
:=
\xx_i^n+\Delta t\,\mathcal V_i^n,
\end{cases}
\end{equation}
where 
\[
\mathbb J
=
\begin{pmatrix}
0&1\\
-1&0
\end{pmatrix},
\qquad
A_i^n
:=
I-\Delta t\,B_i^n\mathbb J.
\]
Together with the discrete Poisson problem, this defines the fully
discrete state map
\[
(\mathbf X^n,\mathbf V^n,\mathbf b^n)
\longmapsto
(\mathbf X^{n+1},\mathbf V^{n+1}).
\]
Using the same time quadrature as in Section~\ref{sec:numerical_optimization},
the corresponding discrete Lagrangian may be written as
\begin{equation}\label{eq:fully_discrete_lagrangian}
\begin{split}
\mathscr L_{N,h,\Delta t}
:={}&
\mathcal J_{N,h,\Delta t}
\\
&+
\sum_{n=0}^{N_t-1}
\sum_{i=1}^{N_s}
\widehat\omega_i
\Big[
\pp_{i,\mathrm{DTO}}^{n+1}
\cdot
\bigl(\xx_i^{n+1}-\mathcal X_i^n\bigr)+
\qq_{i,\mathrm{DTO}}^{n+1}
\cdot
\bigl(\vv_i^{n+1}-\mathcal V_i^n\bigr)
\Big].
\end{split}
\end{equation}
Since
\[
\mathcal X_i^n=\xx_i^n+\Delta t\,\mathcal V_i^n,
\]
we introduce
\begin{equation}\label{eq:zeta_DTO}
\boldsymbol\zeta_i^n
:=
(A_i^n)^{-T}
\left(
\qq_{i,\mathrm{DTO}}^{n+1}
+
\Delta t\,\pp_{i,\mathrm{DTO}}^{n+1}
\right),
\end{equation}
which collects the contribution of the position and velocity adjoints
through the semi-implicit velocity update. Thus, the discrete adjoint writes as
\begin{equation}\label{eq:fully_discrete_adjoint_DTO}
\begin{cases}
\displaystyle
\pp_{i,\mathrm{DTO}}^n
=
\pp_{i,\mathrm{DTO}}^{n+1}
+
\Delta t
\left[
\left(
\nabla_{\xx_i}F_i^{h,n,+}
\right)^T
\boldsymbol\zeta_i^n
-
\nabla_{\xx_i}\mathcal G_i^{N,n}
+
L_i^n
\right],
\\[1em]
\displaystyle
\qq_{i,\mathrm{DTO}}^n
=
\boldsymbol\zeta_i^n
-
\Delta t\,
\nabla_{\vv_i}\mathcal G_i^{N,n},
\end{cases}
\end{equation}
with terminal conditions
\[
\pp_{i,\mathrm{DTO}}^{N_t}
=
-\nabla_{\xx_i}\mathcal G_{T,i}^N,
\qquad
\qq_{i,\mathrm{DTO}}^{N_t}
=
-\nabla_{\vv_i}\mathcal G_{T,i}^N,
\]
and where
\[
F_i^{h,n,+}
:=
\EE_i^{h,n}
+
\vv_i^{n+1}\times B_i^n.
\]
Here $L_i^n$ denotes the contribution associated with the dependence of the
electric field on the particle density, obtained from the same discrete
Poisson-adjoint relation derived in \eqref{eq:app_particle_adjoints} and evaluated at time $t^n$.
Finally, stationarity with respect to $b_k^n$ gives, for the diagonal
regularization \eqref{eq:diagonal_control_penalty},
\begin{equation}\label{eq:gradient_fully_discrete_DTO}
g_{k,n}^{\mathrm{DTO}}
=
\sum_{i=1}^{N_s}
\widehat\omega_i
\left[
\gamma b_k^n
+
\vv_i^{n+1}\times\boldsymbol\zeta_i^n
\right]
\mathcal K_k(\xx_i^n).
\end{equation}
We observe that since 
\[
(A_i^n)^{-T}
=
I-\Delta t\,B_i^n\mathbb J+\mathcal O(\Delta t^2),
\qquad
\vv_i^{n+1}=\vv_i^n+\mathcal O(\Delta t),
\]
and the adjoint discretizations in \eqref{eq:zeta_DTO} and \eqref{eq:gradient_fully_discrete_DTO} are such that
\[
\boldsymbol\zeta_i^n
=
\qq_i^n+\mathcal O(\Delta t),\qquad
\vv_i^{n+1}\times\boldsymbol\zeta_i^n
=
\vv_i^n\times\qq_i^n+\mathcal O(\Delta t).
\]
It follows that the gradient used in
Section~\ref{sec:numerical_optimization} is a first-order consistent
approximation of the fully discrete gradient, indeed we have
\begin{equation}\label{eq:gradient_DTO_OTD_consistency}
g_{k,n}^{\mathrm{DTO}}
=
g_{k,n}
+
\mathcal O(\Delta t),
\end{equation}
where $g_{k,n}$ is the gradient defined in
\eqref{eq:gradient_discrete_time}. 

 \section{Computation of the instantaneous control} \label{app:instantaneous_control} 
In the following, we will briefly describe how to compute the baseline instantaneous feedback control $B_{\mathbf{b}}^0(t,\xx)$, referring to \cite{albi2025instantaneous}-\cite{albi2025robust} for a complete description of the methodology. Accordingly, the initial guess for the open-loop control $B_{\mathbf{b}}(t,\xx)$ in the optimization procedure is chosen to be this instantaneous feedback law. 

We consider the following minimization problem at the particle level: 
\begin{equation}\label{eq:J_inst}
	\min_{\mathbf{b}^0 \in \mathcal{\mathcal{U}}_{ad}} \mathcal{J}(\mathbf{b}^0) = \int_{t^n}^{t^{n+1}} \left[ \mathcal{G}_i^N(t) +  \frac{\gamma}{2} \sum_{k=1}^{N_g}\vert b_k^0(t) \vert^2 \mathcal{K}_k(\xx_i) \right]   ~ dt,   \qquad \text{s.t. } \quad  \eqref{eq:char_curves},
\end{equation}
with $\mathcal{G}_i^N (t)$ and $\mathcal{U}_{ad}$ defined as in \eqref{eq:app_particle_source} and \eqref{eq:Uadm_h} respectively, and $\gamma>0$ a penalization coefficient.

We introduce an explicit first order discretization of the characteristic curves in \eqref{eq:char_curves} over the time interval $[0,T]$, with step size $\Delta t>0$ as follows:
\begin{equation}\label{eq:char_curves_discretized}
	\begin{cases}
	\xx_i^{n+1} = \xx_i^n + \Delta t \vv_i^{n+1},\\
	\vv_i^{n+1} = \vv_i^n + \Delta t (\EE_i^n + \vv_i^n \times B_i^{0,n}), 
	\end{cases}
\end{equation}
for any $i=1,\ldots,N$, and where 
\begin{equation}\label{eq:Bi_inst}
	B_i^{0,n} = \sum_{k=1}^{N_g}  b_{k}^{0,n}  \mathcal{K}_k(\xx_i),
\end{equation}
with  $\mathcal{K}$ denoting the spatial basis functions defined as in \eqref{eq:Phi_control}. 
Following \cite{albi2025instantaneous,albi2025robust}, we introduce the
scaling
\begin{equation}\label{eq:instantaneous_scaling}
    \alpha_{\mathrm{x}}^{\Delta t}
    =
    \frac{\alpha_{\mathrm{x}}}{\Delta t},
    \qquad
    \beta_{\mathrm{x}}^{\Delta t}
    =
    \frac{\beta_{\mathrm{x}}}{\Delta t},
    \qquad
    \gamma_{\Delta t}
    =
    \gamma\Delta t,
\end{equation}
while $\alpha_{\mathrm{v}}$ and $\beta_{\mathrm{v}}$ are kept fixed.
Since the $y$-component of the velocity update satisfies
\[
    v_{y,i}^{n+1}
    =
    v_{y,i}^n
    +
    \Delta t\,E_{y,i}^n
    -
    \Delta t\,v_{x,i}^n B_i^{0,n},
\]
we have
\[
    \frac{\partial v_{y,i}^{n+1}}
         {\partial b_j^{0,n}}
    =
    -\Delta t\,v_{x,i}^n
    \mathcal K_j(\xx_i^n),
    \qquad
    \frac{\partial y_i^{n+1}}
         {\partial b_j^{0,n}}
    =
    -\Delta t^2\,v_{x,i}^n
    \mathcal K_j(\xx_i^n).
\]
Using \eqref{eq:instantaneous_scaling} and considering the limit
$\Delta t\to0$ with $t^n=n\Delta t\to t$, the first-order optimality
conditions yield the instantaneous feedback law
\begin{equation}\label{eq:instantaneous_control}
    H(t)\,\mathbf b^0(t)=F(t),
\end{equation}
where $H(t)$ is a diagonal matrix based on the assumption \eqref{eq:assumption-kernel}, and we have 
\begin{equation}\label{eq:instantaneous_MF}
    H_{jj}(t)
    =
    \gamma
    \sum_{i=1}^{N}
    \omega_i
    \mathcal K_j(\xx_i(t)),
    \qquad
    F_j(t)
    =
    \sum_{i=1}^{N}
    \omega_i
    v_{x,i}(t)
    \mathcal R_i^0(t)
    \mathcal K_j(\xx_i(t)),
\end{equation}
with
\begin{equation}\label{eq:instantaneous_R}
    \mathcal R_i^0(t)
    =
    \alpha_{\mathrm{x}}
    \left(y_i(t)-\widehat y\right)
    +
    \beta_{\mathrm{x}}
    \left(y_i(t)-\overline y^{\,N}(t)\right)
    +
    \alpha_{\mathrm{v}}
    \left(v_{y,i}(t)-\widehat v_y\right)
    +
    \beta_{\mathrm{v}}
    \left(v_{y,i}(t)-\overline v_y^{\,N}(t)\right).
\end{equation}
Equation \eqref{eq:instantaneous_control} provides the unconstrained
instantaneous control, if the constraint in $\mathcal U_{\mathrm{ad}}^h$ is
active, the corresponding constrained quadratic problem is solved instead.
Once determined through this first-order framework, the instantaneous control coefficients $\mathbf b^0(t)$ are applied within the second-order scheme to advance the particle trajectories.
\section*{Acknowledgement}
This work has been written within the activities of GNCS and GNFM groups of INdAM (Italian National Institute of High Mathematics). The research of FF and LP was supported by the Italian Ministry of University and Research (MUR) through the Fondo Italiano per la Scienza 2022–2023 (Bando FIS 2), under Decree No. 1236 of 01/08/2023, within the FIS2023-01334 Advanced Grant “Tackling Complexity: Advanced Numerical Approaches for Multiscale Systems with Uncertainties” (ADAMUS).
The work of LP was supported by the Royal Society
Wolfson Fellowship “Uncertainty quantification, data-driven simulations and learning of multiscale complex systems governed by PDEs”.
  GA and GD acknowledge the support of the PRIN 2026 - 2026PYFEZ9, ``Emergent Structures in Kinetic Theory: from Gases to Plasmas" funded with the contribution of MUR.

\bibliographystyle{abbrv}
\bibliography{biblio.bib}
\end{document}